\RequirePackage{ifpdf}
\ifpdf 
\documentclass[pdftex]{sigma}
\else
\documentclass{sigma}
\fi

\usepackage{amsbsy}
\usepackage{amscd}

\usepackage{mathrsfs}
\usepackage{nicefrac}

\usepackage{youngtab}
\usepackage[all]{xy}
\usepackage{mathdots}

\newenvironment{aenume}{%
  \begin{enumerate}%
  }{\end{enumerate}}

\usepackage{verbatim}
\usepackage{version}
\usepackage{color}

\excludeversion{NB}

\newtheorem{Theorem}[equation]{Theorem}
\newtheorem{Corollary}[equation]{Corollary}
\newtheorem{Lemma}[equation]{Lemma}
\newtheorem{Proposition}[equation]{Proposition}

\theoremstyle{definition}
\newtheorem{Definition}[equation]{Definition}
\newtheorem{Example}[equation]{Example}

\newtheorem{Remark}[equation]{Remark}
\newtheorem{Remarks}[equation]{Remarks}

\numberwithin{equation}{section}

\newcommand{\thmref}[1]{Theorem~\ref{#1}}
\newcommand{\secref}[1]{Section~\ref{#1}}
\newcommand{\lemref}[1]{Lemma~\ref{#1}}
\newcommand{\propref}[1]{Proposition~\ref{#1}}

\newcommand{\subsecref}[1]{Section~\ref{#1}}

\newcommand{\remref}[1]{Remark~\ref{#1}}
\newcommand{\remsref}[1]{Remarks~\ref{#1}}
\newcommand{\lsp}[2]{{\mskip-.3mu}{}^{#1}\mskip-1mu{#2}}
\newcommand{\defeq}{\overset{\operatorname{\scriptstyle def}}{=}}
\newcommand{\C}{{\mathbb C}}
\newcommand{\Z}{{\mathbb Z}}
\newcommand{\Q}{{\mathbb Q}}
\newcommand{\R}{{\mathbb R}}
\newcommand{\proj}{{\mathbb P}}
\newcommand{\CP}{\proj}

\newcommand{\SL}{\operatorname{\rm SL}}

\newcommand{\GL}{\operatorname{GL}}

\newcommand{\algsl}{\operatorname{\mathfrak{sl}}} 

\newcommand{\gl}{\operatorname{\mathfrak{gl}}}

\newcommand{\g}{{\mathfrak g}}
\newcommand{\h}{{\mathfrak h}}

\newcommand{\End}{\operatorname{End}}
\newcommand{\Hom}{\operatorname{Hom}}

\newcommand{\Ker}{\operatorname{Ker}}
\newcommand{\Coker}{\operatorname{Coker}}
\newcommand{\Ima}{\operatorname{Im}}

\newcommand{\tr}{\operatorname{tr}}

\newcommand{\id}{\operatorname{id}}

\newcommand{\vin}[1]{\operatorname{i}(#1)} 
\newcommand{\vout}[1]{\operatorname{o}(#1)} 
\newcommand{\bM}{{\mathbf M}} 
\newcommand{\M}{{\mathfrak M}} 
\newcommand{\Mreg}{\M^{\operatorname{s}}}
\newcommand{\La}{{\mathfrak L}} 
\newcommand{\dslash}{/\!\!/} 
\newcommand{\bv}{{\mathbf v}} 
\newcommand{\bw}{{\mathbf w}} 
\newcommand{\bC}{{\mathbf C}} 
\newcommand{\bA}{{\mathbf A}} 
\newcommand{\bI}{{\mathbf I}} 

\newcommand{\be}{{\mathbf e}}

\newcommand{\bc}{\mathbf c}
\newcommand{\codim}{\operatorname{codim}} 
\newcommand{\gr}{\operatorname{gr}} 
\newcommand{\topdeg}{{\operatorname{top}}} 
\newcommand{\Uq}{{\mathbf U}_q(\mathfrak g)} 
\newcommand{\Un}{{\mathbf U}_q}
\newcommand{\Ul}{{\mathbf U}_q({\mathbf L}{\mathfrak g})} 

\newcommand{\bU}{\mathbf U} 
\newcommand{\Pa}{{\mathfrak P}} 
\newcommand{\HomE}{\operatorname{E}}
\newcommand{\HomL}{\operatorname{L}}

\newcommand{\shfO}{\mathcal O}

\newcommand{\Wedge}{{\textstyle \bigwedge}}

\newcommand{\linf}{\ell_\infty}
\newcommand{\Irr}{\operatorname{Irr}}
\newcommand{\Res}{\operatorname{Res}}
\newcommand{\Lareg}{\mathfrak L^{\operatorname{s}}}
\newcommand{\MR}[1]{}
\newcommand{\aff}{\mathrm{af\/f}}
\newcommand{\affLevi}{(\g_{I_0^0})_{\aff}}

\newcommand{\hf}{\hfill}

\makeatletter
\newdimen\y@inside
\y@inside=\y@boxdim
\advance\y@inside by -3\y@linethick
\definecolor{halfgray}{gray}{0.5}
\newcommand{\rf}{\color{halfgray}
\rule[-.2\y@inside]{\y@inside}{\y@inside}
}
\makeatother
\newcommand{\graysquare}{\color{halfgray}\blacksquare}

\begin{document}
\allowdisplaybreaks

\renewcommand{\thefootnote}{$\star$}

\renewcommand{\PaperNumber}{003}

\FirstPageHeading

\ShortArticleName{Quiver Varieties and Branching}

\ArticleName{Quiver Varieties and Branching\footnote{This paper is a
contribution to the Special Issue on Kac--Moody Algebras and Applications. The
full collection is available at
\href{http://www.emis.de/journals/SIGMA/Kac-Moody_algebras.html}{http://www.emis.de/journals/SIGMA/Kac-Moody{\_}algebras.html}}}

\Author{Hiraku NAKAJIMA~$^{\dag\ddag}$}

\AuthorNameForHeading{H.~Nakajima}

\Address{$^{\dag}$~Department of Mathematics, Kyoto University, Kyoto 606-8502,
Japan}
\Address{$^{\ddag}$~Research Institute for Mathematical Sciences, Kyoto University,
Kyoto 606-8502, Japan}

\EmailD{\href{mailto:nakajima@math.kyoto-u.ac.jp}{nakajima@math.kyoto-u.ac.jp}}

\ArticleDates{Received September 15, 2008, in f\/inal form January 05,
2009; Published online January 11, 2009}

\Abstract{Braverman and Finkelberg  recently proposed the
  geometric Satake correspondence for the af\/f\/ine Kac--Moody group
  $G_\aff$ [Braverman A., Finkelberg M.,  \href{http://arxiv.org/abs/0711.2083}{arXiv:0711.2083}]. They conjecture that intersection cohomology sheaves on
  the Uhlenbeck compactif\/ication of the framed moduli space of
  $G_{\mathrm{cpt}}$-instantons on $\R^4/\Z_r$ correspond to weight
  spaces of representations of the Langlands dual group $G_\aff^\vee$
  at level $r$. When $G = \SL(l)$, the Uhlenbeck compactif\/ication is
  the quiver variety of type $\algsl(r)_\aff$, and their conjecture
  follows from the author's earlier result and I.~Frenkel's level-rank
  duality. They further introduce a~convolution diagram which
  conjecturally gives the tensor product multiplicity [Braverman~A., Finkelberg~M., Private communication, 2008].
  In this paper, we develop the theory for the
  branching in quiver varieties and check this conjecture for $G=\SL(l)$.}


\Classification{17B65; 14D21}

\tableofcontents

\renewcommand{\thefootnote}{\arabic{footnote}}
\setcounter{footnote}{0}

\section{Introduction}

In \cite{Na-quiver,Na-alg} the author showed that the top degree
homology group of a Lagrangian subvariety~$\La$ in a quiver variety
$\M$ has a structure of an integrable highest weight representation of
a Kac--Moody Lie algebra $\g$.
In a subsequent work \cite{Na-qaff} the author showed that the
equivariant $K$-homology group of $\La$ has a structure of an
$\ell$-integrable highest weight representation of the quantum loop
algebra $\Ul$ (e.g., the quantum af\/f\/ine algebra if $\g$ is of f\/inite
type, the quantum toroidal algebra if $\g$ is of af\/f\/ine type). As an
application, the characters of arbitrary irreducible representations
of $\Ul$ were computed in terms of the intersection cohomology (IC for
short) groups associated with graded/cyclic quiver varieties (= the
f\/ixed point set in the quiver variety with respect to $\C^*$/cyclic
group action), and hence analogs of Kazhdan--Lusztig polynomials. This
result cannot be proved by a purely algebraic method. See
\cite{MR1989196} for a survey.

The quiver variety $\M$ is def\/ined as a geometric invariant theory
quotient of an af\/f\/ine varie\-ty~$\mu^{-1}(0)$ by a product of general
linear groups with respect to a particular choice of a stability
condition $\zeta$, a lift of an action to the trivial line bundle over
$\mu^{-1}(0)$. Let us denote~$\M$ by~$\M_\zeta$ hereafter to emphasize
a choice of a stability condition.
We can consider other stability conditions. For example, if we choose
the trivial lift $\zeta = 0$, then we get an af\/f\/ine algebraic variety
$\M_0$. It has been studied already in the literature, and played
important roles. For example, we have a projective morphism
$\pi_{0,\zeta}\colon\M_\zeta\to \M_0$, and $\La$ is the inverse image
$\pi_{0,\zeta}^{-1}(0)$ of a distinguished point $0\in\M_0$. Moreover
$\M_0$ has a natural stratif\/ication parametrized by conjugacy classes
of stabilizers, and their IC complexes give the restriction
multiplicities of the above $\Ul$-module to $\Uq$.

In this paper we study a more general stability condition
$\zeta^\bullet$ whose corresponding variety~$\M_{\zeta^\bullet}$ sits
between $\M_\zeta$ and $\M_0$: the morphism $\pi$ factorizes
\(
   \M_\zeta\xrightarrow{\pi_{\zeta^\bullet, \zeta}}
    \M_{\zeta^\bullet}\xrightarrow{\pi_{0,\zeta^\bullet}}\M_0.
\)
Under a mild assumption, $\M_{\zeta^\bullet}$ is a partial resolution
of singularities of $\M_0$, while $\M_\zeta$ is a full resolution.
Then we show that the top degree cohomology group of
$\pi^{-1}_{0,\zeta^\bullet}(0)$ with coef\/f\/icients in the IC complex of
a stratum of $\M_{\zeta^\bullet}$ gives the restriction multiplicities
of an integrable highest weight representation to a subalgebra
determined by $\zeta^\bullet$.
This follows from a general theory for representations constructed by
the convolution product \cite{BM}, but we identify the subalgebras
with the Levi subalgebra of $\g$ or the af\/f\/ine Lie algebra of the Levi
subalgebra for certain choices of $\zeta^\bullet$ (see
Theorems~\ref{thm:branch}, \ref{thm:branch'}). (For a quiver variety of
f\/inite or af\/f\/ine type, they exhaust all choices up to the Weyl group
action.)

This work is motivated by a recent proposal by Braverman and
Finkelberg \cite{braverman-2007,braverman-prep} on a conjectural
af\/f\/ine Kac--Moody group analog of the geometric Satake correspondence.
(See also a recent paper \cite{bk} for the af\/f\/ine analog of the
classical Satake correspondence.)
The ordinary geometric Satake correspondence says that the category of
equivariant perverse sheaves on the af\/f\/ine Grassmannian $G(\mathcal
K)/G(\mathcal O)$ associated with a f\/inite dimensional complex simple
Lie group $G$ is equivalent to the category of f\/inite dimensional
representations of the Langlands dual group $G^\vee$, so that IC
sheaves of $G(\mathcal O)$-orbits correspond to irreducible
representations.
Here $\mathcal K = \C((s))$, $\mathcal O = \C[[s]]$.
If we would have the af\/f\/ine Grassmannian for the af\/f\/ine Kac--Moody
group $G_\aff$ (hence called the double af\/f\/ine Grassmannian), then it
should correspond to representations of the af\/f\/ine Kac--Moody group
$G_{\aff}^\vee$.
It is not clear whether we can consider IC sheaves on the double
af\/f\/ine Grassmannian $G_{\aff}(\mathcal K)/G_{\aff}(\mathcal O)$ or
not, but Braverman and Finkelberg \cite{braverman-2007,braverman-prep}
propose that the transversal slice of a $G_{\aff}(\mathcal O)$-orbit
in the closure of a bigger $G_{\aff}(\mathcal O)$-orbit should be the
Uhlenbeck partial compactif\/ication of the framed moduli space of
$G_{\mathrm{cpt}}$-instantons\footnote{By the Hitchin--Kobayashi
  correspondence, proved in \cite{Bando} in this setting, the framed
  moduli space of $G_{\mathrm{cpt}}$-instantons on $S^4/\Z_r$ is
  isomorphic to the framed moduli space of holomorphic $G$-bundles on
  $\CP^2/\Z_r = (\C^2\cup\linf)/\Z_r$.} on~$\R^4/\Z_r$ (more precisely
on $S^4/\Z_r = \R^4\cup\{\infty\}/\Z_r$), where $r$ is the level of
the corresponding representation of $G_{\aff}^\vee$.
Here $G_{\mathrm{cpt}}$ is the maximal compact subgroup of $G$.

When $G = \SL(l)$, the Uhlenbeck partial compactif\/ication in question
is the quiver variety~$\M_0$ of the af\/f\/ine type $A_{r-1}^{(1)} =
\algsl(r)_\aff$. As mentioned above, the representation of
$\algsl(r)_\aff$ can be constructed from $\M_0$ in the framework of
quiver varieties. The level of the representation is $l$. Now the
proposed conjectural link to the representation theory of $G_\aff^\vee
= \operatorname{PGL}(l)_\aff$ is provided by the theory of quiver
varieties composed with I.~Frenkel's level rank duality \cite{MR675108}
between representations of $\algsl(r)_\aff$ with level $l$ and of
$\algsl(l)_\aff$ with level $r$:
{\allowdisplaybreaks
\begin{equation*}
  \xymatrix@R=.8pc{
    {\begin{matrix}
    \text{transversal slice in}
    \\
    \text{double af\/f\/ine Grassmannian for $\SL(l)_\aff$}
    \end{matrix}}
    \ar@{.>}[r]^{\text{\cite{braverman-2007,braverman-prep}}}_{??}
    \ar@{=}[d]
    & \text{representations of $\operatorname{PGL}(l)_\aff$ of level
      $r$}
    \ar@{<->}[dd]^{\text{level-rank duality}}
\\
    \text{moduli space of $\operatorname{SU}(l)$-instantons on
      $\R^4/\Z_r$} \ar@{=}[d] &
\\
    \text{quiver variety of type $\algsl(r)_\aff$}
    \ar[r]^{\text{\cite{Na-quiver,Na-alg}}}
    & \text{representations of $\algsl(r)_\aff$ of level $l$}
   }
\end{equation*}
In} fact, the existence of this commutative diagram is one of the
sources of Braverman and Fin\-kelberg's proposal, and has been already
used to check that the intersection cohomology group of the Uhlenbeck
compactif\/ication has dimension equal to the corresponding weight space~\cite{braverman-2007}.

One of the most important ingredients in the geometric Satake
correspondence is the convolution diagram which gives the tensor
product of representations. Braverman and Finkelberg
\cite{braverman-prep} propose that its af\/f\/ine analog is the Uhlenbeck
compactif\/ication of the framed moduli space of
$G_{\mathrm{cpt}}$-instantons on the partial resolution
$X$ of $\R^4/\Z_{r_1+r_2}$ having two singularities of type
$\R^4/\Z_{r_1}$ and $\R^4/\Z_{r_2}$ connected by $\proj^1$.
Now again for $G = \SL(l)$, the Uhlenbeck compactif\/ication is the
quiver variety $\M_{\zeta^\bullet}$ (see \secref{sec:instanton}).
Since the tensor product of a level $r_1$ representation and a~level
$r_2$ representation corresponds to the restriction to
$(\algsl(r_1)\oplus \algsl(r_2))_\aff$ under the level-rank duality,
their proposal can be checked from the theory developed in this paper.

Let us explain several other things treated/not treated in this paper.
Recall that Kashiwara and Saito \cite{KS,Saito} gave a structure of
the crystal on the set of irreducible components of $\La$, which is
isomorphic to the Kashiwara's crystal of the corresponding
representation of $\Uq$.
In \secref{sec:crystal} we study the branching in the crystal when
$\M_{\zeta^\bullet}$ corresponds to a Levi subalgebra~$\g_{I^0}$ of
$\g$ corresponding to a subset $I^0$ of the index set $I$ of simple
roots. From a general theory on the crystal, the $\g_{I^0}$-crystal of
the restriction of a $\Uq$-representation to the subalgebra
$\Un(\g_{I^0})\subset \Uq$ is given by forgetting $i$-arrows with
$i\notin I^0$.
Each connected component contains an element corresponding to a
highest weight vector.
An irreducible component of $\La$ is a highest weight vector in the
$\g_{I^0}$-crystal if and only if it is mapped birationally onto its
image under $\pi_{\zeta^\bullet,\zeta}$ (see
\thmref{thm:mult_crystal}).

One can also express the branching coef\/f\/icients of the restriction
from $\Ul$ to ${\mathbf U}_q({\mathbf L}({\mathfrak g_{I_0}}))$ in
terms of IC sheaves on graded/cyclic quiver varieties, but we omit the
statements, as a reader can write down them rather obviously if he/she
knows \cite{Na-qaff} and understands \thmref{thm:branch}.

In \secref{sec:MV} we check several statements concerning the
af\/f\/ine analogs of Mirkovi\'c--Vilonen cycles for $G = \SL(r)$, proposed
by Braverman--Finkelberg \cite{braverman-prep}.

In Appendix~\ref{sec:level-rank} we review the level-rank duality following
\cite{Hasegawa,NT}. The results are probably well-known to experts,
but we need to check how the degree operators are interchanged.

After the f\/irst version of the paper was submitted, one of the
referees pointed out that the restriction to a Levi subalgebra of $\g$
was already considered by Malkin \cite[\S~3]{Malkin} at least in the
level of the crystal. Thus the result in \secref{sec:crystal} is {\it
  not} new. But the author decided to keep \secref{sec:crystal}, as
it naturally arises as a good example of the theory developed in this
paper.

\section{Partial resolutions}\label{sec:partial}

\subsection{Quiver varieties}\label{subsec:quiver}

Suppose that a f\/inite graph is given. Let $I$ be the set of vertices
and $E$ the set of edges.
Let $\bC = (c_{ij})$ be the Cartan matrix of the graph, namely
\begin{equation*}
   c_{ij} =
   \begin{cases}
     2 - 2 (\text{the number of edges joining $i$ to itself})
     & \text{if $i=j$},
\\
     - (\text{the number of edges joining $i$ to $j$})
     & \text{if $i\neq j$}.
   \end{cases}
\end{equation*}
If the graph does not contain edge loops, it is a symmetric
generalized Cartan matrix, and we have the corresponding symmetric
Kac--Moody algebra $\g$. We also def\/ine
\(
   a_{ij} = 2\delta_{ij} - c_{ij}.
\)
Then $\bA = (a_{ij})$ is the adjacency matrix when there are no edge
loops, but {\it not} in general.

In \cite{Na-quiver,Na-alg} the author assumed that the graph does not
contain edge loops (i.e., no edges joining a vertex with itself), but
most of results (in particular def\/initions, natural morphisms, etc)
hold without this assumption. And more importantly we need to consider
such cases in local models even if we study quiver varieties without
edge loops. See \subsecref{subsec:local}. So we do not assume the
condition.

Let $H$ be the set of pairs consisting of an edge together with its
orientation. So we have $\# H = 2\# E$.
For $h\in H$, we denote by $\vin{h}$ (resp.\ $\vout{h}$) the incoming
(resp.\ outgoing) vertex of~$h$.  For $h\in H$ we denote by $\overline
h$ the same edge as $h$ with the reverse orientation.
Choose and f\/ix an orientation $\Omega$ of the graph,
i.e., a subset $\Omega\subset H$ such that
$\overline\Omega\cup\Omega = H$, $\Omega\cap\overline\Omega = \varnothing$.
The pair $(I,\Omega)$ is called a {\it quiver}.

Let $V = (V_i)_{i\in I}$ be a f\/inite dimensional $I$-graded vector
space over $\C$. The dimension of $V$ is a vector
\[
  \dim V = (\dim V_i)_{i\in I}\in \mathbb Z_{\ge 0}^I.
\]
We denote the $i^{\mathrm{th}}$ coordinate vector by $\be_i$.

If $V^1$ and $V^2$ are $I$-graded vector spaces, we def\/ine vector spaces by
\begin{gather*}
  \HomL(V^1, V^2) \defeq
  \bigoplus_{i\in I} \Hom(V^1_i, V^2_i), \qquad
  \HomE(V^1, V^2) \defeq
  \bigoplus_{h\in H} \Hom(V^1_{\vout{h}}, V^2_{\vin{h}})
  .
\end{gather*}

For $B = (B_h) \in \HomE(V^1, V^2)$ and
$C = (C_h) \in \HomE(V^2, V^3)$, let us def\/ine a multiplication of $B$
and $C$ by
\[
  CB \defeq \left(\sum_{\vin{h} = i} C_h B_{\overline h}\right)_i \in
  \HomL(V^1, V^3).
\]
Multiplications $ba$, $Ba$ of $a\in \HomL(V^1,V^2)$, $b\in\HomL(V^2,
V^3)$, $B\in \HomE(V^2, V^3)$ are def\/ined in the obvious manner. If
$a\in\HomL(V^1, V^1)$, its trace $\tr(a)$ is understood as $\sum_i
\tr(a_i)$.

For two $I$-graded vector spaces $V$, $W$ with $\bv = \dim V$, $\bw =
\dim W$, we consider the vector space given by
\begin{equation*}\label{def:bM}
  \bM \equiv \bM(\bv,\bw) \equiv \bM(V, W) \defeq
  \HomE(V,V)\oplus \HomL(W,V)\oplus \HomL(V,W),
\end{equation*}
where we use the notation $\bM(\bv,\bw)$ when the isomorphism classes
of $I$-graded vector spaces $V$,~$W$ are concerned, and $\bM$ when
$V$, $W$ are clear in the context.
The dimension of $\bM$ is $\lsp{t}{\bv}(2\bw + (2\bI - \bC)\bv)$,
where $\bI$ is the identity matrix.
The above three components for an element of $\bM$ will be denoted by
$B = \bigoplus B_h $, $a = \bigoplus a_i$, $b = \bigoplus b_i$
respectively.

When the graph has no edge loop and corresponds to the symmetric
Kac--Moody Lie algeb\-ra~$\g$, we consider $\bv$, $\bw$ as weights of
$\g$ by the following rule: $\bv = \sum_i v_i \alpha_i$, $\bw = \sum_i
w_i\Lambda_i$, where $\alpha_i$ and $\Lambda_i$ are a simple root and
a fundamental weight respectively.

The orientation $\Omega$ def\/ines a function $\varepsilon\colon H \to
\{ \pm 1\}$ by $\varepsilon(h) = 1$ if $h\in\Omega$, $\varepsilon(h)=
-1$ if $h\in\overline\Omega$. We consider $\varepsilon$ as an element
of $\HomL(V, V)$.
Let us def\/ine a symplectic form $\omega$ on $\bM$ by
\begin{equation*}
        \omega((B, a, b), (B', a', b'))
        \defeq \tr(\varepsilon B\, B') + \tr(a b' - a' b).
\label{def:symplectic}\end{equation*}

Let $G \equiv G_\bv \equiv G_V$ be an algebraic group def\/ined by
\begin{equation*}
   G \equiv G_\bv \equiv G_V \defeq \prod_i \GL(V_i),
\end{equation*}
where we use the notation $G_\bv$ (resp.\ $G_V$) when we want to
emphasize the dimension (resp.\ the vector space). Its Lie algebra is
the direct sum $\bigoplus_i \gl(V_i)$.
The group $G$ acts on $\bM$ by
\begin{equation*}\label{eq:Kaction}
  (B,a,b) \mapsto g\cdot (B,a,b)
  \defeq \big(g B g^{-1}, ga, bg^{-1}\big)
\end{equation*}
preserving the symplectic structure. The space $\bM$ has a factor
\begin{equation}\label{eq:trivial}
   \bM^{\text{el}}\defeq
   \bigoplus_{h : \vout{h} = \vin{h}} \C \id_{V_{\vout{h}}}
\end{equation}
on which $G$ acts trivially. This has a $2$-dimensional space for each
edge loop, and hence has dimension $\sum_i (2 - c_{ii})$ in total.

The moment map vanishing at the origin is given by
\begin{equation}
  \mu(B, a, b) = \varepsilon B\, B + ab \in \HomL(V, V),
\end{equation}
where the dual of the Lie algebra of $G$ is identif\/ied with
$\HomL(V, V)$ via the trace.

We call a point $(B,a,b)$ in $\mu^{-1}(0)$ (or more generally
$\bM(V,W)$) a {\it module}. In fact, it is really a module of a
certain path algebra (with relations) after Crawley-Boevey's trick in
\cite[{the end of Introduction}]{CB} (see also \subsecref{subsec:CB}),
but this view point is not necessary, and a reader could consider this
is simply naming.

We would like to consider a `symplectic quotient' of $\mu^{-1}(0)$
divided by $G$.
However we cannot expect the set-theoretical quotient
to have a good property. Therefore we consider the quotient using the
geometric invariant theory. Then the quotient depends on an additional
parameter $\zeta = (\zeta_i)_{i\in I}\in \Z^I$ as follows:
Let us def\/ine a character of $G$ by
\begin{equation*}
  \chi_\zeta(g) \defeq \prod_{i\in I} \left(\det g_i\right)^{-\zeta_i}.
\end{equation*}
Let $A(\mu^{-1}(0))$ be the coordinate ring of the af\/f\/ine variety
$\mu^{-1}(0)$. Set
\begin{equation*}
   A(\mu^{-1}(0))^{G,\chi_\zeta^n}
   \defeq \{ f\in A(\mu^{-1}(0)) \mid f(g\cdot(B,a,b))
   = \chi_\zeta(g)^n f((B,a,b)) \}.
\end{equation*}
The direct sum with respect to $n\in\Z_{\ge 0}$ is a graded algebra,
hence we can def\/ine
\begin{equation*}
   \M_\zeta \equiv \M_\zeta(\bv,\bw) \equiv \M_\zeta(V,W)
   \defeq \operatorname{Proj}\Big(\bigoplus_{n\ge 0}
   A(\mu^{-1}(0))^{G,\chi_\zeta^n}\Big).
\end{equation*}
This is the {\it quiver variety\/} introduced in \cite{Na-quiver}.
Since this space is unchanged when we replace $\chi$ by a positive power
$\chi^N$ ($N > 0$), this space is well-def\/ined for $\zeta\in \Q^I$.
We call $\zeta$ a {\it stability parameter}.

When $W = 0$, the scalar subgroup $\C^*\operatorname{id}$ acts
trivially on $\bM$, so we choose the parameter $\zeta$ so that
$\sum_i \zeta_i \dim V_i = 0$, and take the `quotient' with respect to
the group $PG \defeq G/\C^*\operatorname{id}$ and the character
$\chi\colon PG\to \C^*$.

Since $G$ acts trivially on the factor \eqref{eq:trivial}, we have the
factorization
\begin{equation}\label{eq:fac}
    \M_\zeta = \bM^{\text{el}} \times \M_\zeta^{\text{norm}},
\end{equation}
where $\M_\zeta^{\text{norm}}$ is the symplectic quotient of the space
of datum $(B,a,b)$ satisfying $\tr(B_h) = 0$ for any $h$ with $\vin{h}
= \vout{h}$.

\subsection{Stability}\label{subsec:stability}

We will describe $\M_\zeta$ as a moduli space. We also introduce
Harder--Narasimhan and Jordan--H\"older f\/iltrations, for which we need
to add an additional variable $\zeta_\infty\in\Q$ to the stability
parameter $\zeta\in \Q^I$. We write $\Tilde\zeta =
(\zeta,\zeta_\infty) \in \Q^{I\sqcup\{\infty\}}$. Let
\begin{gather*}
 \Tilde\zeta(V,W)
   \defeq \sum_{i\in I} \zeta_{i} \dim V_i
   + \zeta_\infty (1-\delta_{W0}),
\\
 \theta_{\Tilde\zeta}(V,W)
   \defeq
   \frac{\Tilde\zeta(V,W)}{1 - \delta_{W0} + \sum_{i\in I} \dim V_i},
\end{gather*}
where $\delta_{W0}$ is $1$ if $W=0$ and $0$ otherwise, and we
implicitly assume $V\neq 0$ or $W\neq 0$ in the def\/inition of
$\theta_{\Tilde\zeta}(V,W)$.

\begin{Definition}\label{def:stable}
  \textup{(1)} Suppose a module $(B, a, b)\in\bM(V,W)$ is given. We
  consider an $I$-graded subspace $V'$ in $V$ such that either of the
  following two conditions is satisf\/ied:
\begin{aenume}\itemsep=0pt
\item $V'$ is contained in $\Ker b$ and $B$-invariant,
\item $V'$ contains $\Ima a$ and is $B$-invariant.
\end{aenume}
In the f\/irst case we def\/ine the {\it submodule\/} $(B,a,b)|_{(V',0)} \in
\bM(V',0)$ by putting $B|_{V'}$ on $\HomE(V',V')$ (and $0$ on
$\HomL(0,V')\oplus \HomL(V',0) = 0$). We def\/ine the {\it quotient
  module\/} $(B,a,b)|_{(V/V',W)}\!\in\! \bM(V/V',W)\!$ by putting
the homomorphisms induced from $B$, $a$, $b$ on
$\HomE(V/V',V/V')$, $\HomL(W,V/V')$, $\HomL(V/V',W)$ respectively.
In the second case we def\/ine the submodule
$(B,a,b)|_{(V',W)}\in \bM(V',W)$ and the quotient module
$(B,a,b)|_{(V/V',0)}\in \bM(V/V',0)$ in a similar way.
We may also say $(V',0)$ \textup(resp.\ $(V',W)$\textup), $(V/V',W)$
\textup(resp.\ $(V/V',0)$\textup) a submodule and a quotient module
of $(B,a,b)$ in the case (a) \textup(resp.\ (b)\textup).
When we want to treat the two cases simultaneously we write
a submodule $(V',\delta W)$ or a quotient module $(V, W)/(V',\delta W)$,
where we mean $\delta W$ is either $0$ or $W$.

\textup{(2)} A module $(B, a, b)\in\bM(V,W)$ is {\it
$\Tilde\zeta$-semistable\/} if we have
\begin{equation}\label{eq:slope}
   \theta_{\Tilde\zeta}(V',\delta W) \le \theta_{\Tilde\zeta}(V,W),
\end{equation}
for any nonzero submodule $(V',\delta W)$ of $(B,a,b)$.

We say $(B,a,b)$ is {\it $\Tilde\zeta$-stable\/} if the strict
inequalities hold unless $(V',\delta W) = (V,W)$.

We say $(B,a,b)$ is {\it $\Tilde\zeta$-polystable\/} if it is a
direct sum of $\Tilde\zeta$-stable modules having the same
$\theta_{\Tilde\zeta}$-value.
\end{Definition}

The function $\theta_{\Tilde\zeta}$ is an analog of the slope of a
torsion free sheaf appearing in the def\/inition of its stability. We
have the following property analogous to one for the slope.
\begin{Lemma}\label{lem:slope}
  Let $(V',\delta W)$ be a submodule of $(B,a,b)$ and
  $(V,W)/(V',\delta W)$ be the quotient. Then
  \begin{equation*}
     \theta_{\Tilde\zeta}(V',\delta W)
     \le (\text{resp.} \ge, =)
     \theta_{\Tilde\zeta}(V,W)
     \Longleftrightarrow
     \theta_{\Tilde\zeta}((V,W)/(V',\delta W))
     \ge (\text{resp.} \le, =)
     \theta_{\Tilde\zeta}(V,W).
  \end{equation*}
\end{Lemma}

The $\Tilde\zeta$-(semi)stability condition is unchanged even if we
shift the stability parameter $\Tilde\zeta$ by a vector
$c(1,1,\dots,1)\in \Q^{I\sqcup\{\infty\}}$ with a constant $c\in \Q$.
Therefore we may normalize so that $\theta_{\Tilde\zeta}(V,W) = 0$
by choosing $c = -\theta_{\Tilde\zeta}(V,W)$.
We then take the component $\zeta\in\Q^I$ after this normalization and
then def\/ine $\M_\zeta$ as in the previous subsection.
Moreover if $W=0$, the additional component $\zeta_\infty$ is clearly
irrelevant, and the normalization condition
$\theta_{\Tilde\zeta}(V,W) = 0$ just means $\sum_i \zeta_i \dim V_i
= 0$. Therefore we can also apply the construction in the previous
subsection.

Conversely when $\zeta\in\Q^I$ and $I$-graded vector spaces $V$, $W$
are given as in the previous subsection, we def\/ine $\Tilde\zeta$ by
the following convention: If $W \neq 0$, take $\zeta_\infty$ so that
$\theta_{\Tilde\zeta}(V,W) = 0$. If $W = 0$, then we have assumed
$\sum \zeta_i \dim V_i = 0$. So we just put $\zeta_\infty = 0$.
Once this convention becomes clear, we say
$\zeta$-(semi)stable instead of $\Tilde\zeta$-(semi)stable.

\begin{Example}\label{ex:stndcham}
(1)
Under this convention and the assumption $\zeta_{i} > 0$ for all $i\in I$,
the inequali\-ty~\eqref{eq:slope} is {\it never\/} satisf\/ied for a
nonzero submodule of the form $(V',0)$. Also
\eqref{eq:slope} is {\it always} satisf\/ied for a submodule
$(V',W)$.
Therefore the $\zeta$-stability is equivalent to the nonexistence of
nonzero $B$-invariant $I$-graded subspaces $V' = \bigoplus V'_i$
contained in $\Ker b$
(and in this case $\zeta$-stability and $\zeta$-semistability
are equivalent). This is the stability condition used in
\cite[{3.9}]{Na-alg}.

(2) Let $\zeta_i = 0$ for all $i$. Then any module is
$\zeta$-semistable. A module is $\zeta$-stable if and only if it is
simple, i.e., has no nontrivial submodules.
\end{Example}

We recall Harder--Narasimhan and Jordan--H\"older f\/iltrations. We need
to f\/ix $\Tilde\zeta\in \Q^{I\sqcup\{\infty\}}$ and do not take the
normalization condition $\theta_{\Tilde\zeta} = 0$ as we want to
compare $\theta_{\Tilde\zeta}$-values for various dimension vectors.

\begin{Theorem}[\cite{Rudakov}]\label{thm:Rudakov}
\textup{(1)} A module $(B,a,b)$ has the unique Harder--Narasimhan
filtration\textup:
a~flag of $I$-graded subspaces
\begin{equation*}
  V = V^0 \supset V^1 \supset \cdots \supset V^N \supset V^{N+1} = 0
\end{equation*}
and an integer $0 \le k_W \le N$ such that
\begin{aenume}\itemsep=0pt
\item $(B,a,b)|_{(V^{k+1}, \delta_{k+1} W)}$ is a submodule of
  $(B,a,b)|_{(V^k,\delta_k W)}$ for $0\le k \le N$,
where $\delta_k W = W$ for $0\le k\le k_W$ and $0$
for $k_W + 1\le k \le N+1$,

\item the quotient module $\gr_k (B,a,b) \defeq
  (B,a,b)|_{(V^k,\delta_k W)/(V^{k+1},\delta_{k+1} W)}$ is
  $\Tilde\zeta$-semistable for $0\le k\le N$ and
\[
    \theta_{\Tilde\zeta}(\gr_0(B,a,b))
    < \theta_{\Tilde\zeta}(\gr_1(B,a,b)) < \cdots
    < \theta_{\Tilde\zeta}(\gr_N(B,a,b)).
\]
\end{aenume}

\textup{(2)} A $\Tilde\zeta$-semistable module $(B,a,b)$ has a
Jordan--H\"older filtration\textup:
a f\/lag of $I$-graded subspaces
\begin{equation*}
  V = V^0 \supset V^1 \supset \cdots \supset V^N \supset V^{N+1} = 0
\end{equation*}
and an integer $0 \le k_W \le N$ such that
\begin{aenume}\itemsep=0pt
\item the same as \textup{(a)} in \textup{(1)},

\item the quotient module $\gr_k (B,a,b)\defeq (B,a,b)|_{(V^k,\delta_k
    W)/(V^{k+1},\delta_{k+1} W)}$ is $\Tilde\zeta$-stable for $0\le
  k\le N$ and
\[
    \theta_{\Tilde\zeta}(\gr_0(B,a,b))
    = \theta_{\Tilde\zeta}(\gr_1(B,a,b)) = \cdots
    = \theta_{\Tilde\zeta}(\gr_N(B,a,b)).
\]
\end{aenume}
Moreover the isomorphism class of $\bigoplus \gr_k(B,a,b)$ is uniquely
determined by that of $(B,a,b)$.
\end{Theorem}

Let $H_{\zeta}^{\operatorname{s}}$ (resp.\
$H_{\zeta}^{\operatorname{ss}}$) be the set of
$\zeta$-stable (resp.\ $\zeta$-semistable) modules in
$\mu^{-1}(0)\subset\bM$.

We say two $\zeta$-semistable modules $(B,a,b)$, $(B',a',b')$
are {\it $S$-equivalent\/} when the closures of orbits intersect in
$H_{\zeta}^{\operatorname{ss}}$.

\begin{Proposition}\label{prop:dim}
  \textup{(1) (\cite{King})} $\M_\zeta$ is a coarse moduli space of
  $\zeta$-semistable modules modulo $S$-equivalences.
  Moreover, there is an open subset $\Mreg_\zeta\subset\M_\zeta$,
  which is a fine moduli space of $\zeta$-stable modules modulo
  isomorphisms.

  \textup{(2) (\cite{King})} Two $\zeta$-semistable modules are
  $S$-equivalent if and only if their Jordan--H\"older filtration
  in Theorem~{\rm \ref{thm:Rudakov}(2)} have the same composition factors. Thus
  $\M_\zeta$ is a coarse moduli space of $\zeta$-polystable modules
  modulo isomorphisms.

  \textup{(3) (See \cite[2.6]{Na-quiver} or \cite[{3.12}]{Na-alg})}
  Suppose $\bw \neq 0$. Then
  $\Mreg_\zeta$, provided it is nonempty, is smooth of dimension
\(
   \lsp{t}{\bv}(2\bw - \bC \bv).
\)
If $\bw = 0$, then $\Mreg_\zeta$, provided it is nonempty, is smooth
of dimension
\(
   2 - \lsp{t}{\bv}\bC \bv.
\)
\end{Proposition}

The $S$-equivalence class of a point $(B,a,b)\in\mu^{-1}(0)\cap
H_{\zeta}^{\operatorname{ss}}$ will be denoted by $[B,a,b]$, and is
considered as a closed point in $\M_\zeta$.

We recall the following result \cite[3.1, 4.2]{Na-quiver}:

\begin{Theorem}\label{thm:affine}
  Suppose $\Mreg_\zeta = \M_\zeta$. Then $\M_\zeta$ is diffeomorphic
  to a nonsingular affine algebraic variety.
\end{Theorem}

\subsection{Face structure on the set of stability conditions}

Fix a dimension vector $\bv = (v_i)$. Let
\begin{gather*}
   R_+  \defeq \{ \theta = (\theta_i)_{i\in I} \in \Z_{\ge 0}^I \mid
     {}^t\theta \bC \theta \le 2\} \setminus \{0\}, \\
   R_+(\bv) \defeq \{ \theta\in R_+ \mid \theta_i \le v_i
       \quad\text{for all $i\in I$}\}, \\
   D_\theta \defeq \{ \zeta = (\zeta_i)\in \Q^I \mid
   \zeta\cdot \theta = 0 \} \quad\text{for $\theta\in R_+$},
\end{gather*}
where $\zeta\cdot \theta = \sum_i \zeta_i\theta_i$. When the graph is
of Dynkin or af\/f\/ine type, $R_+$ is the set of positive roots, and
$D_\theta$ is the wall def\/ined by the root $\theta$
(\cite[Proposition~5.10]{Kac}).
In general, $R_+$ may be an inf\/inite set, but $R_+(\bv)$ is always
f\/inite.

\begin{Lemma}[\protect{cf.\ \cite[2.8]{Na-quiver}}]\label{lem:nondeg}
  Suppose $\Mreg_\zeta(\bv,0)\neq\varnothing$. Then $\bv \in R_+$.
\end{Lemma}

This is clear from the dimension formula \propref{prop:dim}(3) as
$\dim \Mreg_\zeta(\bv,0)$ must be nonnegative.

The set $R_+(\bv)$ def\/ines a system of {\it faces\/}. A subset
$F\subset \R^I$ is a {\it face\/} if there exists a disjoint
decomposition $R_+(\bv) = R_+^0(\bv) \sqcup R_+^+(\bv) \sqcup
R_+^-(\bv)$ such that
\begin{equation*}
    F = \{ \zeta\in \R^I \mid \text{
     $\zeta\cdot \theta = 0$ (resp.\ $> 0$, $< 0$) for
    $\theta\in R_+^0(\bv)$ (resp.\ $R_+^+(\bv)$, $R_+^-(\bv)$)}\}.
\end{equation*}
When we want to emphasize that it depends on the dimension vector
$\bv$, we call it a {\it $\bv$-face}.
A~face is an open convex cone in the subspace $\{ \zeta\in\R^I \mid
\text{$\zeta\cdot\theta = 0$ for $\theta\in R_+^0(\bv)$}\}$. A face is
called a {\it chamber\/} (or {\it $\bv$-chamber\/}) if it is an open
subset in $\R^I$, i.e., $R_+^0(\bv) = \varnothing$. The closure of~$F$
is
\begin{equation*}
    \overline{F} = \{ \zeta\in \R^I \mid \text{
     $\zeta\cdot \theta = 0$ (resp.\ $\ge 0$, $\le 0$) for
    $\theta\in R_+^0(\bv)$ (resp.\ $R_+^+(\bv)$, $R_+^-(\bv)$)}\}.
\end{equation*}

\begin{Lemma}\label{lem:face}
  Suppose $W\neq 0$.

  \textup{(1) (cf.\ \cite[2.8]{Na-quiver})} Suppose $\zeta$ is in a
  chamber. Then $\zeta$-semistability implies $\zeta$-stability, and hence
  $\Mreg_\zeta(\bv,\bw) = \M_\zeta(\bv,\bw)$.

  \textup{(2) (cf.\ \cite[1.4]{Na-ADHM})} If two stability parameters
  $\zeta$, $\zeta'$ are contained in the same face $F$,
  $\zeta$-stability \textup(resp.\ $\zeta$-semistability\textup) is
  equivalent to $\zeta'$-stability \textup(resp.\
  $\zeta'$-semistability\textup).

  \textup{(3)} Suppose stability parameters $\zeta$, $\zeta^\bullet$ are
  contained in $F$ and $F^\bullet$ respectively. If $F^\bullet\subset
  \overline{F}$, then
  \begin{aenume}\itemsep=0pt
    \item a $\zeta$-semistable module is $\zeta^\bullet$-semistable,
    \item a $\zeta^\bullet$-stable module is $\zeta$-stable.
  \end{aenume}
\end{Lemma}

\begin{proof}
(1) For given $V$,$W$, we def\/ine $\Tilde\zeta=(\zeta,\zeta_\infty)$
with $\theta_{\Tilde\zeta}(V,W) = 0$ as before.
Suppose that $(B,a,b)$ is $\Tilde\zeta$-semistable, but not
$\Tilde\zeta$-stable. We take a Jordan--H\"older f\/iltration as in
Theo\-rem~\ref{thm:Rudakov}(2). Since $\delta_k W/\delta_{k+1} W \neq 0$ only
for one $k$ in $0,\dots, N$, we have a $k$ with $\delta_k
W/\delta_{k+1} W = 0$ from the assumption $N\ge 1$. As the quotient
module $\gr_k(B,a,b)\! =\! (B,a,b)|_{(V^k,\delta_k
  W)/\linebreak[2](V^{k+1}, \delta_{k+1} W)}$ is
$\Tilde\zeta$-stable, its dimension vector $\dim V^k - \dim V^{k+1}$
(with the $\infty$-component is $0$) is in $R_+(\bv)$ by
\lemref{lem:nondeg}. Moreover, by the condition in the Jordan--H\"older
f\/iltration, we have
\[
  \theta_{\Tilde\zeta}(\gr_k(B,a,b)) = \theta_{\Tilde\zeta}(V,W).
\]
The right hand side is equal to $0$ by our convention, and the left
hand side is equal to $\zeta\cdot(\dim V^k - \dim V^{k+1})$ up to
scalar. This contradicts with our assumption that $\zeta$ is in a chamber.

(2) We def\/ine $\Tilde\zeta$, $\Tilde\zeta'$ as before.
Suppose that $(B,a,b)\in\mu^{-1}(0)\subset \bM(\bv,\bw)$ is
$\Tilde\zeta$-semistable and is not $\Tilde\zeta'$-semistable.
We take the Harder--Narasimhan f\/iltration for $(B,a,b)$ with respect to
the $\Tilde\zeta'$-stability as in \thmref{thm:Rudakov}(1). We have
$N\ge 1$ from the assumption.

Consider f\/irst the case when there exists $0\le k\le k_W - 1$ with
\(
   \theta_{\Tilde\zeta'}(\gr_k(B,a,b)) < 0.
\)
Then
\(
   \theta_{\Tilde\zeta'}(\gr_l(B,a,b)) < 0
\)
for any $0\le l\le k$.
On the other hand from the $\zeta$-semistability of $(B,a,b)$, we have
\(
   \theta_{\Tilde\zeta}((V,W)/(V^{k},\delta_{k} W)) \ge 0
\)
by \lemref{lem:slope}.
Therefore there exists at least one $l$ in $[0, k]$ with
\(
   \theta_{\Tilde\zeta}(\gr_l(B,a,b)) \ge 0.
\)
We further take a Jordan--H\"older f\/iltration of $\gr_l(B,a,b)$ to f\/ind
$\Tilde\zeta'$-stable representations $V'$ with
\(
  \nicefrac{\zeta'\cdot\dim V'}{\sum \dim V'_i} =
  \theta_{\Tilde\zeta'}(V',0) = \theta_{\Tilde\zeta'}(\gr_l(B,a,b)) < 0.
\)
Note that $\gr_l(B,a,b)$, and hence $V'$ has $0$ in
the $W$-component.
Since
\(
   \theta_{\Tilde\zeta}(\gr_l(B,a,b)) \ge 0,
\)
we can take $V'$ so that
\(
    \nicefrac{\zeta\cdot\dim V'}{\sum \dim V'_i}
    = \theta_{\Tilde\zeta}(V',0) \ge 0.
\)
As $\dim V'\in R_+(\bv)$, this contradicts with the assumption that
$\zeta$ and $\zeta'$ are in the common face.

The same argument leads to a contradiction in the case when
there exists $k_W+1 \le k \le N$ with
\(
   \theta_{\Tilde\zeta'}(\gr_k(B,a,b)) > 0.
\)
Since $N\ge 1$ and
\(
   \theta_{\Tilde\zeta'}(V,W) = 0,
\)
at least one of two cases actually occur. Therefore $(B,a,b)$ is
$\zeta'$-semistable.

Suppose further that $(B,a,b)$ is $\zeta$-stable. We want to show that
it is also $\zeta'$-stable. Assume not, and take a Jordan--H\"older
f\/iltration of $(B,a,b)$ with respect to the $\zeta'$-stability. Either
of $\gr_0(B,a,b)$ and $\gr_N(B,a,b)$ has the $W$-component $0$. Suppose
the f\/irst one has the $W$-component $0$. We have
\(
   \nicefrac{\zeta'\cdot \dim (V/V^1)}{\sum \dim V_i/V^1_i} =
   \theta_{\Tilde\zeta'}(\gr_0(B,a,b)) = 0.
\)
On the other hand, the $\zeta$-stability of $(B,a,b)$ implies
\(
   \nicefrac{\zeta\cdot \dim (V/V^1)}{\sum \dim V_i/V^1_i} =
   \theta_{\Tilde\zeta}(\gr_0(B,a,b)) > 0.
\)
This contradicts with the assumption. The same argument applies
to the case when $\gr_N(B,a,b)$ has the $W$-component~$0$. Therefore
$(B,a,b)$ is $\zeta'$-stable.

(3) Take a sequence $\{ \zeta_n \}$ in $F$ converging to $\zeta^\bullet$.
Take a nonzero submodule $(S,\delta W)$ of $(B,a,b)$.
(a) From the $\zeta$-semistability and (2), we have
\(
   \theta_{\Tilde\zeta_n}(S,\delta W) \le
   \theta_{\Tilde\zeta_n}(V, W)
\)
for any $n$. Taking limit, we get
\(
   \theta_{\Tilde\zeta^\bullet}(S,\delta W) \le
   \theta_{\Tilde\zeta^\bullet}(V,W).
\)
Therefore $(B,a,b)$ is $\zeta^\bullet$-semistable.
(b) We assume $(S,\delta W) \neq (V,W)$.
From the $\zeta^\bullet$-stability, we have
\(
   \theta_{\Tilde\zeta^\bullet}(S,\delta W) <
   \theta_{\Tilde\zeta^\bullet}(V,W).
\)
Then
\(
   \theta_{\Tilde\zeta_n}(S,\delta W) <
   \theta_{\Tilde\zeta_n}(V,W)
\)
for suf\/f\/iciently large $n$. Therefore $(B,a,b)$ is
$\zeta_n$-stable. By (2) it is also $\zeta$-stable.
\end{proof}

\begin{Remark}\label{rem:W=0reg}
  When $W=0$, we need to modify the def\/initions of faces and chambers:
  Replace $\R^I$ by $\{ \zeta \mid \zeta\cdot\bv = 0\}$, and
  $R_+(\bv)$ by $R_+(\bv)\setminus\{\bv\}$. Then the same proof works.
\end{Remark}

From this lemma, we can def\/ine the $\zeta$-(semi)stability for
$\zeta\in \R^I$, not necessarily in $\Q^I$, as it depends only on the
face containing $\zeta$.

\subsection[Nonemptiness of ${\mathfrak M}^{\rm s}_\zeta$]{Nonemptiness of $\boldsymbol{\Mreg_\zeta}$}\label{subsec:CB}

For most of purposes in our paper, \lemref{lem:nondeg} is enough, but
we can use a rotation of the hyper-K\"ahler structure and then apply
Crawley-Boevey's result \cite{CB} to get a necessary and suf\/f\/icient
condition for $\Mreg_\zeta(\bv,\bw)\neq\varnothing$. This will be given
in this subsection.

Suppose that the graph $(I,E)$ does not have edge loops. We thus
associate a Kac--Moody Lie algebra $\g$ to $(I,E)$.
We f\/ix a dimension vector $\bw = (w_i)\in \Z_{\ge 0}^I$, which is
considered as a~dominant weight for $\g$.
Let $(\Tilde{I}, \Tilde{E})$ be a graph obtained from $(I,E)$ by
adding a new vertex $\infty$ and $w_i$ edges between $\infty$ and $i$.
The new graph $(\Tilde{I}, \Tilde{E})$ def\/ines another Kac--Moody Lie
algebra which we denote by $\Tilde{\g}$. We denote the corresponding
Cartan matrix by $\Tilde{\bC}$. This new quiver was implicitly used in
\subsecref{subsec:stability}.

\begin{Lemma}\label{lem:root}
  Let $\alpha_\infty$ denote the simple root corresponding to the
  vertex $\infty$. Let $V$ be the direct sum of root spaces
  $\Tilde{\g}_\alpha$ of $\Tilde{\g}$, where $\alpha$ is of the form
  $\alpha = - \sum_{i\in I} m_i \alpha_i - \alpha_\infty$. Let
  $\g$ act on $V$ by the restriction of the adjoint representation.
  Then $V$ is isomorphic to the irreducible integrable highest weight
  representation $V(\bw)$ with the highest weight vector $f_\infty$ of
  weight $\bw$. In particular, the root multiplicity $\dim
  \Tilde{\g}_\alpha$ is equal to the weight multiplicity of
  $\alpha|_{\mathfrak h}$ in $V(\bw)$, where $\mathfrak h$ is the
  Cartan subalgebra of $\g$.
\end{Lemma}

\begin{proof}
  From the def\/inition we have
  \begin{equation*}
    [e_i, f_\infty] = 0, \quad
    [h_i, f_\infty] = w_i f_\infty \qquad\text{for $i\in I$}.
  \end{equation*}
  Moreover $\Tilde{\g}_\alpha$ is a linear span of the elements of the
  form
\begin{gather*}
  [f_{i_1}, [f_{i_2}, [ \dots\!,  [f_{i_{k-1}}\!, [ f_\infty,
  \underbrace{[f_{i_k}, [\dots
  [ f_{i_{s-1}}\!, f_{i_s}]\dots]]}_{\ =x}]] \dots]]]
  = - [f_{i_1}, [f_{i_2}, [ \dots\!, [ f_{i_{k-1}}\!, [x, f_\infty]] \dots ]]]
\end{gather*}
for $i_1, i_2, \ldots\in I$. This means that $V$ is a highest
weight module.
We also know $V$ is integrable by \cite[3.5]{Kac}. Now the assertion
follows from \cite[10.4]{Kac}.
\end{proof}

\begin{Theorem}\label{thm:CB}
  \textup{(1)} Suppose $W = 0$ and consider a dimension vector $\bv$
  with $\zeta\cdot \bv = 0$. Then $\Mreg_\zeta(\bv,0)\neq\varnothing$ if
  and only if the following holds\textup:
\begin{itemize}\itemsep=0pt
\item $\bv$ is a positive root, and $p(\bv) > \sum_{t=1}^r
  p(\beta^{(t)})$ for any decomposition $\bv = \sum_{t=1}^r
  \beta^{(t)}$ with $r \ge 2$ and $\beta^{(t)}$ a positive root with
  $\zeta\cdot \beta^{(t)} = 0$ for all $t$,
\end{itemize}
where $p(x) = 1 - \frac12 {}^t x \bC x$.

\textup{(2)} Suppose $\bw\neq 0$. Then
$\Mreg_\zeta(\bv,\bw)\neq\varnothing$ if and only if the following
holds:
\begin{itemize}\itemsep=0pt
\item $\bw - \bv$ is a weight of the integrable highest weight
  representation $V(\bw)$ of the highest weight $\bw$, and
  $\lsp{t}{\bv}(\bw - \frac12 \bC \bv) > \lsp{t}{\bv^0}(\bw - \frac12
  \bC \bv^0) + \sum_{t=1}^r p(\beta^{(t)})$ for any decomposition $\bv
  = \bv^0 + \sum_{t=1}^r \beta^{(t)}$ with $r \ge 1$, $\bw - \bv^0$ is
  a weight of $V(\bw)$, and $\beta^{(t)}$ a positive root with
  $\zeta\cdot \beta^{(t)} = 0$ for all $t$.
\end{itemize}
\end{Theorem}

\begin{Remark}
For more general pairs of parameters $(\zeta,\zeta_\C)$, it seems
natural to expect that the existence of $\zeta$-stable module $B$
with $\mu(B,a,b) = \zeta_\C$ is equivalent to the above condition with
$\zeta\cdot \beta^{(t)} = 0$ replaced by
`$\zeta\cdot \beta^{(t)} = 0$ and $\zeta_\C\cdot\beta^{(t)} = 0$'.
\end{Remark}

\begin{proof}
(1) A $\zeta$-stable module with $\mu = 0$ corresponds to a point
satisfying the hyper-K\"ahler moment map $(\mu_\R,\mu) = (\zeta,0)$,
which correspond a $0$-stable module (i.e., a simple module) with $\mu
= \zeta$ after a rotation of the complex structures. Then by
\cite[Theorem~1.2]{CB} such a simple module exists if and only if the
above condition holds.

(2) For the given $\bw$ we construct the quiver $(\Tilde{I},\Tilde{E})$ as
above. We consider the dimension vector $\Tilde{\bv} = \bv + \alpha_\infty$
and def\/ine the stability parameter $\Tilde{\zeta} = (\zeta,\zeta_\infty)$ by
\(
\zeta\cdot \dim \bv + \zeta_\infty = 0.
\)
Then $\Mreg_\zeta(\bv,\bw)$ is isomorphic to the quiver variety
$\Mreg_{\Tilde{\zeta}}(\Tilde{\bv},0)$ associated with
$(\Tilde{I},\Tilde{E})$ \cite[the end of Introduction]{CB}.
Therefore we can apply the criterion in (1) together with the
observation \lemref{lem:root}. We f\/irst note that
\[
  p(\Tilde{\bv}) = 1 - \tfrac12 \lsp{t}{\Tilde{\bv}} \Tilde{\bC} \Tilde{\bv}
  = \lsp{t}{\bv}\left(\bw - \tfrac12 \bC \bv\right).
\]
Also if we have a decomposition $\Tilde{\bv} = \sum \beta^{(t)}$, one
of $\beta^{(t)}$ has $1$ in the entry $\infty$, and other
$\beta^{(t)}$'s have $0$ in the entry $\infty$. We rewrite the former
as $\bv^0 + \alpha_\infty$, and identify the latter with positive
roots for $(I,E)$. Now the assertion follows from (1).
\end{proof}

\subsection{Partial resolutions}\label{subsec:partial}

Let $\zeta,\zeta^\bullet\in\R^I$ as in \lemref{lem:face}(3). Then we have a
morphism
\begin{equation*}
   \pi_{\zeta^\bullet,\zeta} \colon \M_\zeta \to \M_{\zeta^\bullet}
\end{equation*}
thanks to \lemref{lem:face}(3)(a). By \lemref{lem:face}(3)(b) it is an
isomorphism on the preimage $\pi_{\zeta^\bullet,\zeta}^{-1}(\Mreg_{\zeta^\bullet})$.

If further $\zeta^\circ$, $\zeta^\bullet$ are as in \lemref{lem:face}(3),
$\zeta$, $\zeta^\bullet$ are also, and we have
\(
  \pi_{\zeta^\bullet,\zeta} = \pi_{\zeta^\bullet,\zeta^\circ}
  \circ \pi_{\zeta^\circ,\zeta}.
\)

Since $0$ is always in the closure of any face, we always have a
morphism
\(
  \pi_{0,\zeta^\bullet}\colon \M_{\zeta^\bullet}\to \M_0
\)
for any $\zeta^\bullet$. There always exists a {\it chamber\/} $\mathcal C$
containing $\zeta^\bullet$ in the closure. If we take a parameter $\zeta$
from $\mathcal C$, we have
\(
  \pi_{\zeta^\bullet,\zeta}\colon \M_{\zeta}\to \M_{\zeta^\bullet}.
\)
And we have $\pi_{0,\zeta} =
\pi_{0,\zeta^\bullet}\circ\pi_{\zeta^\bullet,\zeta}$. Since
$\zeta$ is in a chamber, we have $\M_{\zeta} = \Mreg_{\zeta}$, and
hence $\M_{\zeta}$ is nonsingular. It is known that
$\M_{\zeta}$ is a resolution of singularities of $\M_0$ in many cases
(e.g., if $\Mreg_{\zeta}\neq\varnothing$ (see
\cite[Theorem~4.1]{Na-quiver})). In these cases, $\M_{\zeta^\bullet}$ is a~{\it
  partial} resolution of singularities of $\M_0$.

When the stability parameters are apparent from the context, we denote
the map $\pi_{\zeta^\bullet,\zeta}$ simply by $\pi$.

\subsection{Stratum}\label{subsec:stratum}

We recall the stratif\/ication on $\M_\zeta$ considered in
\cite[Section~6]{Na-quiver}, \cite[Section~3]{Na-alg}.

Suppose that $(B,a,b)$ is $\zeta$-polystable. Let us decompose it as
\begin{gather}
     V \cong V^0 \oplus (V^1)^{\oplus {\widehat v}_1}\oplus \cdots \oplus
   (V^r)^{\oplus {\widehat v}_r},
\nonumber\\
     (B,a,b) \cong (B^0,a^0,b^0) \oplus
   (B^1)^{\oplus {\widehat v}_1} \oplus \cdots \oplus
   (B^r)^{\oplus {\widehat v}_r},\label{eq:decomp}
\end{gather}
where $(B^0,a^0,b^0)\in\mu^{-1}(0)\cap\bM(V^0,W)$ is the unique factor
having $W\neq 0$, and $B^k\in \mu^{-1}(0)\cap\HomE(V^k,V^k)$
($k=1,\dots,r$) are pairwise non-isomorphic $\zeta$-stable modules
and ${\widehat v}_k$ is its multiplicity in $(B,a,b)$. (See
\cite[6.5]{Na-quiver}, \cite[3.27]{Na-alg}.)

The stabilizer $\widehat G$ of $(B,a,b)$ is conjugate to
\[
   \{ \id_{V^0} \} \times
   \prod_{k=1}^r \left(\GL(\widehat{v}_k)\otimes \id_{V^k}\right).
\]
Conversely if the stabilizer is conjugate to this subgroup, the module
has the decomposition above.

We thus have a stratif\/ication by $G$-orbit type:
\begin{equation*}
   \M_\zeta = \bigsqcup_{(\widehat G)} (\M_\zeta)_{(\widehat G)},
\end{equation*}
where $(\M_\zeta)_{(\widehat G)}$ consists of modules whose
stabilizers are conjugate to a subgroup $\widehat G$ of $G$.

A list of strata can be given by \thmref{thm:CB} in principle, but to
use this result, we need to know all roots, and it is not so easy in
general.
The following def\/inition was considered in \cite[2.6.4]{Na-qaff} to
avoid this dif\/f\/iculty and concentrate only on the stratum with $\bw
\neq 0$.

\begin{Definition}
  A stratum $(\M_\zeta)_{(\widehat G)}$ is {\it regular\/} if it is of
  the form $\Mreg_\zeta(\bv',\bw)\times\{ \bigoplus_{i\in I}
  S_i^{\oplus (v_i - v'_i)}\}$ for a $\bv' = (v'_i)$, where $S_i$ is
  the module with $\C$ on the vertex $i$, and $0$ on the other
  vertices and $B = 0$. Here we assume $\zeta_i = 0$.

  A point $x\in\M_\zeta(\bv,\bw)$ is {\it regular\/} if it is contained in
  a regular stratum.
\end{Definition}

If the graph is of f\/inite type, all strata in $\M_0(\bv,\bw)$ are
regular. This was proved in \cite[6.7]{Na-quiver}, but it also follows
from \thmref{thm:CB}, as $p(x) = 0$ for a real root $x$.

\subsection{Local structure}\label{subsec:local}

Let $\zeta$, $\zeta^\bullet$ be as in \lemref{lem:face}(3).
We examine the local structure of $\M_{\zeta^\bullet}$ (resp.\ $\M_{\zeta}$)
around a~point $x$ (resp.\ $\pi_{\zeta^\bullet,\zeta}^{-1}(x)$) in this
subsection. This has been done in \cite[Section~6]{Na-quiver},
\cite[Section~3.2]{Na-qaff}, but we explain the results in a slightly
dif\/ferent form.

Let us take a representative $(B,a,b)$ of $x$, which is
$\zeta^\bullet$-polystable. We consider the complex (see
\cite[(3.11)]{Na-alg})
\begin{gather}
   \mathscr{C}^\bullet :
   \HomL(V,V) \xrightarrow{\alpha}
   \HomE(V,V)\oplus\HomL(W,V) \oplus \HomL(V,W) \xrightarrow{\beta}
   \HomL(V,V),
\nonumber\\
   \alpha(\xi) = (B\xi-\xi B)\oplus (-\xi a) \oplus (b\xi), \qquad
   \beta(C,d,e) = \varepsilon B C + \varepsilon C B + a e + d b,\label{eq:tancpx}
\end{gather}
where $\alpha$ is the inf\/initesimal action of the Lie algebra of $G_V$
on $\bM$, and $\beta$ is the dif\/ferential of the moment map $\mu$ at
$(B,a,b)$.

Suppose that $(B,a,b)$ is decomposed into a direct sum of
$\zeta^\bullet$-stable modules as in~\eqref{eq:decomp}.
Then $\Ker\beta/\Ima\alpha$ is the symplectic normal denoted by
$\widehat\bM$ in \cite[Section~6]{Na-quiver}, \cite[Section~3.2]{Na-qaff}.
The complex $\mathscr{C}^\bullet$ decomposes as $\mathscr{C}^\bullet =
\bigoplus_{k,l=0}^r (\mathscr{C}_{k,l}^\bullet)^{\oplus {\widehat v}_k {\widehat
    v}_l}$ with
\begin{equation*}
   \mathscr{C}_{k,l}^\bullet :
   \HomL(V^k, V^l)
   \xrightarrow{\alpha}
   \HomE(V^k, V^l)
   \oplus \HomL(W, V^l)^{\oplus \delta_{k0}}
   \oplus \HomL(V^k, W)^{\oplus \delta_{0l}}
   \xrightarrow{\beta}
   \HomL(V^k, V^l)
\end{equation*}
where $\HomL(W, V^l)$ appears in case $k=0$ and $\HomL(V^k, W)$ in
case $l=0$. We also put $\widehat v_0 = 1$.
Then it is easy to show that $\Ker\alpha = 0$ unless $k = l \neq 0$
and $\Ker\alpha = \C \id$ for $k=l\neq 0$, and the similar statement
for $\Coker\beta$: Note that $\Ker\alpha$ is the space of
homomorphisms between $\zeta^\bullet$-stable modules, and $\Coker\beta$ is
its dual.
Then a standard argument, comparing
$\theta_{\Tilde\zeta^\bullet}$-values of $\Ker\xi$ and $\Ima\xi$ of a
homomorphism $\xi$ and using \lemref{lem:slope}, implies the
assertion.
We remark that a complex used often in \cite{Na-alg} (see
\eqref{eq:crystalcpx}) is an example of $\mathscr{C}_{kl}^\bullet$
where $V^k$ is $S_i$, a module with $\C$ on the vertex $i\in I$, and
$0$ on the other vertices and $B = 0$.

We construct a new graph with $\widehat I = \{1, \dots, r\}$ with the
associated Cartan matrix $\widehat{\bC} = (\widehat{c}_{kl})$ by
\(
  \widehat{c}_{kl} \defeq 2\delta_{kl} - \dim \Ker\beta/\Ima\alpha
\)
for the complex $\mathscr{C}^\bullet_{kl}$.
This is equal to the alternating sum of the dimensions of terms,
i.e.,
\(
  = {}^t \bv^k \bC \bv^l
\)
by the above discussion. Note $\widehat{a}_{kl} = \widehat{a}_{lk}$.
We also put
\begin{equation*}
   \hat V_k \defeq \C^{\widehat{v}_k}, \qquad
   \hat W_k \defeq \text{$\Ker\beta/\Ima\alpha$ \quad
     for \   $\mathscr{C}_{k0}^\bullet$},
\end{equation*}
and consider $\bM(\hat V,\hat W)$ def\/ined for the new graph
with the $\widehat I$-graded vector spaces $\hat V$, $\hat W$.
The stabilizer $\widehat G$ of $(B,a,b)$ is naturally isomorphic to
$\prod_{k\in\widehat I} \GL(\hat V_k)$. It acts on $\bM(\hat V,\hat
W)$. We have the moment map $\widehat\mu\colon\bM(\hat V,\hat W)\to
\HomL(\hat V,\hat V) \cong \operatorname{Lie}(\widehat G)^*$. We
consider the quotient
$\widehat\M_0(\Hat{V},\Hat{W}) = \widehat\mu^{-1}(0)\dslash \widehat
G$, where the stability parameter is $0$.
We also consider $\widehat\M_\zeta(\Hat{V},\Hat{W})$, where the
stability parameter $\zeta$ is considered as a $\R$-character
$\chi_{\zeta}$ of $\widehat G$ through the inclusion $\widehat
G\subset G$. There is a morphism
$\widehat\pi\colon \widehat\M_\zeta(\Hat{V},\Hat{W})
\to \widehat\M_0(\Hat{V},\Hat{W})$.
Then \cite[3.2.1]{Na-qaff} $\M_{\zeta^\bullet}(V,W)$, around $x =
[B,a,b]$ a point corresponding to $(B,a,b)$, is locally isomorphic to
a neighbourhood of $0$ in
$\widehat\M_0^{\text{norm}}(\Hat{V},\Hat{W})\times T$, where
$\widehat\M_0^{\text{norm}}(\Hat{V},\Hat{W})$ is as in \eqref{eq:fac},
$T$ is the product of the trivial factor
$\widehat{\bM}^{\text{el}}(\Hat{V},\Hat{W})$ (see \eqref{eq:trivial})
and $\Hat{W}_0 = \Ker\beta/\Ima\alpha$ for the complex
$\mathscr{C}^\bullet_{00}$ which is isomorphic to the tangent space
\linebreak[2] $T_{[(B^0,a^0,b^0)]}\M_{\zeta^\bullet}(V^0,W)$.
We also have a lifting of the local isomorphism to $\M_\zeta(V,W)$ and
$\widehat\M_{\zeta}(\Hat{V},\Hat{W})$, and hence the following
commutative diagram:
\begin{equation}\label{eq:diagram}
    \begin{matrix}
   \M_{\zeta}(V,W) & \supset & \pi^{-1}(U) & \xrightarrow[\cong]{\tilde\Phi}
    & \widehat\pi^{-1}(\Hat{U})\times U_0  & \subset &
    \widehat\M_\zeta^{\text{norm}}(\Hat{V},\Hat{W})\times T
\\
   && \vcenter{%
     \llap{$\scriptstyle{\pi}$}}
   \Big\downarrow && \Big\downarrow
   \vcenter{%
     \rlap{$\scriptstyle{\widehat\pi}\times\id$}}
   &&
\\
   \M_{\zeta^\bullet}(V,W) & \supset & U & \xrightarrow[\cong]{\Phi}
    & \Hat{U}\times U_0 & \subset &
    \widehat\M_0^{\text{norm}}(\Hat{V},\Hat{W})\times
    T
\\
   && \rotatebox{90}{$\in$} & & \rotatebox{90}{$\in$} &&
\\
   && x & \mapsto & 0 &&
  \end{matrix}
\end{equation}
Note that the factor $T$ is the tangent space to the stratum
$(\M_{\zeta^\bullet}(V,W))_{(\widehat G)}$ containing $x$.
(Remark: \cite[3.2.1]{Na-qaff} states this result for $\zeta^\bullet =
0$, but the same proof works. To see that the stability parameter for
$\widehat\M_{\zeta^\bullet}$ becomes $0$, we note also that the
restriction of the character $\chi_{\zeta^\bullet}$ to~$\widehat G$ is
trivial.)

The following formula will be useful:
\begin{equation}\label{eq:weights}
\begin{split}
   & \widehat{w}_k - \sum_{l\neq 0} \widehat{c}_{kl} \widehat{v}_l
   = \lsp{t}{\bv}^k (\bw - \bC \bv).
\end{split}
\end{equation}

\begin{Remark}\label{rem:affine}
  The Cartan matrix $\widehat{\bC} = (\widehat{c}_{kl})$ is given by
  the formula $\widehat{c}_{kl} = \lsp{t}{\bv^k} \bC \bv^l$, where $\bv^k =
  \dim V^k$. Suppose that the original graph is of af\/f\/ine type. Then
  $\bC$ is positive semidef\/inite and the kernel is spanned by
  $\delta$. Therefore $\bv^k$ with $\lsp{t}{\bv^k} \bC \bv^k = 0$ is a multiple
  of $\delta$, hence $\widehat{c}_{kl} = 0$ for any $l\in \widehat I$.
  This means that a connected component of the graph $(\widehat
  I,\widehat\bC)$ is either a graph of af\/f\/ine or f\/inite type (without
  edge loops), or a graph with a single vertex and a single edge loop
  (i.e., the Jordan quiver). If the original graph is of f\/inite type,
  then we only get a graph of f\/inite type as $\bC$ is positive
  def\/inite.
\end{Remark}

As an application, we obtain the following:

\begin{Proposition}[cf.~\protect{\cite[6.11]{Na-quiver},
\cite[10.4, 10.11]{Na-alg}}]\label{prop:semismall}
Let $\zeta$, $\zeta^\bullet$ be as above and consider
$\pi\colon\M_\zeta(V,\linebreak[3]W)\to \M_{\zeta^\bullet}(V,W)$.
  We further assume that $\zeta$ is in a chamber so that
  $\M_\zeta(V,W)$ is nonsingular.

\textup{(1)} Take a point $x$ in a stratum $(\M_{\zeta^\bullet})_{(\widehat
  G)}$. If $\pi^{-1}(x)\neq \varnothing$, then we have
\begin{equation*}
  \dim \pi^{-1}(x)
  \le \tfrac12
  \big( \dim \M_\zeta(V,W) - \dim (\M_{\zeta^\bullet})_{(\widehat G)}\big).
\end{equation*}

\textup{(2)} Replace the target of $\pi$ by the image, hence consider
$\pi\colon \M_{\zeta}(V,W)\to \pi(\M_\zeta(V,W))$. It is semismall
with respect to the stratification $\pi(\M_\zeta(V,W)) = \bigsqcup
(\M_{\zeta^\bullet})_{(\widehat G)}$ where $\widehat G$ runs over the
conjugacy classes of subgroups of $G_V$ such that
$(\M_{\zeta^\bullet})_{(\widehat G)}$ is contained in $\pi(\M_\zeta(V,W))$.
\end{Proposition}

\begin{proof}
  (1) This was proved under the assumption
  $\Mreg_{\zeta^\bullet}(V,W)\neq\varnothing$ in \cite[6.11]{Na-quiver},
  but its proof actually gives the above. The point of the proof was
  that the f\/iber $\widehat\pi^{-1}(0)$ is a subvariety in
  $\widehat\M_\zeta(\Hat{V},\Hat{W})$ (divided by the trivial factor
  \eqref{eq:trivial}), which is an af\/f\/ine algebraic manifold by
  \thmref{thm:affine}.

(2) By \eqref{eq:diagram} and the subsequent remark, for each stratum
$(\M_{\zeta^\bullet})_{(\widehat G)}$, the restriction of $\pi$ to the
preimage $\pi^{-1}((\M_{\zeta^\bullet})_{(\widehat G)})$ is a f\/iber bundle.
Moreover, by (1), there exists a stratum $(\M_{\zeta^\bullet})_{(\widehat
  G)}$ such that $\dim (\M_{\zeta^\bullet})_{(\widehat G)} = \dim
\M_\zeta(V,W)$. Therefore we have $\dim \pi(\M_\zeta(V,W)) = \dim
\M_\zeta(V,W)$. Hence (1) implies the assertion.
\end{proof}

\begin{Remark}\label{rem:relevant}
  In \cite[10.11]{Na-alg} we claimed a {\it stronger\/} statement
  that all strata are relevant in the statement (2) of
  \propref{prop:semismall} when the graph is of f\/inite type. This
  follows from the above observation (1) and the fact that
  $\widehat\pi^{-1}(0)$ is exactly half-dimensional in
  $\widehat{\M}_\zeta(\Hat{V},\Hat{W})$ if the graph does not contains
  edge loops \cite[5.8]{Na-quiver}.

  Now consider the case when the graph is of af\/f\/ine type. By the above
  observation (1), it is enough to show that the f\/iber
  $\widehat\pi^{-1}(0)$ in
  $\widehat{\M}_\zeta^\text{norm}(\Hat{V},\Hat{W})\to
  \widehat{\M}_0^\text{norm}(\Hat{V},\Hat{W})$ is exactly
  half-dimensional. It is clearly enough to prove this assertion for
  each connected component of the graph $(\widehat I,\widehat\bC)$.
  Since it is known when the component has no edge loops, it is enough
  to consider the case when the component is the Jordan quiver.
  Since the trivial factor
  $\widehat{\bM}^{\text{el}}(\Hat{V},\Hat{W})$ is $2$-dimensional
  (except for the trivial case $\Hat{V} = 0$), this means $\dim
  \widehat\pi^{-1}(0) = \dim \widehat{\M}_\zeta(\Hat{V},\Hat{W})/2 -
  1$. This is known. (See \cite[Exercise~5.15]{Lecture} and the
  references therein.)
  Therefore all strata are relevant also in the af\/f\/ine case.

  For general quivers, the author does not know whether the same
  result holds or not.
\end{Remark}

\subsection{Example: Levi factors of parabolic
  subalgebras}\label{subsec:Levi}

We give an example of a face, which will be related to the restriction
to the Levi factor of a parabolic subalgebra.

Consider the chamber $\mathcal C = \{ \zeta\in \R^I \mid
\text{$\zeta_i > 0$ for all $i\in I$}\}$. This corresponds to the
stability condition $\zeta$ in Example~\ref{ex:stndcham}(1). A face
$F$ contained in the closure $\overline{\mathcal C}$ is of a form
\begin{equation*}
    F = \{ \zeta^\bullet\in \R^I \mid \text{
     $\zeta_i^\bullet = 0$ (resp.\ $> 0$) for
    $i \in I^0$ (resp.\ $I^+$)}\}
\end{equation*}
for a disjoint decomposition $I = I^0 \sqcup I^+$. We allow the cases
$I^0 = \varnothing$, i.e., $F = \mathcal C$, $\zeta^\bullet = \zeta$ and
$I^0=I$, i.e., $\zeta^\bullet = 0$.
We have
\begin{equation*}
   \M_\zeta(\bv,\bw)\xrightarrow{\pi_{\zeta^\bullet,\zeta}}
   \M_{\zeta^\bullet}(\bv,\bw)\xrightarrow{\pi_{0,\zeta^\bullet}}
   \M_0(\bv,\bw).
\end{equation*}

Let us describe strata of $\M_{\zeta^\bullet}(\bv,\bw)$ in this
example.

\begin{Proposition}\label{prop:Levi}
  \textup{(1)} The strata are of the form
\begin{equation}\label{eq:Levistratum}
   \Mreg_{\zeta^\bullet}(\bv^0,\bw) \times (\M_0)_{(\widehat G)},
\end{equation}
where $\bv - \bv^0$ is supported on $I^0$, and
$(\M_0)_{(\widehat G)}$ is a stratum of $\M_0(\bv-\bv^0,0)$ for
the subgraph $I^0$, extended to the whole graph by $0$.

\textup{(2)} If $\Mreg_{\zeta^\bullet}(\bv^0,\bw)\neq\varnothing$, we
have
\(
   \langle h_i, \bw - \bv^0\rangle \ge 0
\)
for $i\in I^0$.
\end{Proposition}

We say a weight $\lambda$ is {\it $I^0$-dominant\/} if
\(
\langle h_i, \lambda\rangle \ge 0
\)
for all $i\in I^0$.

\begin{proof}
  (1) A closed point in $\M_{\zeta^\bullet}(\bv,\bw)$ is represented
  by a $\zeta^\bullet$-polystable module. Let
  $B\in\Mreg_{\zeta^\bullet}(\bv',0)$ be its component with $W = 0$.
  As we have the constraint $\zeta^\bullet\cdot \bv' = 0$, our choice
  of $\zeta^\bullet$ implies that $\bv'$ is supported on $I^0$. Thus
  $\zeta^\bullet$-stable modules are $0$-stable (i.e., simple) modules
  for the subgraph $I^0$, extended by $0$.

(2) is implicitly proved in \subsecref{subsec:local}. Let us make it
more apparent. Consider the complex $\mathscr{C}^{\bullet}_{kl}$ for a point
$[B,a,b]\in \Mreg_{\zeta^\bullet}(\bv^0,\bw)$ and $S_i$, a module with
$\C$ on the vertex $i\in I^0$, and $0$ on the other vertices and $B =
0$. Explicitly it is given by (see \cite[(4.2)]{Na-alg}):
\begin{gather}\label{eq:crystalcpx}
 V_i
  \xrightarrow{\alpha_i}
  \bigoplus_{h: \vin{h} = i} V_{\vout{h}} \oplus W_i
  \xrightarrow{\beta_i} V_i,
\\
   \alpha_i\defeq
    \bigoplus_{\vin{h}=i} B_{\overline h}\oplus b_i, \qquad
   \beta_i\defeq
    \begin{bmatrix}
    \displaystyle\sum_{\vin{h}=i} \varepsilon(h) B_h & a_i
    \end{bmatrix}.\tag*{\qed}
\end{gather}\renewcommand{\qed}{}
\end{proof}

\begin{Lemma}[\protect{cf.\ \cite[4.7]{Na-alg}}]\label{lem:betasurj}
  Let $\zeta^\bullet$ as above. Fix $i\in I^0$. If $(B,a,b)$ is
  $\zeta^\bullet$-stable, then either of the following hold\textup:
  \begin{aenume}\itemsep=0pt
  \item $W = 0$ \textup(hence $a = b = 0$\textup), and $B$ is the simple
    module $S_i$.

  \item The map $\alpha_i$ is injective and $\beta_i$ is
    surjective.
  \end{aenume}
\end{Lemma}

Since $W = 0$ is excluded in (b), this gives the proof of the proposition.

\begin{proof}
  This is a standard argument. The kernel of $\alpha_i$ gives a
  homomorphism from $S_i$ to $(B,a,b)$. Since both are
  $\zeta^\bullet$-stable and have the same
  $\theta_{\widehat\zeta^\bullet}$ (both are $0$), the homomorphism is
  \mbox{either}~$0$ or an isomorphism. Similarly the cokernel of $\beta_i$ is
  the dual of the space of homomorphisms from $(B,a,b)$ to $S_i$. Thus
  we have the same assertion.
\end{proof}

Let us consider the restriction of $\pi_{0,\zeta^\bullet}$ to the
closure of the stratum $\Mreg_{\zeta^\bullet}(\bv^0,\bw) \times
(\M_0)_{(\widehat G)}$ (in $\M_{\zeta^\bullet}(\bv,\bw)$). It f\/its in
the commutative diagram
\begin{equation}\label{eq:basechange}
\begin{CD}
   \overline{\Mreg_{\zeta^\bullet}(\bv^0,\bw)
     \times (\M_0)_{(\widehat G)}} @<{\kappa}<<
    \M_{\zeta^\bullet}(\bv^0,\bw) \times \overline{(\M_0)_{(\widehat G)}}
\\
   @V{\pi_{0,\zeta^\bullet}}VV @VV{\pi_{0,\zeta^\bullet}\times\id}V
\\
   \M_0(\bv,\bw) @<{\kappa'}<<
   \M_0(\bv^0,\bw) \times \overline{(\M_0)_{(\widehat G)}}
\end{CD}
\end{equation}
where the horizontal arrows are given by $([B,a,b], [B'])\mapsto
[(B,a,b)\oplus B']$. They are f\/inite morphisms. The upper horizontal
arrow $\kappa$ is an isomorphism on the preimage of
$\Mreg_{\zeta^\bullet}(\bv^0,\bw) \times (\M_0)_{(\widehat G)}$. This
means that a study of a stratum can be reduced to the special case
when it is of the form $\Mreg_{\zeta^\bullet}(\bv^0,\bw)$.

For a general graph, we can use \thmref{thm:CB} in principle, but to
apply this result, we need to know all roots, and it is not so easy in
general. Therefore we restrict ourselves to the case when the graph is
of af\/f\/ine type from now until the end of this section.

\begin{Proposition}[cf.~\protect{\cite[10.5, 10.8]{Na-alg}}]
\label{prop:affinestrata}
Suppose the graph is of affine type and $\bw\neq 0$.

\textup{(1)} Suppose $I^0 = I$ and $\lsp{t}{\delta}\bw = 1$.
Then $\Mreg_{\zeta^\bullet}(\bv,\bw) = \varnothing$ unless $\bv = 0$.

\textup{(2)} Suppose otherwise.
Then $\Mreg_{\zeta^\bullet}(\bv,\bw) \neq \varnothing$ if and only if
$\bw-\bv$ is an $I^0$-dominant weight
of the irreducible integrable highest weight module $V(\bw)$.
\end{Proposition}

\begin{proof}
(1) We have $\zeta^\bullet = 0$ in this case.

Let us f\/irst give a proof based on the interpretation of
$\Mreg_0(\bv,\bw)$ as moduli of vector bundles.
Let $\Gamma$ be the f\/inite subgroup of $\SL_2(\C)$ corresponding to
the graph via the McKay correspondence.
Then $\Mreg_0(\bv,\bw)$ parametrizes framed $\Gamma$-equivariant
vector bundles over $\proj^2$~\cite{Lecture}. Since $\bw$ is
$1$-dimensional, it parametrizes line bundles. As the existence of the
framing implies $c_1 = 0$, the line bundles are trivial. In
particular, we have $\bv = 0$.

Let us give another proof by using \thmref{thm:CB}. The weights $\bw -
\bv$ of the level $1$ representation $V(\bw)$ are of the form
\(
  \bw - \bv^0 - n\delta,
\)
where $\bw - \bv^0\in W\cdot \bw$, $n\in\Z_{\ge 0}$.
(Here $W$ is the af\/f\/ine Weyl group.)
Then
\[
   \lsp{t}{\bv}\left(\bw - \tfrac12 \bC \bv\right)
   = \lsp{t}{\bv^0}\left(\bw - \tfrac12 \bC \bv^0\right) + n
   = \lsp{t}{\bv^0}\left(\bw - \tfrac12 \bC \bv^0\right) + n p(\delta).
\]
This violates the inequality in \thmref{thm:CB} unless $n=0$. But we
also know that $\bw - \bv^0$ is dominant by \lemref{lem:betasurj}.
Hence $\bw - \bv^0 = \bw$, i.e., $\bv^0 = 0$.

(2) This is proved in \cite[10.5, 10.8]{Na-alg} in the special case
$I^0 = I$. The argument works in general. We will give another
argument based on \thmref{thm:CB} in a similar situation later
(\propref{prop:affinemodulistrata}), we omit the detail here.
\end{proof}

\section{Instantons on ALE spaces}\label{sec:instanton}

Quiver varieties were originally introduced by generalizing the ADHM
description of instantons on ALE spaces by Kronheimer and the author
\cite{KN}. In this section we go back to the original description and
explain partial resolutions in terms of instantons (or sheaves) on
(possibly singular) ALE spaces.

\subsection{ALE spaces}\label{subsec:ALE}

We review Kronheimer's construction~\cite{Kr} of ALE spaces brief\/ly in
our terminology.

We consider the untwisted af\/f\/ine Lie algebra of type $ADE$.
Let $0\in I$ be the vertex corresponding to the simple root, which is
the negative of the highest weight root of the corresponding simple
Lie algebra. Let $I_0 \defeq I\setminus \{0\}$. Let $\delta$ be the
vector in the kernel of the af\/f\/ine Cartan matrix whose $0$-component
is equal to $1$. Such a vector is uniquely determined.
Let $G_\delta$ be the complex Lie group corresponding to $\delta$ as
in \subsecref{subsec:quiver}. Choose a parameter $\zeta^\circ \in
\R^I$ from the level $0$ hyperplane $\{ \zeta\in\R^I \mid
\zeta\cdot \delta = 0\}$.
Let us denote the corresponding quiver variety~$\M_{\zeta^\circ}(\delta,0)$
for the parameter $\zeta^\circ$ by $X_{\zeta^\circ}$.
This space is called an {\it ALE space} in the literature.
We have a morphism
\begin{equation}\label{eq:minimalresol}
   \pi_{0,\zeta^\circ}\colon X_{\zeta^\circ} \to X_0
\end{equation}
from the construction in \subsecref{subsec:partial}.

If we take $\zeta^\circ$ from an open face $F$ in the level $0$ hyperplane,
i.e., it is not contained in any real root hyperplane $D_\theta$, then
$X_{\zeta^\circ}$ is nonsingular by \remref{rem:W=0reg}. Kronheimer~\cite{Kr}
showed
\begin{aenume}\itemsep=0pt
\item $X_0$ is isomorphic to $\C^2/\Gamma$, where $\Gamma$ is the
  f\/inite subgroup of $\SL_2(\C)$ associated to the af\/f\/ine Dynkin
  graph,

\item \eqref{eq:minimalresol} is the minimal resolution of
  $\C^2/\Gamma$, if $\zeta^\circ$ is taken as above.
\end{aenume}

For a later purpose we take a specif\/ic face $F^\circ$ with the above
property and a parameter $\zeta^\circ$ as
\begin{equation*}
  F^\circ \defeq \{ \zeta^\circ\in \R^I \mid
  \text{$\zeta^\circ\cdot \delta = 0$ and $\zeta^\circ_i > 0$ for $i\neq 0$}\}.
\end{equation*}
By the Weyl group action, any open face in the level $0$ hyperplane is
mapped to $F^\circ$.

A face $F^\bullet$ in the closure of $F^\circ$ is of the form
\begin{equation*}
    F^\bullet = \{ \zeta^\bullet\in \R^I \mid \text{
     $\zeta^\bullet\cdot \delta = 0$,
     $\zeta^\bullet_i = 0$ (resp.\ $> 0$) for
    $i \in I_0^0$ (resp.\ $I_0^+$)}\}
\end{equation*}
for a disjoint decomposition $I_0 = I_0^0 \sqcup I_0^+$.
We allow the cases $I_0^0 = \varnothing$, i.e., $\zeta^\bullet =
\zeta^\circ$ or $I_0^0=I_0$, i.e., $\zeta^\bullet = 0$.
From the construction in \subsecref{subsec:partial} we have
\begin{equation*}
  X_{\zeta^\circ} \xrightarrow{\pi_{\zeta^\bullet,\zeta^\circ}}
  X_{\zeta^\bullet}\xrightarrow{\pi_{0,\zeta^\bullet}}
  X_0 \cong \C^2/\Gamma.
\end{equation*}
Kronheimer also showed that (see \cite[Lemma~3.3]{Kr})
\begin{aenume}\itemsep=0pt
\setcounter{enumi}{2}
\item $X_{\zeta^\bullet}$ is a partial resolution of $\C^2/\Gamma$
  having singularities of type $\C^2/\Gamma'$ for some dif\/ferent
  $\Gamma'\subset\SL_2(\C)$.
\end{aenume}

In fact, we can describe singularities of $X_{\zeta^\bullet}$
explicitly. Recall that the exceptional set of the minimal resolution~\eqref{eq:minimalresol} consists of an union of the complex projective
line. By \cite{Na-proc} each irreducible component $C_i$ naturally
corresponds to a nonzero vertex $i\in I_0$ as follows: $C_i$ consists
of data $B$ having a quotient isomorphic to $S_i$, a module with $\C$
on the vertex $i$, and $0$ on the other vertices.
It is clear that $B$ is $\zeta^\bullet$-semistable, but not
$\zeta^\bullet$-stable.
Then the curves $C_i$ with $i\in I_0^0$ are contracted under the
morphism
\(
   X_{\zeta^\circ}
   \xrightarrow{\pi_{\zeta^\bullet,\zeta^\circ}} X_{\zeta^\bullet}.
\)

In order to show that $\pi_{\zeta^\bullet,\zeta^\circ}$ does not
contract further curves, we give the following more precise
classif\/ication:

\begin{Lemma}\label{lem:bulletstable}
  A $\zeta^\bullet$-stable point $B\in\mu^{-1}(0)\subset\bM(\bv,0)$
  satisfying the normalization condition $\zeta^\bullet\cdot\bv = 0$
  is of one of the following three forms\textup:
  \begin{aenume}\itemsep=0pt
  \item a point in $X_{\zeta^\circ}\setminus \bigcup_{i\in I_0^0}
    C_i$,
  \item $\bv$ is a coordinate vector $\be_i$ for an $i\in
    I_0^0$, or
  \item $\bv$ is $\delta - \alpha_h$, where $\alpha_h$ is the highest
    root vector of a connected component of the sub-Dynkin diagram
    $I_0^0$.
  \end{aenume}
  In cases \textup{(b)}, \textup{(c)} the corresponding
  $\zeta^\bullet$-stable point is unique up to isomorphisms.
\end{Lemma}

Let $\mathcal C$ be the set of components of the sub-Dynkin diagram
$I_0^0$. Let $B_c$ denote (the isomorphism class of) a
$\zeta^\bullet$-stable representation corresponding to a component
$c\in \mathcal C$ in (c). Then the $S$-equivalence class $x_c$ of
\(
   B_c\oplus \bigoplus_{i\in c} S_i^{\oplus (\alpha_h)_i}
\)
def\/ines a point in $X_{\zeta^\bullet}$, where $(\alpha_h)_i$ is the
$i^{\mathrm{th}}$-entry of $\alpha_h$. If $c$ is a Dynkin diagram of
type $ADE$, $X_{\zeta^\bullet}$ has a singularity of the corresponding type
around $x_c$. This is an example of the description in
\subsecref{subsec:local}. The new graph has vertices
$\widehat I = c\sqcup\{ 0\}$ with the Cartan matrix
\(
   \widehat{c}_{ij} = c_{ij}
\)
if $i$, $j\in c$,
\(
   \widehat{c}_{0i} = - \lsp{t}{(\delta - \alpha_h)}\bC \be_i
   = \lsp{t}{\alpha_h} \bC \be_i
\)
if $i\in c$,
and
$\widehat{c}_{00} = 2$. This is a graph of af\/f\/ine type.

\begin{proof}[Proof of \lemref{lem:bulletstable}]
  Consider the criterion in \thmref{thm:CB}. Since our graph is of
  af\/f\/ine type, we have $p(x) = 1$ if $x$ is an imaginary root (i.e.,
  $x = m\delta$ with $m\ge 1$) and $p(x) = 0$ if $x$ is a real root.
  Therefore if $\bv$ is imaginary root, the condition is equivalent to
  $\bv = \delta$. This is the case (a). If $\bv$ is a real root, then
  the condition is that $\bv$ cannot have a nontrivial decomposition
  $\bv = \sum_t \beta^{(t)}$ with $\zeta^\bullet\cdot \beta^{(t)} = 0$. Then
  it is clear that we have either (b) or (c).

Next show that a $\zeta^\bullet$-stable point is unique in case (b),
(c). It is clear in the case (b). Consider the case (c). We apply an
argument used by Mukai for the case of rigid sheaves on a $K3$
surface:
Let $B$, $B'$ be $\zeta^\bullet$-stable points with the underlying
$I$-graded vector space $V$ of dimension vector $\delta-\alpha_h$. We
consider the complex
\begin{gather*}
  \HomL(V,V) \xrightarrow{\alpha} \HomE(V,V) \xrightarrow{\beta}
  \HomL(V,V),
\\
  \alpha(\xi) = B\xi - \xi B', \qquad
  \beta(C) = \varepsilon (B C - C B').
\end{gather*}
Since the alternating sum of dimensions is $2$, we cannot have both
$\alpha$ is injective and $\beta$ is surjective. If $\alpha$ is not
injective, we have a nonzero homomorphism $\xi$ from the module $B$ to
$B'$. Then considering submodules $\Ker\xi\subset B$ and
$\Ima\xi\subset B'$, the $\zeta^\bullet$-stability of $B$ and $B'$
imply $\Ker\xi = 0$ and $\Ima\xi = B'$. Therefore $B$ and $B'$ are
isomorphic. If $\beta$ is not surjective, we consider
$\beta^*\in \HomL(V^*,V^*) \cong \HomL(V,V)$ instead, we get the same
assertion.

Finally consider the case $\bv = \delta$.
By \lemref{lem:face}(3b) and the subsequent remark $B$ is
$\zeta^\circ$-stable. So $B$ corresponds to a point $[B]$ in
$X_{\zeta^\circ}$. Therefore it only remains to determine which point
in $X_{\zeta^\circ}$ is $\zeta^\bullet$-stable.

If $B$ is not $\zeta^\bullet$-stable, we consider a
Jordan--H\"older f\/iltration to f\/ind a proper quotient which is
$\zeta^\bullet$-stable. By the above discussion, the quotient must
have the dimension vector $\be_i$ with $i\in I_0^0$ or
$\delta-\alpha_h$. But in the latter case, $B$ contains a submodule
with the dimension vector $\be_i$ with $i\in I_0^0$. This is
impossible for the $\zeta^\circ$-stability, as $\zeta^\circ \cdot
\be_i > 0$.
Thus $B$ has the quotient with the dimension vector $\be_i$ with
$i\in I_0^0$. In particular, $\beta$ above is not surjective, and
hence $B$ is contained in the exceptional component $C_i$ by
\cite[Proof of 5.10]{Na-proc}.
On the other hand, if $B$ is in $C_i$, $\beta$ is not surjective,
and it has a quotient with $\dim = \be_i$. It is not
$\zeta^\bullet$-stable.
Therefore $\zeta^\bullet$-stable points are precisely in
$X_{\zeta^\circ}\setminus \bigcup_{i\in I_0^0} C_i$.
\end{proof}

\begin{Remark}
  It is not so dif\/f\/icult to prove the criterion in \thmref{thm:CB}
  directly without referring to a general result \cite{CB} in this
  particular example. The detail is left for a reader as an exercise.
\end{Remark}

\subsection{Sheaves and instantons on ALE spaces}\label{subsec:sheaf}

Let $\zeta^\circ$, $\zeta^\bullet$ be as in the previous subsection.
We consider the corresponding quiver varieties
$\M_{\zeta^\circ}(\bv,\bw)$, $\M_{\zeta^\bullet}(\bv,\bw)$ for
$\bw\neq 0$.
By the main result of \cite{KN}, the former space
$\M_{\zeta^\circ}(\bv,\bw)$ is the Uhlenbeck (partial)
compactif\/ication of framed instantons on $X_{\zeta^\circ}$. Since we
do not need to recall what {\it instantons\/} or their {\it framing\/}
mean in this paper, we refer the def\/initions to the original paper~\cite{KN}.
Rather we use a dif\/ferent (but closely related) description in
\cite{Na-ADHM}. We f\/irst compactify $X_{\zeta^\circ}$ to an orbifold
$\overline{X}_{\zeta^\circ}$ by adding $\linf \defeq \proj^1/\Gamma$.
It is obtained by resolving the (isolated) singularity at $0$ of
$\proj^2/\Gamma$. Then $\Mreg_{\zeta^\circ}(\bv,\bw)$ is a f\/ine moduli
space of framed holomorphic orbifold vector bundles $(E,\Phi)$, where
$\Phi$ is an isomorphism between $E|_{\linf}$ and
$(\rho\otimes\shfO_{\proj^1})/\Gamma$ for a~f\/ixed representation
$\rho$ of $\Gamma$. Here $\rho$ corresponds to the vector $\bw$ via
the McKay correspondence: $w_i$ is the multiplicity of the irreducible
representation $\rho_i$ in $\rho$. And $\bv$ corresponds to Chern
classes of $E$. Since the explicit formula is not relevant here, we
refer to \cite[(1.9)]{Na-ADHM} for an interested reader.

Let us consider $\M_{\zeta^\circ}(\bv,\bw)$ which is a partial
compactif\/ication of $\Mreg_{\zeta^\circ}(\bv,\bw)$. Its closed point
is represented by an $S$-equivalence class of
$\zeta^\circ$-semistable points, and hence a direct sum of
$\zeta^\circ$-stable point all having $\theta_{\zeta^\circ} = 0$.
By \cite[Proposition~9.2(ii)]{KN} (or \cite[4.3]{Na-ADHM} for a detailed
argument) such a point is of a form
\(
    (B,a,b) = (B',a',b')\oplus x_1 \oplus x_2 \oplus \cdots \oplus x_k,
\)
where $(B',a',b')\in\Mreg_{\zeta^\circ}(\bv^0,\bw)$ with $\bw\neq 0$ and
$x_i\in\M_{\zeta^\circ}(\delta,0)$ corresponds to a point in the ALE
space $X_{\zeta^\circ}$.
It is also a special case of \lemref{lem:bulletstable} with $I_0^0 =
\varnothing$.
Therefore $\M_{\zeta^\circ}(\bv,\bw)$ is written as
\begin{equation}\label{eq:Uhzetacirc}
  \M_{\zeta^\circ}(\bv,\bw)
= \bigsqcup_{k\ge 0} \Mreg_{\zeta^\circ}(\bv-k\delta,\bw)
  \times S^k X_{\zeta^\circ},
\end{equation}
where $S^k X_{\zeta^\circ}$ is the $k^{\mathrm{th}}$ symmetric product
of $X_{\zeta^\circ}$.
As $\Mreg_{\zeta^\circ}(\bv-k\delta,\bw)$ is the framed moduli space
of orbifold holomorphic vector bundles on $\overline{X}_{\zeta^\circ}$
with a smaller second Chern number, the above description means that
$\M_{\zeta^\circ}(\bv,\bw)$ is the Uhlenbeck (partial)
compactif\/ication of $\Mreg_{\zeta^\circ}(\bv,\bw)$.

Let us apply the results in \subsecref{subsec:local} in this
situation. In order to give a decomposition~\eqref{eq:decomp} we need
to introduce a f\/iner stratif\/ication for the symmetric power:
\begin{equation*}
   S^k X_{\zeta^\circ} = \bigsqcup_{|\lambda| = k} S^k_\lambda X_{\zeta^\circ},
\end{equation*}
where $\lambda = (\lambda_1 \ge \lambda_2\ge \cdots\ge \lambda_l > 0)$
is a partition of $k$ and $S^k_\lambda X_{\zeta^\circ}$ consists of
those conf\/igurations of form $\sum_{n=1}^l \lambda_n [x_n]$ for
$x_n$ distinct. If we take a module from a stratum
$\Mreg_{\zeta^\circ}(\bv-k\delta,\bw) \times S^k_\lambda
X_{\zeta^\circ}$, the corresponding graph is
\(
   \widehat I = \{ 1,\dots, l = l(\lambda)\},
\)
with
\(
   \widehat c_{kk'} = 0.
\)
This is the disjoint union of $l = l(\lambda)$ copies of the edge loop.
We have $\Hat{W}_k = \C^{\lsp{t}{\bw} \delta}$. Therefore
$\M_{\zeta^\circ}(\bv,\bw)$ is locally isomorphic to
\(
   \C^{\dim \Mreg_{\zeta^\circ}(\bv-k\delta,\bw)}\times
   \prod_{k=1}^l M_0(\C^{\lambda_k}, \C^{\lsp{t}{\bw} \delta}),
\)
where $M_0(\ ,\ )$ is the quiver variety associated with the Jordan
quiver, a single vertex and a single edge loop, which is known to be
the Uhlenbeck partial compactif\/ication for $\C^2$~\cite{Lecture}.

Next consider the stability parameter $\zeta^\bullet$. By \lemref{lem:bulletstable} we have
\begin{equation}\label{eq:Uhbullet}
  \M_{\zeta^\bullet}(\bv,\bw)
= \bigsqcup \Mreg_{\zeta^\bullet}(\bv^0,\bw)
  \times S^{|\lambda|}_\lambda
  \Big[X_{\zeta^\circ}\setminus \bigcup_{i\in I_0^0} C_i \Big]
  \times \left\{
    \bigoplus_{c\in\mathcal C} \left(
      B_c^{\oplus n_c} \oplus
    \bigoplus_{i\in c}  S_i^{\oplus m_i} \right)
    \right\},
\end{equation}
where $c$ runs over the set $\mathcal C$ of connected components of
the sub-Dynkin diagram $I_0^0$ with the highest root $\alpha^c_h$, and
$B_c$ is a $\zeta^\bullet$-stable module with the dimension vector
$\delta - \alpha^c_h$, which is unique up to isomorphisms.
(See \lemref{lem:bulletstable}.)
We have an obvious constraint
\(
   \bv = \bv^0 + |\lambda| \delta + \sum_{i\in I_0^0} m_i \be_i
   + \sum_c n_c (\delta - \alpha^c_h).
\)

This space can be considered as the partial Uhlenbeck compactif\/ication
of the framed moduli space of instantons on the orbifold
$X_{\zeta^\bullet}$ by a simple generalization of the result in
\cite{KN}. The factor $\Mreg_{\zeta^\bullet}(\bv^0,\bw)$ parametrizes
holomorphic orbifold vector bundles (i.e., ref\/lexive sheaves), and
$S^k \big[X_{\zeta^\circ}\setminus \bigcup_{i\in I_0^0} C_i\big]$ are
unordered points with multiplicities. These are usual factors
appearing in the Uhlenbeck compactif\/ication. The second factor is the
length of a $0$-dimensional sheaf $Q$ supported on
$X_{\zeta^\circ}\setminus \bigcup_{i\in I_0^0} C_i$.
The last factor
\(
    B_c^{\oplus n_c} \oplus
    \bigoplus_{i\in c}  S_i^{\oplus m_i}
\)
is new, and corresponds to a~representation of the local fundamental
group around the singular point $x_c$, which is a f\/inite subgroup
$\Gamma_c$ corresponding to the sub-Dynkin diagram $c$ via the McKay
correspondence,
i.e., $i\in c$ corresponds to a nontrivial irreducible representation,
and $B_c$ corresponds to the trivial representation, which usually
corresponds to the $0$-vertex.
It corresponds to a $0$-dimensional sheaf $Q$
supported at $x_c$, but we encode not only its length, but also the
$\Gamma_c$-module structure.

From the argument in \cite{Na-ADHM} (after modif\/ied to the case of
$X_{\zeta^\bullet}$), we see that the morphism
$\pi_{\zeta^\bullet,\zeta^\circ}$ is given in each stratum of
\eqref{eq:Uhzetacirc} by
\begin{gather*}
  \Mreg_{\zeta^\circ}(\bv^0,\bw)\ni
   (E,\Phi) \mapsto
   ((\pi_{\zeta^\bullet,\zeta^\circ})_*(E)^{\vee\vee},\Phi)
   \in\Mreg_{\zeta^\bullet}((\bv^0)',\bw),
\\
 S^k X_{\zeta^\bullet} \ni \sum \lambda_n [x_n]
   \mapsto
   \sum_{\xi_n\notin \bigcup C_i} \lambda_n [\pi_{\zeta^\bullet,\zeta^\circ}(x_n)],
\end{gather*}
and the remaining factor $\bigoplus_{c\in\mathcal C} \left(
      B_c^{\oplus n_c} \oplus
    \bigoplus_{i\in c}  S_i^{\oplus m_i} \right)$ is determined so
that the map preserve the dimension vector $\bv$.

Let us apply the local description in \subsecref{subsec:local} in this
situation. Take a point $x$ from the stratum of the above form. Then
the graph $\widehat{I}$ is the disjoint union of $l = l(\lambda)$
copies of Jordan quivers and the af\/f\/ine graphs corresponding to the
connected components $c$ of $\mathcal C$ (i.e., we add the $0$-vertex
to $c$).  The dimension vectors $\widehat{\bv}$, $\widehat{\bw}$ have
$\lambda_i$, $\lsp{t}{\bw}\delta$ in the components for Jordan
quivers. The components for the af\/f\/ine graph attached to $c$ are
\begin{equation}\label{eq:hatdim}
   \widehat{\bv} : \ (n_c, (m_i)_{i\in c}), \qquad
   \widehat{\bw} : \ \big( \lsp{t}{(\delta-\alpha^c_h)}(\bw - \bC\bv^0),
   \left.(\bw - \bC\bv^0)\right|_c\big),
\end{equation}
where the f\/irst components are the entries for the $0$-vertex, and
$\left.(\ )\right|_c$ means taking the components in $c$.
In particular the entries of $\widehat\bw$ must be nonnegative. This
follows from the same argument as in the proof of
\propref{prop:Levi}(2) and \lemref{lem:betasurj}. We consider the
complex $\mathscr{C}^\bullet_{kl}$ in \subsecref{subsec:local} for
$[B,a,b]\in \Mreg_{\zeta^\bullet}(\bv^0, \bw)$ and $B_c$, $S_i$.

Following \cite{Na-ADHM} we take the chamber $\mathcal C$ containing
$\zeta$ with
\begin{equation}\label{eq:zetasheaf}
   \zeta_i = \zeta_i^\circ \quad \text{for $i\neq 0$}, \qquad
   \text{$\zeta\cdot\delta$ is a suf\/f\/iciently small negative number}.
\end{equation}
Then $\mathcal C$ contains $\zeta^\circ$ in its closure. Therefore we
have a morphism $\pi_{\zeta^\circ,\zeta}\colon \M_{\zeta}\to
\M_{\zeta^\circ}$. The main result of \cite{Na-ADHM} says that
$\M_{\zeta}$ is a f\/ine moduli space of framed torsion free sheaves
$(E,\Phi)$, where $\Phi$ is as above. (In particular, $E$ is locally
free on $\linf$.) By its proof the morphism $\pi_{\zeta^\circ,\zeta}$
is given by the association
\[
   (E,\Phi) \mapsto ((E^{\vee\vee},\Phi), \operatorname{len}(E^{\vee\vee}/E)),
\]
where $E^{\vee\vee}$ is the double dual of $E$, which is locally free
as $X_{\zeta^\circ}$ is a nonsingular surface, and
$\operatorname{len}(E^{\vee\vee}/E)$ is the length of $E^{\vee\vee}/E$
considered as a conf\/iguration of unordered points in $X_{\zeta^\circ}$
counted with multiplicities.
This $\M_{\zeta}$ is called the {\it Gieseker partial
  compactification\/} of the framed moduli $\Mreg_{\zeta^\circ}$ of
locally free sheaves on the ALE space $X_{\zeta^\circ}$.

Strictly speaking we cannot take $\zeta$ independently from $\bv$.
When $\bv$ becomes larger, we need to take $\zeta$ closer and closer
to $\zeta^\circ$. In particular, we cannot specify $\zeta$ when we
move $\bv$ (as we will do in \secref{sec:convolution}). Since it is
cumbersome to use dif\/ferent notation for $\zeta$ for each $\bv$, we
simply denote all parameters by $\zeta$.

In summary we have four spaces and morphisms between them:
\begin{equation*}
   \M_{\zeta} \xrightarrow{\pi_{\zeta^\circ,\zeta}}
   \M_{\zeta^\circ}\xrightarrow{\pi_{\zeta^\bullet,\zeta^\circ}}
   \M_{\zeta^\bullet}\xrightarrow{\pi_{0,\zeta^\bullet}}
   \M_0,
\end{equation*}
where $\M_{\zeta}$ is the Gieseker partial compactif\/ication on
$X_{\zeta^\circ}$, and $\M_{\zeta^\circ}$, $\M_{\zeta^\bullet}$,
$\M_0$ are the Uhlenbeck partial compactif\/ication on
$X_{\zeta^\circ}$, $X_{\zeta^\bullet}$, $X_0 = \C^2/\Gamma$ respectively.

\begin{Proposition}[cf.~\protect{\cite[10.5, 10.8]{Na-alg}}]
\label{prop:affinemodulistrata}
\textup{(1)} Suppose $\lsp{t}{\delta}\bw = 1$.
Then $\Mreg_{\zeta^\bullet}(\bv,\bw) \neq\varnothing$ if and only if
$\bw - \bv$ is in the Weyl group orbit of $\bw$ and satisfies
\[
  \lsp{t}{(\delta-\alpha^c_h)}(\bw - \bC\bv)\ge 0
  \quad\text{for $c\in\mathcal C$}, \qquad
  \lsp{t}{\be_i}(\bw - \bC\bv)\ge 0 \quad\text{for $i\in I_0^0$}.
\]

\textup{(2)} Suppose $\lsp{t}{\delta}\bw \ge 2$.
Then $\Mreg_{\zeta^\bullet}(\bv,\bw) \neq \varnothing$ if and only if
$\bw-\bv$ is a weight of the irreducible integrable highest weight
module $V(\bw)$ and the above inequalities hold.
\end{Proposition}

\begin{proof}
One can give a proof along arguments in \cite[10.5, 10.8]{Na-alg}, but
we use \thmref{thm:CB} here.

If $\Mreg_{\zeta^\bullet}(\bv,\bw)\neq \varnothing$, then $\bw - \bv$ is
  a weight of the irreducible integrable highest weight module by
\thmref{thm:CB}(2).
Moreover two inequalities hold, as we have
  explained why $\widehat\bw$ has nonnegative entries.
  If $\lsp{t}{\delta}\bw = 1$, the second proof of
  \propref{prop:affinestrata} shows $\bw-\bv\in W\cdot\bw$ is also
  necessary.

  For the converse, we check the criterion in \thmref{thm:CB}. From
  our choice of $\zeta^\bullet$, a positive root $\beta^{(t)}$ with
  $\zeta^\bullet\cdot \beta^{(t)} = 0$ is one of the following:
  \begin{aenume}\itemsep=0pt
    \item $m\delta$ for $m > 0$,
    \item $m\delta + \alpha$ for $m\ge 0$, $\alpha\in
      \Delta(I_0^0)_+$, or
    \item $m\delta - \alpha$ for $m > 0$, $\alpha\in \Delta(I_0^0)_+$,
  \end{aenume}
  where $\Delta(I_0^0)_+$ is the set of positive roots of the subroot
  system $I_0^0$. We have
  \begin{equation*}
    \lsp{t}{(m\delta)} (\bw - \bC \bv) = m\lsp{t}\delta\bw \ge 1, \qquad
    \lsp{t}{(m\delta + \alpha)} (\bw - \bC \bv) \ge 0, \qquad
    \lsp{t}{(m\delta - \alpha)} (\bw - \bC \bv) \ge 0
  \end{equation*}
  from the assumption.
  Then
  \begin{gather*}
   \lsp{t}\bv\left(\bw - \tfrac12 \bC \bv\right)
     - \lsp{t}{\bv^0}\left(\bw - \tfrac12\bC\bv^0\right) - \sum p(\beta^{(t)})
\\
 \qquad{}    =    \sum_t \lsp{t}{\beta^{(t)}}(\bw - \bC\bv)
      - \# \{ t \mid \beta^{(t)} \in \Z_{>0}\delta\}
      + \frac12 \left(\sum \lsp{t}{\beta^{(t)}}\right) \bC
       \left(\sum {\beta^{(t)}}\right)
  \end{gather*}
  is nonnegative. Suppose that this is $0$. Then we must have
  $\lsp{t}\delta\bw = 1$ and all $\beta^{(t)}$ are $\delta$. This case
  is excluded in (1) from the assumption that $\bw - \bv$ is the Weyl
  group orbit of $\bw$. (See the proof of \propref{prop:affinestrata}.)
\end{proof}

\section{Crystal and the branching}\label{sec:crystal}

In this section\footnote{As we mentioned in the introduction, the
  results of this section was already obtained by Malkin
  \cite{Malkin}.}, we assume the graph has no edge loops. Then it
corresponds to a symmetric Kac--Moody Lie algebra $\mathfrak g$.
Let $\zeta$, $\zeta^\bullet$ as in \subsecref{subsec:Levi}. We have the
decomposition $I = I^0\sqcup I^+$. We have the Levi subalgebra
$\mathfrak g_{I^0}\subset\mathfrak g$ corresponding to the subdiagram
$I^0$.
Take and f\/ix $\bw\neq 0$. This is identif\/ied with a dominant integral
weight $\sum_i w_i\Lambda_i$ of $\mathfrak g$.

Let $\La_\zeta(\bv,\bw)\subset\M_\zeta(\bv,\bw)$ be the subvariety
$\pi_{0,\zeta}^{-1}(0)$. It is a Lagrangian subvariety
\cite[5.8]{Na-quiver} if $\M_\zeta(\bv,\bw)\neq\varnothing$. Let
$\Irr\La_\zeta(\bv,\bw)$ be the set of irreducible components of
$\La_\zeta(\bv,\bw)$ and let $\Irr\La_\zeta(\bw)$ be their disjoint
union $\bigsqcup_\bv \Irr\La_\zeta(\bv,\bw)$.
Kashiwara and Saito \cite{KS,Saito}, based on an earlier construction
due to Lusztig \cite{Lu-crystal}, constructed a $\mathfrak g$-crystal
structure on $\Irr\La_\zeta(\bw)$ which is isomorphic to the crystal of the
irreducible representation $V(\bw)$ of the quantum enveloping algebra
$\Uq$ with the highest weight $\bw$. (See \cite{Na-tensor} for a
dif\/ferent proof.)

We do not recall the def\/inition of the crystal and the construction in
\cite{KS,Saito} except two key ingredients, which are the weight and
a function, usually denoted by $\varepsilon_i$: If
$Y\in\Irr\La_\zeta(\bw)$, then
\(
  \operatorname{wt}(Y) = \bw - \bv.
\)
We consider the complex \eqref{eq:crystalcpx}
for a generic element $(B,a,b)$ in $Y$. Then we have
$\varepsilon_i(Y) = \codim \Ima\beta_i$. Note also that
$\langle h_i, \operatorname{wt}(Y)\rangle$ is the alternating sum of
the dimensions of terms in \eqref{eq:crystalcpx}, where the middle
term contribute in $+$.

From a general theory on the crystal, the $\g_{I^0}$-crystal of the
restriction of $V(\bw)$ to the subalgebra $\Un(\g_{I^0})\subset \Uq$
is given by forgetting $i$-arrows with $i\notin I^0$. Each connected
component, which is isomorphic to the crystal of an irreducible
highest weight representation of $\Un(\g_{I^0})$ has the unique
element $Y$ corresponding to the highest weight vector. It is
characterized by the property $\varepsilon_i(Y) = 0$ for any $i\in
I^0$. We will give its geometric characterization.

Let
\(
   \Lareg_{\zeta^\bullet}(\bv,\bw) \defeq
   \Mreg_{\zeta^\bullet}(\bv,\bw) \cap \pi_{0,\zeta^\bullet}^{-1}(0).
\)
Since $\pi_{\zeta^\bullet,\zeta}$ is an isomorphism on the preimage of
$\Mreg_{\zeta^\bullet}(\bv,\bw)$ (see \subsecref{subsec:partial}),
this can be identif\/ied with
$\pi_{\zeta^\bullet,\zeta}^{-1}(\Mreg_{\zeta^\bullet}(\bv,\bw))\linebreak[2]
\cap\La_\zeta(\bv,\bw)$. The latter is an open subvariety in
$\La_\zeta(\bv,\bw)$. Hence $\Lareg_{\zeta^\bullet}(\bv,\bw)$ is
of pure dimension with $\dim \Lareg_{\zeta^\bullet}(\bv,\bw)
= \dim \Mreg_{\zeta^\bullet}(\bv,\bw)/2$.
Let
\(
   \Irr \Lareg_{\zeta^\bullet}(\bv,\bw)
\)
be the set of irreducible components of~$\Lareg_{\zeta^\bullet}(\bv,\bw)$.
From what is explained above, this is a subset of
$\Irr\La_\zeta(\bv,\bw)$ consisting of those~$Y$ which intersect with the
open subset $\pi_{\zeta^\bullet,\zeta}^{-1}(\Mreg_{\zeta^\bullet}(\bv,\bw))$.

\begin{Theorem}\label{thm:mult_crystal}
  The set of irreducible components $Y\in\Irr\La_\zeta(\bv,\bw)$ with
  $\varepsilon_i(Y) = 0$ for any $i\in I^0$ is identified with
  $\Irr\Lareg_{\zeta^\bullet}(\bv,\bw)$.
\end{Theorem}

\begin{Corollary}
  The multiplicity of the irreducible highest weight module
  $V_{I^0}(\bw')$ of $\Un(\mathfrak g_{I^0})$ in the restriction of
  $V(\bw)$ is equal to the number of
  $Y\in\bigsqcup_{\bv}
  \Irr\Lareg_{\zeta^\bullet}(\bv,\bw)$ such that the restriction of
  $\operatorname{wt}(Y)$ to $I^0$ is $\bw'$.
\end{Corollary}

Note that the restriction of $\operatorname{wt}(Y)$ to $I^0$ is
dominant for $Y\in\Irr\Lareg_{\zeta^\bullet}(\bv,\bw)$ thanks to Lem\-ma~\ref{lem:betasurj}.

The theorem follows from the following:

\begin{Lemma}
  Suppose that $x = [B,a,b]\in \M_\zeta(V,W)$ is regular. Then $\Ima
  \beta_i = V_i$ for all $i\in I^0$ if and only if $(B,a,b)$ is
  $\zeta^\bullet$-stable.
\end{Lemma}

This was proved in \cite[2.9.4]{Na-qaff} (which was essentially a
collection of arguments in \cite{Na-alg}) in the special case $I =
I^0$. The same proof works in general. We reproduce the proof for the
sake of a~reader.

\begin{proof}
The `if' part follows from \lemref{lem:betasurj} as the case (a) is
excluded as we assumed $W\neq 0$.

The `only if' part: Since $(B,a,b)$ is $\zeta^\bullet$-semistable by
\lemref{lem:face} we take its Jordan--H\"older f\/iltration in
\thmref{thm:Rudakov}. We must have $k_W = N$, as the submodule
$V^{k_W+1}$ violates the $\zeta$-stability otherwise. Therefore if $N
\neq 0$, $\gr_0(B,a,b)$ has the $W$-component $0$, therefore it is
$0$-stable module for the subgraph $I^0$ extended to the whole graph
by $0$. But the regularity assumption implies that it must be $S_i$
for $i\in I^0$. In particular, $\beta_i$ is not surjective.
\end{proof}

Note that $Y\in\Irr\La_\zeta(\bv,\bw)$ is mapped into the closure of
the stratum
\begin{equation*}
   \Mreg_{\zeta^\bullet}(\bv',\bw) \times \left\{ \bigoplus_{i\in I^0}
   S_i^{\oplus \varepsilon_i(Y)} \right\} \qquad\text{with}\quad
   \bv' = \bv - \sum_{i\in I^0} \varepsilon_i(Y) \be_i
\end{equation*}
by the morphism $\pi_{\zeta^\bullet,\zeta}$.

\section{Convolution algebra and partial resolution}\label{sec:convolution}

We assume the graph has no edge loops and the stability parameters are
either $\zeta$, $\zeta^\bullet$ as in~\subsecref{subsec:Levi} or as in
\subsecref{subsec:sheaf} (and assume the graph is of an af\/f\/ine type).

\subsection{General results}

We f\/ix $\bw$ and consider {\it unions\/} of quiver varieties
$\M_{\zeta}(\bv,\bw)$, $\M_{\zeta^\bullet}(\bv,\bw)$, $\M_0(\bv,\bw)$
over va\-rious~$\bv$'s.
For $\M_0$ and $\M_{\zeta}$, they were introduced in
\cite[Section~2.5]{Na-qaff}:
\begin{equation*}
   \M_{0}(\bw) \defeq \bigcup_{\bv} \M_{0}(\bv,\bw),\qquad
   \M_{\zeta}(\bw) \defeq \bigsqcup_{\bv} \M_{\zeta}(\bv,\bw),
\end{equation*}
where we take the union in $\M_0(\bv,\bw)$ with respect to the closed
immersion $\M_0(V,W)\subset\M_0(V\oplus V',\bw)$ induced by the
extension by $0$ to the component $V'$. (In \cite[Section~2.5]{Na-qaff}
this was denoted by $\M_0(\infty,\bw)$.) These are inf\/inite unions of
varieties (of various dimensions), but there is no dif\/f\/iculty as we
can work on f\/initely many $\bv$'s in any of later constructions.

For $\M_{\zeta^\bullet}$ we mimic this construction. For
$\zeta^\bullet$ in \subsecref{subsec:Levi} we consider the closed
immersion $\M_{\zeta^\bullet}(V,W)\subset\M_{\zeta^\bullet}(V\oplus
V',\bw)$ with an $I^0$-graded (instead of $I$-graded) vector space
$V'$ and take the union. For $\zeta^\bullet$ in
\subsecref{subsec:sheaf} we similarly consider
$\M_{\zeta^\bullet}(V,W)\subset\M_{\zeta^\bullet}(V\oplus V',\bw)$
where $V'$ is a direct sum of various $S_i$ and $B_c$'s in
\eqref{eq:Uhbullet}.
This union is compatible with the morphism $\pi_{0,\zeta^\bullet}$ in
either cases. We have the induced morphism
$\pi_{0,\zeta^\bullet}\colon \M_{\zeta^\bullet}(\bw)\to \M_0(\bw)$. We
also have $\pi_{\zeta^\bullet,\zeta}\colon
\M_{\zeta}(\bw)\to\M_{\zeta^\bullet}(\bw)$.

Let
\(
  Z_{\zeta^\bullet,\zeta}(\bw) \defeq
  \M_{\zeta}(\bw)\times_{\M_{\zeta^\bullet}(\bw)} \M_\zeta(\bw)
\)
and similarly for $Z_{0,\zeta}(\bw)$. These are union of various
subvarieties $Z_{\zeta^\bullet,\zeta}(\bv^1,\bv^2;\bw)$
or $Z_{0,\zeta}(\bv^1,\bv^2;\bw)$
in $\M_\zeta(\bv^1,\bw)\times\M_\zeta(\bv^2,\bw)$, where the f\/iber
product is def\/ined over a space $\M_{\zeta^\bullet}(\bv,\bw)$
or $\M_0(\bv,\bw)$ with some large $\bv$ compared with $\bv^1$, $\bv^2$.
By \propref{prop:semismall} $Z_{0,\zeta}(\bv^1,\bv^2;\bw)$,
$Z_{\zeta^\bullet,\zeta}(\bv^1,\bv^2;\bw)$ are at most
half dimensional in $\M_\zeta(\bv^1,\bw)\times\M_\zeta(\bv^2,\bw)$.
Moreover $Z_{\zeta^\bullet,\zeta}(\bv^1,\bv^2;\bw)$ is a subvariety of
$Z_{0,\zeta}(\bv^1,\bv^2;\bw)$.
The latter space was introduced in \cite{Na-alg} as an analog of the
Steinberg variety appearing in a geometric construction of the Weyl group.

We consider the top degree Borel--Moore homology groups
$H_\topdeg(Z_{\zeta^\bullet,\zeta}(\bw))$ and
$H_\topdeg(Z_{0,\zeta}(\bw))$. More precisely the former is the subspace
\[
\prod_{\bv^1,\bv^2}' H_{\dim_\C \M_\zeta(\bv^1,\bw)\times\M_\zeta(\bv^2,\bw)}
(Z_{\zeta^\bullet,\zeta}(\bv^1,\bv^2;\bw),\Q),
\]
of the direct products consisting elements $(F_{\bv,\bv'})$ such
that{\samepage
\begin{enumerate}\itemsep=0pt
\item[1)] for f\/ixed $\bv^1$, $F_{\bv^1,\bv^2} = 0$ for all but f\/initely
many choices of $\bv^2$,
\item[2)] for f\/ixed $\bv^2$, $F_{\bv^1,\bv^2} = 0$  for all but f\/initely
many choices of $\bv^1$.
\end{enumerate}
Similarly for the latter.}

Since $\M_\zeta(\bv,\bw)$ is smooth, we have an associative algebra
structure on $H_\topdeg(Z_{\zeta^\bullet,\zeta}(\bw))$ and
$H_\topdeg(Z_{0,\zeta}(\bw))$ with the unit given by the sum of
diagonals in $\M_\zeta(\bv,\bw)\times\M_\zeta(\bv,\bw)$ for va\-rious~$\bv$.
We have an injective algebra homomorphism
\begin{equation}\label{eq:subalg}
    H_\topdeg(Z_{\zeta^\bullet,\zeta}(\bw)) \to
    H_\topdeg(Z_{0,\zeta}(\bw)).
\end{equation}

We can analyze irreducible representations of these algebras and the
branching rules with respect to the above homomorphism by a general
theory in \cite{CG} based on \cite{BM}. We prepare several notations
and concepts, and then state the result.

For a point $x$ in a stratum $\M_0(\bw)_{(\widehat G)}$ of
$\M_0(\bw)$ let $d_x(\bv,\bw)\! =\! \dim_\C \M_\zeta(\bv,\bw)
-\linebreak[2] \dim_\C \M_0(\bw)_{(\widehat G)}$ and
$H_\topdeg(\pi_{0,\zeta}^{-1}(x))$ be the direct sum of the cohomology
groups of the f\/ibers $\pi_{0,\zeta}^{-1}(x)\cap\M_\zeta(\bv,\bw)$ of
the degree $d_x(\bv,\bw)$.  By the convolution product
$H_\topdeg(\pi_{0,\zeta}^{-1}(x))$ is a
$H_\topdeg(Z_{0,\zeta}(\bw))$-module. The cohomology vanishes in
degree above $d_x(\bv,\bw)$ and has a basis by
$d_x(\bv,\bw)$-dimensional irreducible components in degree
$d_x(\bv,\bw)$ by \propref{prop:semismall}, and does not vanish in
$d_x(\bv,\bw)$ if the graph is of f\/inite or af\/f\/ine type by
\remref{rem:relevant}.

For a simple local system $\phi$ on $\M_0(\bw)_{(\widehat G)}$ (i.e.,
an irreducible representation of the fundamental group of
$\M_0(\bw)_{(\widehat G)}$), let $IC(\M_0(\bw)_{(\widehat G)}, \phi)$
be the corresponding intersection cohomology complex.

The fundamental group of $\M_0(\bw)_{(\widehat G)}$ acts on
$H_\topdeg(\pi_{0,\zeta}^{-1}(x))$ by monodromy. It is a permutation
of the above basis elements. We have a decomposition
\(
  H_\topdeg(\pi_{0,\zeta}^{-1}(x)) = \bigoplus \phi\otimes V_{(\widehat{G}),\phi}
\)
into the direct sum of simple local system $\phi$ tensored with the
multiplicity vector space $V_{(\widehat G),\phi}$.
A pair $(\M_0(\bw)_{(\widehat G)},\phi)$ of a stratum
$\M_0(\bw)_{(\widehat G)}$ and a simple local system $\phi$ on it is
{\it relevant\/} for $\pi_{0,\zeta}$ if $V_{(\widehat G),\phi}\neq 0$.
Then we have the decomposition theorem
\begin{equation*}
  (\pi_{0,\zeta})_*(\C_{\M_\zeta(\bw)}[\dim])
  = \bigoplus IC(\M_0(\bw)_{(\widehat G)},\phi) \otimes V_{(\widehat
    G),\phi},
\end{equation*}
where $\C_{\M_\zeta(\bw)}[\dim]$ is the direct sum
$\bigoplus_\bv \C_{\M_\zeta(\bv,\bw)}[\dim \M_\zeta(\bv,\bw)]$.
We also know \cite[8.9.8]{CG} that
\begin{equation}\label{eq:convalg}
  H_{\topdeg}(Z_{0,\zeta}(\bw)) \cong \bigoplus \End_\C(V_{(\widehat G),\phi}).
\end{equation}
Thus $\{ V_{(\widehat G,\phi)} \mid \text{$((\M_0)_{(\widehat
    G)},\phi)$ is relevant for $\pi_{0,\zeta}$}\}$ is the set of
isomorphism classes of irreducible representations of
$H_{\topdeg}(Z_{0,\zeta}(\bw))$.

We have similar formulas for $\pi_{\zeta^\bullet,\zeta}\colon
\M_{\zeta}(\bw)\to \M_{\zeta^\bullet}(\bw)$. We denote by $(\widehat H)$
a conjugacy class of $\GL(V)$ corresponding to a stratum of
$\M_{\zeta^\bullet}(\bw)$.
Similar to above we can def\/ine
$H_{\topdeg}(\pi_{\zeta^\bullet,\zeta}^{-1}(y))$, which is a
$H_\topdeg(Z_{\zeta^\bullet,\zeta}(\bw))$-module for
$y\in \M_{\zeta^\bullet}(\bw)_{(\widehat H)}$.

If we decompose the direct image of
$IC(\M_{\zeta^\bullet}(\bw)_{(\widehat H)},\psi)$ as
\begin{equation}\label{eq:bullet}
   (\pi_{0,\zeta^\bullet})_*(IC(\M_{\zeta^\bullet}(\bw)_{(\widehat
     H)},\psi))
   = \bigoplus IC(\M_0(\bw)_{(\widehat G)},\phi)\otimes
     V_{(\widehat G),\phi}^{(\widehat H),\psi}
\end{equation}
with the multiplicity vector space $V_{(\widehat G),\phi}^{(\widehat
  H),\psi}$, we have the following branching rule from the double
decomposition formula (see \cite[1.11]{BM})
\begin{equation}\label{eq:branching}
   \Res V_{(\widehat G),\phi} = \bigoplus
   V_{(\widehat G),\phi}^{(\widehat H),\psi} \otimes V_{(\widehat H),\psi},
\end{equation}
where $\Res$ is the restriction functor from
$H_{\topdeg}(Z_{0,\zeta}(\bw))$-modules to
$H_{\topdeg}(Z_{\zeta^\bullet,\zeta}(\bw))$-modules via the injective
homomorphism \eqref{eq:subalg}. Since $V_{(\widehat H),\psi}$ are
pairwise non-isomorphic, $V_{(\widehat G),\phi}^{(\widehat H),\psi}$
is the multiplicity space
\(
   \Hom_{H_{\topdeg}(Z_{\zeta^\bullet,\zeta}(\bw))}
   (V_{(\widehat H),\psi}, \Res V_{(\widehat G),\phi}).
\)

Note also that $V_{(\widehat G),\phi}^{(\widehat H),\psi}$ is the
$\phi$-isotypical component of
\(
   H^{\topdeg}(\pi_{0,\zeta^\bullet}^{-1}(x),
   IC(\M_{\zeta^\bullet}(\bw)_{(\widehat H)},\psi))[d_y]
\)
the top degree cohomology of $\pi_{0,\zeta^\bullet}^{-1}(x)$ with the
coef\/f\/icient in the shifted IC sheaf, where
$d_y = \dim \M_\zeta(\bw) - \dim \M_{\zeta^\bullet}(\bw)_{(\widehat H)}$.
(More precisely we consider the degree `$\topdeg$' and $d_y$ for each
$\bv$ separately.) See \cite[1.10]{BM}.

\subsection{The restriction to a Levi factor}\label{subsec:rest_Levi}

Let us consider the case when $\zeta$, $\zeta^\bullet$ are as in
\subsecref{subsec:Levi} in this subsection.

In \cite{Na-alg} we constructed an algebra homomorphism
\begin{equation*}
   \bU(\g) \to H_{\topdeg}(Z_{0,\zeta}(\bw)).
\end{equation*}

The homomorphism was given on generators by
\begin{equation*}
   e_i \mapsto \sum_{\bv^2} [\Pa(\bv^2,\bw)],\quad
   f_i \mapsto \sum_{\bv^2} \pm [\omega\Pa(\bv^2,\bw)], \quad
   h \mapsto \sum_{\bv} \langle h, \bw - \bv\rangle [\Delta(\bv,\bw)],
\end{equation*}
where $\Pa(\bv^2,\bw)$ is the `Hecke correspondence' parametrizing
pairs $([B^1,a^1,b^1],[B^2,a^2,b^2])$ such that $(B^1,a^1,b^1)$ is a
submodule of $(B^2,a^2,b^2)$ with the quotient isomorphic to $S_i$, a
module with $\C$ on the vertex $i$, and $0$ on the other vertices,
$\omega$ is the exchange of factors of $\M_\zeta(\bv^1,\bw)\times
\M_\zeta(\bv^2,\bw)$ and $\Delta(\bv,\bw)$ is the diagonal in
$\M_\zeta(\bv,\bw)\times \M_\zeta(\bv,\bw)$. The $\pm$-sign in the
def\/inition of $f_i$ is not important in the discussion below, so its
precise def\/inition is omitted.
From this def\/inition it is clear that we have an algebra homomorphism
\begin{equation*}
   \bU(\g_{I^0}+\mathfrak h) \to H_{\topdeg}(Z_{\zeta^\bullet,\zeta}(\bw)),
\end{equation*}
where $\mathfrak h$ is the Cartan subalgebra of $\g$ and $\g_{I^0}$ is
the Levi factor corresponding to the subset $I^0\subset I$.

We must be careful when we apply the result in the previous
subsection, as this homomorphism is {\it not\/} an isomorphism.

\begin{Remark}\label{rem:degree}
  Later we consider the case when $\g$ is an af\/f\/ine Lie algebra. Our
  af\/f\/ine Lie algebra, as in \cite{Kac}, contains the degree operator
  $d$. It is mapped to
\(
   \sum_\bv \langle d, \bw - \bv\rangle [\Delta(\bv,\bw)]
   = - \sum_{\bv} v_0[\Delta(\bv,\bw)],
\)
where $0$ is the special vertex of the af\/f\/ine graph as in
\subsecref{subsec:ALE}.
\end{Remark}

Let us denote by $IC(\M_{\zeta^\bullet}(\bw)_{(\widehat H)})$ the IC
sheaf corresponding to the {\it trivial\/} local system.

\begin{Theorem}\label{thm:branch}
  Let $V^\g(\lambda)$ \textup(resp.\ $V^{(\g_{I^0}+\mathfrak
    h)}(\lambda)$\textup) denote the irreducible integrable highest
  weight module of $\g$ \textup(resp.\ $\g_{I^0}+\mathfrak h$\textup)
  with the highest weight $\lambda$. Then we have
  \begin{gather*}
    (\pi_{0,\zeta^\bullet})_*(IC(\Mreg_{\zeta^\bullet}(\bv^0,\bw)))
   =
    \bigoplus_{\bv'}
    \Hom_{\g_{I^0}+\mathfrak
    h}(V^{(\g_{I^0}+\mathfrak h)}(\bw \!-\! \bv^0), V^\g(\bw\!-\!\bv'))
  \!\otimes\! IC(\Mreg_0(\bv',\bw))
  \\
  \phantom{(\pi_{0,\zeta^\bullet})_*(IC(\Mreg_{\zeta^\bullet}(\bv^0,\bw)))=}{} \oplus
  \bigoplus (\text{\rm
    $IC$ sheaves associated with {\it non\/}-regular strata}).
  \end{gather*}
\end{Theorem}

\begin{proof}
Let $\M_0(\bw)_{(\widehat G)}$ be a stratum of $\M_0(\bw)$. Take a point
$x\in \M_0(\bw)_{(\widehat G)}$. The decomposition
$H_{\topdeg}(\pi_{0,\zeta}^{-1}(x)) = \bigoplus
\phi\otimes V_{(\widehat G),\phi}$ in the previous subsection is a
decomposition of a $\g$-module.
Suppose $\M_0(\bw)_{(\widehat G)}$ is regular, i.e., $=
\Mreg_0(\bv',\bw)$ for some $\bv'$.
Then it is known that $H_{\topdeg}(\pi_{0,\zeta}^{-1}(x))$ is an
irreducible integrable representation of $\g$ with the highest weight
vector $[x] \in H_{\topdeg}(\pi_{0,\zeta}^{-1}(x)\cap
\M_\zeta(\bv',\bw))$ \cite[10.2]{Na-alg}, where
$x$ can be considered as a point in $\M_\zeta(\bv',\bw)$ as
$\pi_{0,\zeta}\colon \M_\zeta(\bv',\bw)\to \M_0(\bv',\bw)$ is an
isomorphism on the preimage of $\Mreg_0(\bv',\bw)$.
In particular, $V_{(\widehat G),\phi} = 0$ unless $\phi$ is the
trivial local system in this case.

Similarly if a stratum $\M_{\zeta^\bullet}(\bw)_{(\widehat H)}$ is of
the form $\Mreg_{\zeta^\bullet}(\bv^0,\bw)$ for some $\bv^0$, then
$H_{\topdeg}(\pi_{\zeta^\bullet,\zeta}^{-1}(y))$ is the irreducible
integrable representation $V^{(\g_{I^0}+\mathfrak h)}(\bw-\bv^0)$ if
$y\in \Mreg_{\zeta^\bullet}(\bv^0,\bw)$.
The highest weight vector is $[y]$ as above.

Consider a stratum $\M_{\zeta^\bullet}(\bw)_{(\widehat H)}$ which
contains $IC(\Mreg_0(\bv',\bw))$ in the decomposition
\eqref{eq:bullet}. We can write it as
$\Mreg_{\zeta^\bullet}(\bv^0,\bw)\times (\M_0)_{(\widehat H')}$ as in
\eqref{eq:Levistratum}. We consider the diagram
\eqref{eq:basechange}. Since $\kappa$ is a f\/inite birational morphism,
we have
\begin{equation*}
  IC(\overline{\Mreg_{\zeta^\bullet}(\bv^0,\bw)
     \times (\M_0)_{(\widehat H')}})
   = \kappa_*\big(
     IC(\M(\bv^0,\bw))\boxtimes IC(\overline{(\M_0)_{(\widehat H')}})
   \big),
\end{equation*}
and hence
\begin{equation*}\label{eq:nonregular}
  (\pi_{0,\zeta^\bullet})_* IC(\overline{\Mreg_{\zeta^\bullet}(\bv^0,\bw)
     \times (\M_0)_{(\widehat H')}}) = \kappa'_*\big(
   (\pi_{0,\zeta^\bullet})_*IC(\M(\bv^0,\bw))\boxtimes
   IC(\overline{(\M_0)_{(\widehat H')}})\big).
\end{equation*}
Since $\kappa'$ is also a f\/inite morphism, this contains only IC
sheaves on nonregular strata unless the second factor
$(\M_0)_{(\widehat H')}$ is trivial.
Hence \eqref{eq:branching} now becomes
\begin{equation}\label{eq:restriction}
   \operatorname{Res} V^{\g}(\bw - \bv')
   = \bigoplus_{\bv^0} V_{\bw-\bv'}^{\bw-\bv^0}
   \otimes V^{(\g_{I^0}+\mathfrak h)}(\bw - \bv^0),
\end{equation}
where $V_{\bw-\bv'}^{\bw-\bv^0}$ is the multiplicity of
$IC(\Mreg_0(\bv',\bw))$ in the decomposition of\linebreak[2]
$(\pi_{0,\zeta^\bullet})_*\linebreak[2](IC(\Mreg_{\zeta^\bullet}(\bv^0,
\linebreak[2]\bw)))$ in
\eqref{eq:bullet}. Thus
$V_{\bw-\bv'}^{\bw-\bv^0} =
\Hom_{\g_{I^0}+\mathfrak h}
(V^{(\g_{I^0}+\mathfrak h)}(\bw - \bv^0), V^\g(\bw-\bv'))$.
Note here that $V^{(\g_{I^0}+\mathfrak h)}(\bw - \bv^0)$ are
non-isomorphic for dif\/ferent $\bv^0$.
\end{proof}

We assume the graph is of af\/f\/ine type until the end of this
subsection.
This is because we know the structure of $H_{\topdeg}(\pi_{0,\zeta}^{-1}(x))$
for a point $x$ from a nonregular stratum only in af\/f\/ine type.
We further suppose $I^0\neq I$ to avoid the trivial
situation $\pi_{0,\zeta^\bullet} = \id$. We thus have
\begin{equation}\label{eq:zetastratum}
  \M_{\zeta^\bullet}(\bw) = \bigsqcup_{\bv}
  \Mreg_{\zeta^\bullet}(\bv,\bw),
\end{equation}
where $\bv$ is such that $\bw - \bv$ is an $I^0$-dominant weight of
$V(\bw)$.

Let us consider $\pi_{0,\zeta}\colon \M_\zeta(\bw) \to \M_0(\bw)$.
The stratif\/ication of $\M_0(\bw)$ induced from \eqref{eq:Uhbullet} is
\begin{equation}\label{eq:0stratum}
   \M_0(\bw) = \bigsqcup \Mreg_0(\bv',\bw)\times
   S^{|\lambda|}_\lambda(\C^2\setminus\{0\}/\Gamma).
\end{equation}
The criterion for $\Mreg_0(\bv',\bw)\neq\varnothing$ was given in
\propref{prop:affinestrata}.
Take a point $x$ from the stratum $\Mreg_0(\bv',\bw)\times
S^{|\lambda|}_\lambda(\C^2\setminus\{0\}/\Gamma)$ with $\lambda =
(\lambda_1,\dots, \lambda_l)$ ($l=l(\lambda)$), and consider the
inverse image under $\pi_{0,\zeta}\colon \M_\zeta(\bv,\bw)\to
\M_0(\bv,\bw)$.
Then $\widehat{I}$ is the disjoint union of $I$ and $l$ copies of the
Jordan quiver. The dimension vectors are
\begin{equation*}
  \widehat{\bv} = (\bv - \bv' - |\lambda|\delta,
  \lambda_1,\dots,\lambda_l),
\qquad
  \widehat{\bw} = (\bw - \bC\bv',
  \lsp{t}{\delta}{\bw},\dots,\lsp{t}{\delta}{\bw}),
\end{equation*}
where the f\/irst parts are the $I$-components and the remaining $l$
entries are the Jordan quiver components. In particular
$\pi_{0,\zeta}^{-1}(x)$ is isomorphic to the product of
$\pi_{0,\zeta}^{-1}(0)$ ($\subset\M_\zeta(\bw-\bC\bv')$) and
$\pi_{0,\zeta}^{-1}(0)$ ($\subset
M_\zeta^{\text{norm}}(\lambda_i,\lsp{t}{\delta}\bw)$) ($i=1,\dots, l$),
where the latter are the quiver varieties associa\-ted with Jordan
quivers. It is known that the latter's are irreducible and
half-dimensional in
$M_\zeta^{\text{norm}}(\lambda_i,\lsp{t}{\delta}\bw)$. (See
\cite[Exercise~5.15]{Lecture} and the references therein.)  Thus
$H_{\topdeg}(\pi^{-1}_{0,\zeta}(x))$ is isomorphic to the tensor
product of $H_{\topdeg}(\pi^{-1}_{0,\zeta}(0))$ and
\(
   \C[\prod_i \pi_{0,\zeta}^{-1}(0)(\subset
   M_\zeta^{\text{norm}}(\lambda_i,\lsp{t}{\delta}\bw))].
\)
The af\/f\/ine Lie algebra $\g$ acts trivially on the second factor (or
more precisely the degree operator $d$ acts by $-\lambda_i$ by
\remref{rem:degree}), so $H_{\topdeg}(\pi^{-1}_{0,\zeta}(x))$ is
isomorphic to $V^\g(\bw - \bv' - |\lambda|\delta)$. In particular, we
get
\begin{equation*}
  (\pi_{0,\zeta})_*(\C_{\M_\zeta(\bw)}[\dim])
  = \bigoplus_{\bv',\lambda} IC(\Mreg_0(\bv',\bw)\times
  S^{|\lambda|}_\lambda(\C^2\setminus\{0\}/\Gamma))
  \otimes V^\g(\bw - \bv' - |\lambda|\delta).
\end{equation*}
Note that the argument also shows that nontrivial local systems do not
appear in the direct summand as for regular strata.
Moreover, the closure of
$S^{|\lambda|}_\lambda(\C^2\setminus\{0\}/\Gamma))$ has only f\/inite
quotient singularities, and hence is rationally smooth:
$IC(S^{|\lambda|}_\lambda(\C^2\setminus\{0\}/\Gamma))) =
\C_{\overline{S^{|\lambda|}_\lambda(\C^2\setminus\{0\}/\Gamma))}}[\dim]$,
where $\dim$ means $\dim
S^{|\lambda|}_\lambda(\C^2\setminus\{0\}/\Gamma))$.

Hence we have
\begin{Theorem}
Suppose the graph is of affine type and $I^0\neq I$. Then

\textup{(1)} The strata of $\M_{\zeta^\bullet}(\bw)$ and $\M_0(\bw)$
are given by \eqref{eq:zetastratum}, \eqref{eq:0stratum} respectively.
The criterion of the nonemptiness of $\Mreg_{\zeta^\bullet}(\bv,\bw)$,
$\Mreg_0(\bv,\bw)$ is given in Proposition~{\rm \ref{prop:affinestrata}}.

\textup{(2)} We have
\begin{gather*}
  (\pi_{0,\zeta^\bullet})_*(IC(\Mreg_{\zeta^\bullet}(\bv^0,\bw)))
 =
   \bigoplus_{\bv,\lambda}
  IC(\Mreg_0(\bv',\bw)) \boxtimes
\C_{\overline{S^{|\lambda|}_\lambda(\C^2\setminus\{0\}/\Gamma))}}[\dim]
\\
\phantom{(\pi_{0,\zeta^\bullet})_*(IC(\Mreg_{\zeta^\bullet}(\bv^0,\bw)))=}{} \otimes_\C
  \Hom_{(\g_{I^0}+\mathfrak h)}
  (V^{(\g_{I^0}+\mathfrak h)}(\bw - \bv^0), V^\g(\bw-\bv'-|\lambda|\delta)).
\end{gather*}
\end{Theorem}

\begin{Remarks}\label{rem:check}
(1) Suppose $I^0=\varnothing$ and hence $\zeta = \zeta^\bullet$,
$\Mreg_{\zeta}(\bv^0,\bw) = \M_{\zeta}(\bv^0,\bw)$.
Then the above implies that the restriction of
$H_*(\pi_{0,\zeta}^{-1}(0))$ to a $\g$-module decomposes as
\begin{equation*}
  \bigoplus_{\bv,\lambda}
  H^*(i_0^! IC(\Mreg_0(\bv,\bw))) \otimes V^\g(\bw-\bv-|\lambda|\delta),
\end{equation*}
where $i_0\colon\{0\}\to \M_0(\bw)$ is the inclusion. This is nothing
but \cite[Theorem~5.2]{MR1989196}.

(2) It is known that $H_*(\pi_{0,\zeta}^{-1}(0))$ is the tensor
product of $\g$-modules corresponding to ones with $\bw = \Lambda_{\mu_p}$
for various $\mu_p\in [0,r-1]$. (See the proof of
\cite[14.1.2]{Na-qaff}.) If we
further assume $\g$ is of af\/f\/ine type $A_{r-1}^{(1)} =
\algsl(r)_\aff$, then those $\g$-modules can be computed by various
means. For example, we use \propref{prop:affinestrata} to f\/ind
$\Mreg_0(\bv,\Lambda_{\mu_p}) = \varnothing$ unless $\bv=0$. Therefore
$H_*(\pi_{0,\zeta}^{-1}(0))$ is isomorphic to
\(
   \bigoplus_\lambda V^\g(\Lambda_{\mu_p}-|\lambda|\delta).
\)
The direct sum of $1$-dimensional spaces for each partition $\lambda$
is isomorphic to the Fock space of the Heisenberg algebra. Thus this
is isomorphic to the restriction of the highest weight representation
$V^{\gl(r)_\aff}(\Lambda_i)$ of $\gl(r)_\aff$ to $\algsl(r)_\aff$.
(See \secref{sec:level-rank}.) Combined with \eqref{eq:duality'} we
f\/ind that $H^*(i_0^! IC(\Mreg_0(\bv,\bw)))$ is isomorphic to
$V^{\algsl(l)_\aff}(\overline\lambda)_{\overline\mu}$. This
observation was used in \cite{braverman-2007} to conf\/irm a (weak form)
of the conjecture proposed there in type $A$.
\end{Remarks}

\subsection[Restriction to the affine Lie algebra of a Levi factor]{Restriction to the af\/f\/ine Lie algebra of a Levi factor}

Let us consider the case in \subsecref{subsec:sheaf}. The analysis is
almost the same as in the previous subsection.

The stratif\/ication of $\M_0(\bw)$ was already given in
\eqref{eq:0stratum}.
The f\/iber $\pi_{0,\zeta}^{-1}(x)$ of a point $x$ in the stratum
$\Mreg_0(\bv',\bw)\times
S^{|\lambda|}_\lambda(\C^2\setminus\{0\}/\Gamma)$ (see
\eqref{eq:0stratum}) is isomorphic to the integrable highest weight
representation $V^\g(\bw - \bv' - |\lambda|\delta)$ of the af\/f\/ine Lie
algebra $\g$.
Strictly speaking, this result does not follows directly from
\cite{Na-alg} as the generic stability parameter $\zeta$ used there is
dif\/ferent from the one used here. Rather it is the one in
\subsecref{subsec:Levi}. Therefore
$H_{\topdeg}(\pi_{0,\zeta}^{-1}(x))$, {\it a priori\/}, only has an
$\mathfrak h$-module structure.
Let us denote the generic parameter $\zeta$ in the previous subsection
by~$\zeta^{+}$.

\begin{Lemma}
  Let us fix $\bv$ and consider $\pi_{0,\zeta}^{-1}(x)\subset
  \M_{\zeta}(\bv,\bw)$ and $\pi_{0,\zeta^+}^{-1}(x)\subset
  \M_{\zeta^+}(\bv,\bw)$. There is a graded vector space isomorphism
\(
   H_{*}(\pi_{0,\zeta}^{-1}(x)
   )
   \cong
   H_{*}(\pi_{0,\zeta^+}^{-1}(x)
   )
\)
\textup(which is canonical in the way explained during the
proof\textup). In particular, $H_{\topdeg}(\pi_{0,\zeta}^{-1}(x))$ is
isomorphic to $V^\g(\bw - \bv' - |\lambda|\delta)$ as an $\mathfrak
h$-module.
\end{Lemma}

By \eqref{eq:convalg} $H_{\topdeg}(Z_{0,\zeta}(\bw))$ has an induced
homomorphism from $\bU(\g)$.

\begin{proof}
  We know that $\pi_{0,\zeta}^{-1}(x)\subset \M_{\zeta}(\bv,\bw)$ is
  isomorphic to the product of $\pi_{0,\zeta}^{-1}(0)\subset
  \M_{\zeta}(\bv-\bv'$, $\bw-\bc\bv')$ and $\pi_{0,\zeta}^{-1}(0)\subset
  M_\zeta^{\text{norm}}(\lambda_i,\lsp{t}{\delta}\bw)$ ($i=1,\dots,
  l$) as in the previous subsection. The same holds for $\zeta^+$. The
  latter factors are the same for our $\zeta$ and $\zeta^+$, so we
  need to worry about the f\/irst factor. We know that the inclusion
  $\pi_{0,\zeta}^{-1}(0)\subset \M_\zeta(\bv-\bv',\bw-\bC\bv')$ and
  one for $\zeta^+$ are homotopy equivalences \cite[5.5]{Na-quiver},
  and the $C^\infty$-structure of $\M_\zeta(\bv-\bv',\bw-\bC\bv')$ is
  independent of the choice of generic parameter $\zeta$. In
  particular, we get the f\/irst assertion. Since this is the weight
  space, we have the second assertion.
\end{proof}

There is another way to construct an isomorphism
$H_{*}(\pi_{0,\zeta}^{-1}(x))\cong
H_{*}(\pi_{0,\zeta^+}^{-1}(x))$.
Let us sketch the construction. We introduce the complex parameter
$\zeta_\C\in\C^I$ and consider the quotients of $\mu^{-1}(\C^I)$ by
$G$ associated with the stability parameters $\zeta$, $\zeta^+$ and
$0$. Let us denote the corresponding variety by $\mathfrak
N_{\zeta}(\bv,\bw)$, etc.
The morphisms
\(
 \pi_{0,\zeta}\colon\mathfrak N_{\zeta}(\bv,\bw) \to
  \mathfrak N_0(\bv,\bw),
\)
\(
 \pi_{0,\zeta^+}\colon\mathfrak N_{\zeta^+}(\bv,\bw) \to
  \mathfrak N_0(\bv,\bw)
\)
become {\it small\/}, and the pushforward of the constant sheaves
$\C_{\mathfrak N_\zeta(\bv,\bw)}[\dim]$, $\C_{\mathfrak N_{\zeta^+}(\bv,\bw)}[\dim]$
are both the $IC$ sheaf of the closure of $\mathfrak
N_0^{\operatorname{s}}(\bv,\bw)$.
Therefore both $H_{*}(\pi_{0,\zeta}^{-1}(x))$ and
$H_{*}(\pi_{0,\zeta^+}^{-1}(x))$ are the stalk of
$IC(\mathfrak N_0^{\operatorname{s}}(\bv,\bw))$ at $x$, hence are
isomorphic.
Moreover as $H_{\topdeg}(Z_{0,\zeta}(\bw))$ is isomorphic to
$\End_{D^b(\M_0(\bw))}((\pi_{0,\zeta})_*(\C_{\M_\zeta(\bw)}[\dim])$
(\cite[8.9.6]{CG}), this construction also gives a natural isomorphism
$H_{\topdeg}(Z_{0,\zeta}(\bw))\cong H_{\topdeg}(Z_{0,\zeta^+}(\bw))$
respecting the convolution product.
This kind of construction is well-known in the geometric construction
of the Springer correspondence (see~\cite{Lu-Green}). And
$\mathfrak N_0(\bv,\bw)$ is an analog of the Beilinson--Drinfeld
Grassmannian in the double af\/f\/ine Grassmannian.

\begin{Remark}\label{rem:Frenkel-Kac}
  It is probably possible to def\/ine a homomorphism $\bU(\g) \to
  H_{\topdeg}(Z_{0,\zeta}(\bw))$ also for this $\zeta$ more directly.
  When $\bw$ is of level $1$, this can be done by combining the
  Heisenberg algebra representation on the cohomology group of Hilbert
  schemes of points \cite{Lecture} and the Frenkel--Kac construction
  \cite[\S~14]{Kac} (see \cite{Grojnowski}). Furthermore, it def\/ines
  the same representation on
  $H_{\operatorname{mid}}(\M_\zeta(\bw))\cong
  H_{\operatorname{mid}}(\M_{\zeta^+}(\bw))$ as one in \cite{Na-alg}
  if the graph is of af\/f\/ine type $A$ \cite{Nagao-2007}.
\end{Remark}

Since $V^\g(\bw - \bv' - |\lambda|\delta)$ is not irreducible as an
$\mathfrak h$-module, the previous argument showing that $V_{(\widehat
  G),\phi} = 0$ unless $\phi$ is the trivial local system does {\it
  not\/} work. But we can deduce it from the corresponding vanishing
for $\zeta^+$ in the previous subsection. This can be done (at least)
in two ways. One is to use the ref\/lection functor \cite{Na-reflect} to
construct an isomorphism $\M_\zeta(\bv,\bw)\to\M_{\zeta^+}(\bv',\bw)$
(for an appropriate $\bv'$). It is compatible with $\pi_{0,\zeta}$,
$\pi_{0,\zeta^+}$, so the assertion follows. The second one is to use
$\mathfrak N_{\zeta}$ as above. Then both
$(\pi_{0,\zeta})_*(\C_{\M_\zeta(\bv,\bw)}[\dim])$,
$(\pi_{0,\zeta^+})_*(\C_{\M_{\zeta^+}(\bv,\bw)})$ are obtained from
the above IC sheaf $IC(\mathfrak N_0^{\operatorname{s}}(\bv,\bw))$ by
the restriction. Therefore they are isomorphic.

Next consider $\pi_{\zeta^\bullet,\zeta}$.
The stratif\/ication of $\M_{\zeta^\bullet}(\bw)$ induced from
\eqref{eq:Uhbullet} is{\samepage
\begin{equation}\label{eq:bulletstratum}
  \M_{\zeta^\bullet}(\bw)
= \bigsqcup \Mreg_{\zeta^\bullet}(\bv^0,\bw)
  \times S^{|\mu|}_\mu
  \Big[X_{\zeta^\circ}\setminus \bigcup_{i\in I_0^0} C_i \Big],
\end{equation}
where $\mu$ is a partition.}

Take a point $y$ from the stratum of the form in
\eqref{eq:bulletstratum}. The f\/iber
$\pi_{\zeta^\bullet,\zeta}^{-1}(y)$ is isomorphic to the product of
$\pi_{0,\zeta}^{-1}(0)\subset \M_\zeta(\widehat\bw_c)$ ($c\in \mathcal
C$) and $\pi_{0,\zeta}^{-1}(0)\subset
M_\zeta^{\text{norm}}(\mu_i,\lsp{t}{\delta}\bw)$ ($i=1,\dots, l$).
Here $\widehat\bw_c$ is the entries corresponding to $c$ in
\eqref{eq:hatdim}, and $\M_\zeta(\widehat\bw_c)$ is the quiver variety
for the af\/f\/ine graph corresponding to $c$.
Let $\affLevi$ be the af\/f\/ine Lie algebra of the Levi factor
$\g_{I_0^0}$ of the f\/inite dimensional Lie algebra $\g_{I^0}$.
Then $H_{\topdeg}(\pi_{\zeta^\bullet,\zeta}^{-1}(y))$ is the
integrable highest weight representation
$V^{\affLevi+\mathfrak h}(\bw - \bv^0 - |\mu|\delta)$.
More precisely it is the direct product of the integrable highest
weight representation for the af\/f\/ine Lie algebra corresponding to $c$
with the highest weight $\widehat\bw_c$ for $c\in\mathcal C$.
Note that each factor has level $\lsp{t}{\delta}\bw$, independent of
$c$, thus the af\/f\/ine Lie algebra $\affLevi$ contains only
the one dimensional central extension even if $\g_{I^0_0}$ has several
simple factors.

As in the case of $\pi_{0,\zeta}$ we need an explanation of the
$(\affLevi+\mathfrak h)$-module structure on
$H_{\topdeg}(\pi_{\zeta^\bullet,\zeta}^{-1}(y))$. It can be done in
two ways as above. We still need to compare it with the $\g$-module
structure considered for $\pi_{0,\zeta}$.
It is clear that the $\h$-module structure is compatible in either
ways. But the author does not know how to compare the $\g$-module and
$\affLevi$-module structure.
Therefore the equality \eqref{eq:restriction} holds only as $\mathfrak
h$-modules. Practically it is enough for our purpose as integrable
highest weight representations are determined by their characters.
Thus we get

\begin{Theorem}\label{thm:branch'}
\textup{(1)} The strata of $\M_{\zeta^\bullet}(\bw)$ and $\M_0(\bw)$
are given by \eqref{eq:bulletstratum}, \eqref{eq:0stratum} respectively.
The criterion of the nonemptiness of $\Mreg_{\zeta^\bullet}(\bv,\bw)$,
$\Mreg_0(\bv,\bw)$ is given in Proposition~{\rm \ref{prop:affinemodulistrata}}.

\textup{(2)} If we define $V^{\bv^0,\mu}_{\bv',\lambda}$ by
\begin{gather*}
  (\pi_{0,\zeta^\bullet})_*\left(IC(\Mreg_{\zeta^\bullet}(\bv^0,\bw))
  \boxtimes \C_{\overline{S^{|\mu|}_\mu
  \big[X_{\zeta^\circ}\setminus \bigcup_{i\in I_0^0} C_i \big]}}[\dim]
  \right)
\\
\qquad{} =
   \bigoplus_{\bv',\lambda}
   IC(\Mreg_0(\bv',\bw)) \boxtimes
   \C_{\overline{S^{|\lambda|}_\lambda(\C^2\setminus\{0\}/\Gamma))}}[\dim]
   \otimes V^{\bv^0,\mu}_{\bv',\lambda},
\end{gather*}
we have
\begin{equation*}
  \bigoplus_{\bv^0,\mu} V^{\bv^0,\mu}_{\bv',\lambda}
  \otimes V^{(\affLevi+\mathfrak h)}(\bw - \bv^0-|\mu|\delta)
  \cong V^\g(\bw - \bv' - |\lambda|\delta).
\end{equation*}
Therefore
\begin{equation}\label{eq:mu}
  \Hom_{\affLevi+\mathfrak h}
  \big(V^{(\affLevi+\mathfrak h)}(\bw - \bv^0),
  V^\g(\bw-\bv'-|\lambda|\delta)\big)
  = \bigoplus_\mu V^{\bv^0-|\mu|\delta,\mu}_{\bv',\lambda}.
\end{equation}
\end{Theorem}

\begin{Remarks}
  (1) By the same argument as in the proof of \thmref{thm:branch} we
  have $V^{\bv^0,\mu}_{\bv',\varnothing} = 0$ unless
  $\mu=\varnothing$. Therefore the right hand side of \eqref{eq:mu}
  contains only the single term $V^{\bv^0,\varnothing}_{\bv',\lambda}$ if
  $\lambda = \varnothing$.

  (2) We can determine the individual $V^{\bv^0,\mu}_{\bv',\lambda}$
  from the branching coef\/f\/icients. First we show that that it is
  enough to consider the case $\mu=\varnothing$.
For a general $\mu$, we use the diagram \eqref{eq:basechange} and
argue as in the proof of \thmref{thm:branch}. We need to compute
\begin{equation*}
   \kappa'_*\left((\pi_{0,\zeta^\bullet})_*
   IC(\Mreg_{\zeta^\bullet}(\bv^0,\bw))\boxtimes
   (\pi_{0,\zeta^\bullet})_*
   \C_{\overline{S^{|\mu|}_\mu \big[X_{\zeta^\circ}\setminus
       \bigcup_{i\in I_0^0} C_i \big]}}[\dim]\right),
\end{equation*}
where the f\/irst $\pi_{0,\zeta^\bullet}$ is
$\pi_{0,\zeta^\bullet}\colon \M_{\zeta^\bullet}(\bv^0,\bw)
\to \M_0(\bv^0,\bw)$ and the second one is
$\pi_{0,\zeta^\bullet}\colon S^{|\mu|} X_{\zeta^\bullet}\to
S^{|\mu|} \C^2/\Gamma$. The f\/irst pushforward is already known, as
it is the case $\mu=\varnothing$. The second pushforward is easily
computed as $\pi_{0,\zeta^\bullet}$ is semismall, both $S^{|\mu|}
X_{\zeta^\bullet}$, $S^{|\mu|} \C^2/\Gamma$ are rationally smooth, and
f\/ibers can be described. Finally we can compute
\begin{equation*}
   \kappa'_*\left(IC(\Mreg_{\zeta^\bullet}(\bv^0,\bw))\boxtimes
   \C_{\overline{S^{|\lambda|}_\lambda(\C^2\setminus\{0\}/\Gamma)}}
   \boxtimes
   \C_{\overline{S^{|\lambda'|}_{\lambda'}(\C^2\setminus\{0\}/\Gamma)}}\right)
\end{equation*}
by calculating the degree of the map
\(
  \kappa'\colon
  S^{|\lambda|}_\lambda(\C^2\setminus\{0\}/\Gamma) \times
  S^{|\lambda'|}_{\lambda'}(\C^2\setminus\{0\}/\Gamma)
  \to S^{|\lambda|+|\lambda'|}_{\lambda\cup\lambda'}
  (\C^2\setminus\{0\}/\Gamma).
\)

When we consider a string `$\bv^0+\Z\delta$', there is a minimal
element $\bv^0_{\min}$ such that
\(
   \Mreg_{\zeta^\bullet}(\bv^0_{\min}-m\delta,\bw) = \varnothing
\)
for $m > 0$. Then the right hand side only has a single summand
$V^{\bv^0_{\min},\varnothing}_{\bv',\lambda}$ if $\bv^0 =
\bv^0_{\min}$. Thus it is determined by the branching coef\/f\/icient.
For general $\bv^0 = \bv^0_{\min} + m\delta$, we can determine
\(
   V_{\bv',\lambda}^{\bv^0,\varnothing}
\)
by induction on $m$, as summands appearing in the right hand side of~\eqref{eq:mu} are smaller.

(3)
  Suppose that $\g$ is of af\/f\/ine type $A_{r-1}^{(1)} = \algsl(r)_\aff$
  and $I_0^0 = I_0 \setminus \{r_1\}$ for some $1 \le r_1 < r$. Then we
  have $\g_{I_0^0} \!\cong\! \algsl(r_1)\oplus \algsl(r_2)$ where $r_2 \!= \!r
 \! - \!r_1$. By \eqref{eq:branchtensor} the factor
\(
  \Hom_{(\affLevi{+}\mathfrak h)}
  (V^{(\affLevi{+}\mathfrak h)}$ $(\bw - \bv^0),
  V^\g(\bw-\bv'))
\)
gives the tensor product multiplicity of
$V^{\algsl(l)_\aff}(\overline{\lambda})$ in
$V^{\algsl(l)_\aff}(\overline{\lambda_1})\otimes
V^{\algsl(l)_\aff}(\overline{\lambda_2})$ for appropriate
$\overline{\lambda}$, $\overline{\lambda_1}$, $\overline{\lambda_2}$.
This conf\/irms the conjecture \cite{braverman-prep} that the
convolution diagram for the double af\/f\/ine Grassmannian realizing the
tensor product is def\/ined by the Uhlenbeck partial compactif\/ication of
the framed moduli space of instantons on the orbifold
$X_{\zeta^\bullet}$.
However it remains to be clarif\/ied how we should interpret summands
$\lambda$, $\mu\neq\varnothing$ to make a `categorical' statement as in
the usual geometric Satake correspondence.

(4) In view of the level-rank duality it is natural to expect that
$H_*(\pi_{0,\zeta}^{-1}(0))$ has a structure of the Heisenberg algebra
module commuting with the $\g$-action. Since $\M_\zeta(\bw)$ is the
framed moduli space of torsion free sheaves on $X_{\zeta^\circ}$, such
a structure was constructed by the author \cite[Chapter~8]{Lecture}
(see \cite{Baranovsky} for a higher rank generalization). When $\bw$
is of level $1$, one can check that the Heisenberg algebra action
commutes with the $\g$-action mentioned in \remref{rem:Frenkel-Kac},
but it is still open in higher level.
\end{Remarks}

\section[MV cycles for the double affine Grassmannian of type
  $A$]{MV cycles for the double af\/f\/ine Grassmannian of type
  $\boldsymbol{A}$}\label{sec:MV}

The Mirkovi\'c--Vilonen cycles (MV cycles in short) are certain
subvarieties in the af\/f\/ine Grassmannian and are natural geometric
basis elements of a weight space of an irreducible f\/inite dimensional
representation of a f\/inite dimensional simple Lie algebra \cite{MV}.
Their conjectural double af\/f\/ine Grassmannian analogs are proposed
by Braverman and Finkelberg \cite{braverman-prep}.

Via the level-rank duality (see \eqref{eq:duality'}) their conjectural
basis for $V^{\algsl(l)_\aff}(\overline{\lambda})_{\overline{\mu}}$
also should give a basis of the tensor product multiplicity space for
$\gl(r)_\aff$.
Recall that the author introduced a~Lagrangian subvariety
$\widetilde{\mathfrak Z}(\bw)$ in $\bigsqcup_{\bv} \M_\zeta(\bv,\bw)$
such that the set of its irreducible components have a crystal
structure isomorphic to the tensor product for $\algsl(r)_\aff$
\cite{Na-tensor} (see also \cite{Malkin}).
A simple modif\/ication gives the tensor product for $\gl(r)_\aff$. Then
we show that those irreducible components of (the modif\/ied version of)
$\widetilde{\mathfrak Z}(\bw)$ intersecting with the open subvariety
$\bigcup \pi^{-1}(\Mreg_0(\bv,\bw))$ are exactly highest weight
vectors. In particular, the number of such irreducible components is
equal to the weight multiplicity. Note also that those irreducible
components are identif\/ied with those of a Lagrangian subvariety in
$\Mreg_0(\bv,\bw)$. The def\/inition makes sense for the Uhlenbeck
compactif\/ication of the framed moduli of $G$-bundles for {\it any\/}
$G$.

\begin{Remark}
  Together with the theory of the crystal base, we have the actual
  highest weight vectors in the tensor product representation (instead
  of the tensor product crystal) parametrized naturally by those
  irreducible components. But the fundamental classes of those
  irreducible components are not necessarily highest weight vectors,
  when we realize the tensor product representation as the homology
  group of $\widetilde{\mathfrak Z}(\bw)$ \cite{Na-tensor}. Note also
  that $\widetilde{\mathfrak Z}(\bw)$ makes sense only for $G = SL(l)$
  (or $GL(l)$). We do not know a satisfactory natural way to remedy
  this f\/law.
\end{Remark}

Let $(I,E)$ be the graph of the af\/f\/ine type $A_{r-1}$. We number the
vertices in the cyclic order as usual starting from $0$ to $r-1$. We
choose the cyclic orientation $\Omega$, i.e., $0\to 1$, $1\to 2$,
\dots, $r-2\to r-1$, $r-1 \to 0$. The corresponding Lie algebra $\g$
is $\algsl(r)_\aff$.

We choose the stability parameter $\zeta$ as in
Example~\ref{ex:stndcham}(1). For $I$-graded vector spaces $V$, $W$
with $\bv = \dim V$, $\bw = \dim W$ we
consider quiver varieties $\M_\zeta(\bv,\bw)$, $\M_0(\bv,\bw)$. We further
choose a decomposition
$W = W^1\oplus W^2\oplus \cdots \oplus W^l$ into $1$-dimensional
subspaces. (At the end the def\/inition of MV cycles depends only on the
f\/lag $0\subset W^1\subset W^1\oplus W^2 \subset \dots \subset W$.)
In particular, $W^p$ is concentrated at a degree $\mu_p\in I$, i.e.,
$\dim W^p = \be_{\mu_p}$.
We def\/ine a $1$-parameter subgroup $\rho_0\colon\C^*\to G_W$ by
\begin{equation*}
  \rho_0(t) = t^{m_1} \id_{W^1} \oplus t^{m_2} \id_{W^2} \oplus\cdots
  \oplus t^{m_l}\id_{W^l}
\end{equation*}
with $m_1 \ll m_2 \ll \cdots \ll m_l$. We def\/ine a $\C^*$-action on
$\bM(V,W)$ and the induced actions on $\M_\zeta(\bv,\bw)$, $\M_0(\bv,\bw)$ by
\begin{equation*}
   B_h \mapsto
   \begin{cases}
     B_h & \text{if $h\in\Omega$},
     \\
     t B_h & \text{if $h\in\overline\Omega$},
   \end{cases}
\qquad
   a \mapsto a \rho_0(t)^{-1},
\qquad
   b \mapsto t\rho_0(t)b.
\end{equation*}
Let
\begin{equation*}
  \mathfrak Z_\zeta(\bv,\bw) \defeq
  \left\{ [B,a,b]\in\M_\zeta(\bv,\bw) \big|\,
      \text{$\displaystyle\lim_{t\to \infty} t\cdot[B,a,b]$ exists}
       \right\}.
\end{equation*}
By \cite[\S~8]{Na-tensor} $\mathfrak Z_\zeta(\bv,\bw)$ is a Lagrangian
subvariety in $\M_\zeta(\bv,\bw)$.
Let $\Irr\mathfrak Z_\zeta(\bv,\bw)$ denote the set
of irreducible components of $\mathfrak Z_\zeta(\bv,\bw)$.
By \cite[\S~8]{Na-tensor}, its disjoint union
\(
   \bigsqcup_\bv \Irr\mathfrak Z_\zeta(\bv,\bw)
\)
has a structure of $\algsl(r)_\aff$-crystal. Moreover it is
isomorphic to the tensor product of the crystal
\(
   \bigsqcup_{\bv^p} \Irr\mathfrak Z_\zeta(\bv^p,\Lambda_{\mu_p})
\)
by the same proof as \cite[4.6]{Na-tensor}.
Here $\mathfrak Z_\zeta(\bv^p,\Lambda_{\mu_p})$ is def\/ined in the same
way as above applied to $\bw = \Lambda_{\mu_p}$. ($\rho(t)$ does not
matter in this case.)

Thanks to the following lemma, we have isomorphisms:
\begin{equation}\label{eq:tensor_Fock}
\begin{split}
  \bigsqcup_{\bv} \Irr\mathfrak Z_\zeta(\bv,\bw)
  \cong \bigotimes_{p=1}^l \left(
    \bigoplus_{\lambda^p}
  \mathscr B^{{\algsl}(r)_\aff}(\Lambda_{\mu_p})
  \otimes T_{-|\lambda^p|\delta}
  \right),
\end{split}
\end{equation}
where $\lambda^p$ runs over partitions. Here $\mathscr
B^{{\algsl}(r)_\aff}(\Lambda_{\mu_p})$ denotes the crystal of the
integrable highest weight representation
$V^{{\algsl}(r)_\aff}(\Lambda_{\mu_p})$, and $T_{-|\lambda^p|\delta}$
is the crystal consisting of a single element with weight
$-|\lambda^p|\delta$.

\begin{Lemma}\label{lem:Fock}
\(
   \bigsqcup_{\bv^p} \Irr\mathfrak Z_\zeta(\bv^p,\Lambda_{\mu_p})
\)
is isomorphic to
\(
  \bigoplus_{\lambda^p}
  \mathscr B^{{\algsl}(r)_\aff}(\Lambda_{\mu_p})
  \otimes T_{-|\lambda^p|\delta}.
\)
\end{Lemma}

We now introduce the analog of MV cycles:
\begin{equation*}
  \mathfrak Z_0^{\operatorname{s}}(\bv,\bw) \defeq
  \left\{ [B,a,b]\in\Mreg_0(\bv,\bw) \big|\,
      \text{$\displaystyle\lim_{t\to \infty} t\cdot[B,a,b]$ exists}
      \right\}.
\end{equation*}
Since $\pi = \pi_{0,\zeta}$ is projective, we have
$\mathfrak Z_0^{\operatorname{s}}(\bv,\bw)
= \pi(\mathfrak Z_\zeta(\bv,\bw))\cap\Mreg_0(\bv,\bw)$. Since $\pi$ is an
isomorphism on $\pi^{-1}(\Mreg_0(\bv,\bw))$ (see
\subsecref{subsec:partial}), this can be further identif\/ied with
$\mathfrak Z_\zeta(\bv,\bw)\cap\pi^{-1}(\Mreg_0(\bv,\bw))$. Since
$\pi^{-1}(\Mreg_0(\bv,\bw))$ is an open subset,
$\mathfrak Z_0^{\operatorname{s}}(\bv,\bw)$ is of pure dimension with
$\dim = \dim\M_\zeta(\bv,\bw)/2$. Its irreducible components are naturally
identif\/ied with irreducible components of
$\mathfrak Z_\zeta(\bv,\bw)$ intersecting with
$\pi^{-1}(\Mreg_0(\bv,\bw))$.

\begin{Theorem}\label{thm:MV}
  \textup{(1)} $\mathfrak Z_0^{\operatorname{s}}(\bv,\bw) \cong
  \mathfrak Z_\zeta(\bv,\bw)\cap\pi^{-1}(\Mreg_0(\bv,\bw))$ is of pure
  dimension with $\dim = \dim\M_\zeta(\bv,\bw)/2$.

  \textup{(2)} Let $Y_0$ be an an irreducible component of $\mathfrak
  Z_0^{\operatorname{s}}(\bv,\bw)$ and $\mathscr B(Y_0)$ be the
  connected component containing the closure of
  $\pi_{\zeta,0}^{-1}(Y_0)$. Then we have
  \begin{equation*}
   \bigsqcup_{\bv} \Irr\mathfrak Z_\zeta(\bv,\bw)
   \cong
   \bigsqcup_{\substack{Y_0\in\bigsqcup_{\bv} \Irr\mathfrak
       Z_0^{\operatorname{s}}(\bv,\bw) \\
       \text{\rm $\lambda$: partition}}}
     \mathscr B(Y_0)\otimes T_{-|\lambda|\delta}.
  \end{equation*}
  Furthermore $\mathscr B(Y_0)$ is isomorphic to the crystal of the
  integrable highest weight module
  $V^{\algsl(r)_\aff}\linebreak[2](\operatorname{wt} Y_0)$ so that $Y_0$ is the
  highest weight vector.
\end{Theorem}

Since the crystal of $V^{\gl(r)_\aff}(\operatorname{wt} Y_0)$ is
isomorphic to $\sqcup_{\lambda} \mathscr B(Y_0)\otimes
T_{-|\lambda|\delta}$, we conclude that $\Irr
\mathfrak Z_0^{\operatorname{s}}(\bv,\bw)$ parametrizes a basis of
\(
  \Hom_{\gl(r)_\aff} (V^{\gl(r)_\aff}(\bw - \bv),
  \bigotimes_{p=1}^l
  V^{{\gl}(r)_\aff}(\Lambda_{\mu_p}))
\)
in \eqref{eq:tensor_Fock}.

\begin{proof}[Proof of \lemref{lem:Fock}]
We consider $\M_0(\bv,\Lambda_{\mu_p})$. Thanks to
\propref{prop:affinestrata}(1), the stratif\/ica\-tion~\eqref{eq:Uhbullet}
is very simple in this case:
\begin{equation*}
  \M_0(\bv,\Lambda_{\mu_p}) \cong
  \bigsqcup_{\lambda,m_i} S^{|\lambda|}_\lambda (\C^2\setminus\{0\}/\Gamma)
  \times \left\{\bigoplus_i S_i^{\oplus m_i}\right\},
\end{equation*}
where the summation runs over partitions $\lambda$ and nonnegative
integers $m_i$ with $\bv = |\lambda|\delta + \sum m_i \be_i$. In this
description, the $\C^*$-action is the induced action from the action
$t\cdot(x,y)\mapsto (x,ty)$ on $\C^2$. Therefore the points
where the limit exists when $t\to\infty$ are
\begin{equation*}
  \bigsqcup_{\lambda,m_i} S^{|\lambda|}_\lambda
  ((\text{$x$-axis}) \setminus\{0\}/\Gamma)
  \times \left\{\bigoplus_i S_i^{\oplus m_i}\right\}.
\end{equation*}
From this description and the def\/inition of the crystal structure, we
know that the highest weight vectors, i.e., those irreducible
components $Y$ with $\varepsilon_i(Y) = 0$ for all $i\in I$, are
closures of $\pi^{-1}(S^{|\lambda|}_\lambda ((\text{$x$-axis})
\setminus\{0\}/\Gamma)$ in
$\M_\zeta(|\lambda|\delta,\Lambda_{\mu_p})$.
Note that this is irreducible as the punctual Hilbert scheme is
irreducible.
Moreover the component of the crystal containing this highest weight
vector is isomorphic to the crystal of $\bigsqcup_{\bv^p} \Irr
\La_\zeta(\bv^p,\Lambda_{\mu_p})\otimes T_{-|\lambda|\delta}$.
Thus we have
\begin{equation*}
   \bigsqcup_{\bv^p} \Irr\mathfrak Z_\zeta(\bv^p,\Lambda_{\mu_p})
   \cong \bigoplus_{\lambda^p}
   \left(\bigsqcup_{\bv^p} \Irr \La_\zeta(\bv^p,\Lambda_{\mu_p})\right)
   \otimes T_{-|\lambda|\delta},
\end{equation*}
where $\La_\zeta(\bv^p,\Lambda_{\mu_p}){=} \pi^{-1}(0)$.
Since $\bigsqcup_{\bv^p} \Irr \La_\zeta(\bv^p,\Lambda_{\mu_p})$ is
isomorphic to the crystal of $V^{\algsl(r)_\aff}(\Lambda_{\mu_p})$
\cite{KS,Saito,Na-tensor}, we have the assertion.
\end{proof}

\begin{proof}[Proof of \thmref{thm:MV}]
(1) is already proved.

(2) Let us consider the stratif\/ication of
$\M_0(\bv,\bw)$ in \eqref{eq:Uhbullet}:
\begin{equation*}
  \M_0(\bv,\bw) = \bigsqcup
  \Mreg_0(\bv^0,\bw)
  \times S^{|\lambda|}_\lambda(\C^2\setminus\{0\}/\Gamma)
  \times \left\{
    \bigoplus_{i\in I} S_i^{\oplus m_i}
    \right\}.
\end{equation*}
The points where the limit exists when $t\to \infty$ are
\begin{equation*}
  \bigsqcup
  \mathfrak Z^{\operatorname{s}}_0(\bv^0,\bw)
  \times S^{|\lambda|}_\lambda((\text{$x$-axis})\setminus\{0\}/\Gamma)
  \times \left\{
    \bigoplus_{i\in I} S_i^{\oplus m_i}
    \right\}
\end{equation*}
as in the proof of \lemref{lem:Fock}. The highest weight vectors are
of the closures of
\begin{equation*}
   \pi^{-1}(
   Y_0 \times
   S^{|\lambda|}_\lambda((\text{$x$-axis})\setminus\{0\}/\Gamma))
\end{equation*}
in $\M_\zeta(\bv^0+|\lambda|\delta,\bw)$ where $Y_0$ is an irreducible
component of $\mathfrak Z^{\operatorname{s}}_0(\bv^0,\bw)$.
From the def\/inition of the crystal structure, the connected component
containing the above vector is the tensor product
\(
   \mathscr B(Y_0) \otimes T_{-|\lambda|\delta}.
\)
This shows the f\/irst statement.

Again from the def\/inition of the crystal structure, $\mathscr B(Y_0)$
is isomorphic to the crystal $\bigsqcup_\bv \linebreak[3]\Irr\La_\zeta(\bv,\bw -
\bv^0)$, which is known to be isomorphic to the crystal of
$V^{\algsl(r)_\aff}(\bw - \bv^0)$ \cite{KS,Saito,Na-tensor}. This shows
  the second statement.
\end{proof}

\appendix
\section{Level-rank duality}\label{sec:level-rank}

Our formulation follows \cite{Hasegawa,NT}.

We denote the central extension of the loop Lie algebra by
$\widehat{\algsl}(r)$ while the af\/f\/ine Lie algebra is denoted by
$\algsl(r)_\aff$. The latter contains the degree operator $d$. Our
notation is slightly dif\/ferent from one in \cite{braverman-2007} and
distinguishes weights of $\algsl(r)_\aff$ and $\gl(r)_\aff$.

\subsection{Weight multiplicities}

Let $X$, $Y$ be f\/inite dimensional vector spaces of dimensions $l$,
$r$ respectively. Let
\(
  \mathscr L(X\otimes Y) = X\otimes Y\otimes t^{1/2} \C[t,t^{-1}].
\)
Let $F \equiv \Wedge^{\infty/2}\mathscr L(X\otimes Y)$ be the
semi-inf\/inite wedge space, or the fermionic Fock space. If we take a
basis $\{ x_p \}$ ($p=1,\dots, l$) of $X$ and $\{ y_i \}$ ($i=1,\dots,
r$) of $Y$, $\mathscr L(X\otimes Y)$ has a~basis
\(
   \{ x_p\otimes y_i \otimes t^n \}
\)
($p=1,\dots, l$, $i=1,\dots, r$, $n\in \Z + 1/2$). We put a total
ordering on the basis elements by the lexicographic ordering, f\/irst
read $n$, then $p$, f\/inally $i$. Then $F$ has a~basis
\begin{equation}\label{eq:basis}
   v_1\wedge v_2 \wedge v_3 \wedge \cdots,
\end{equation}
where $v_k = x_{p_k}\otimes y_{i_k} \otimes t^{n_k}$ is a basis
element in $\mathscr L(X\otimes Y)$ and we require that $v_1 > v_2 >
v_3 > \cdots$ and $v_{k+1}$ is the next element of $v_k$ for $k\gg 0$.
We have the fermion operators $\psi^{ip}(n)$ and their conjugate
operators $\psi_{ip}(n)$ acting on $F$ satisfying the Clif\/ford
algebra relations
\begin{equation*}
   \left\{ \psi^{ip}(n), \psi^{jq}(m) \right\} = 0
   = \left\{ \psi_{ip}(n), \psi_{jq}(m) \right\},
\qquad
   \left\{ \psi^{ip}(n), \psi_{jq}(m) \right\}
   = \delta_{ij}\delta_{pq} \delta_{m+n,0}.
\end{equation*}
The vacuum vector $|0\rangle\in F$ is the basis vector
\eqref{eq:basis} where $v_k$ runs over all $x_p\otimes y_i\otimes t^n$ with
$n < 0$.

The basis vectors above are parametrized by Maya diagrams $M$ of
$(l\times r)$-components:
\begin{equation*}
   M = \left\{ m_{ip}(n) \left|
  \begin{aligned}[m]
  & m_{ip}(n) = \text{$\square$ or $\graysquare$}\
  (1\le p\le l, 1\le i\le r, n\in\Z+1/2)
\\
  & \text{$m_{ip}(n) = \graysquare$
    (resp.\ $\square$) for
    $n \ll 0$ (resp.\ $n \gg 0$)}
    \end{aligned}
   \right.\right\}.
\end{equation*}
It can be visualized as
\begin{equation*}
   \text{\scriptsize $l$ boxes }\Big\{
   \Yvcentermath1
   \cdots
   \underbrace{\overset{-3/2}{\young(\rf\rf\rf,\rf\rf\rf)}}_{\text{$r$ boxes}}
   \,
   \overset{-1/2}{\young(\hf\rf\hf,\rf\rf\rf)}
   \,
   \overset{1/2}{\young(\hf\rf\hf,\hf\hf\rf)}
   \,
   \overset{3/2}{\young(\hf\hf\hf,\hf\hf\hf)}
   \cdots.
\end{equation*}
The corresponding basis element \eqref{eq:basis} is determined so that
$v_1$ corresponds to the f\/irst $\graysquare$ reading from the right,
$v_2$ corresponds to the second, etc.

We def\/ine the {\it degree operator} $d$ acting on $F$ by{\samepage
\begin{equation*}
   d(M) = \left(
     - \sum_{\substack{n > 0\\ m_{i,p}(n) = \graysquare}} n
     + \sum_{\substack{n < 0\\ m_{i,p}(n) = \square}} n
   \right) M.
\end{equation*}
We have $[d, \psi^{ip}(n)] = n \psi^{ip}(n)$, $[d, \psi_{ip}(n)] = n
\psi_{ip}(n)$.}

We have commuting actions of $\widehat{\algsl}(X)_r$,
$\widehat{\algsl}(Y)_l$ (the scripts indicate the levels) and the
Heisenberg algebra $\widehat{a}$ given by {\allowdisplaybreaks
\begin{alignat*}{2}
   & \widehat{\algsl}(X)_r & \qquad
   & J^p_q(n) \ (p\neq q), \quad
     J^p_p(n) - J_{p+1}^{p+1}(n) \ (p=1,\dots, l-1),
\\
   &&&
   J_p^q(n) \defeq \sum_{i=1}^l \sum_m :\!
   \psi^{ip}(n-m)\psi_{iq}(m)\!:,
\\
   & \widehat{\algsl}(Y)_l & \qquad
   & J^i_j(n) \ (i\neq j), \quad
     J^i_i(n) - J_{i+1}^{i+1}(n) \ (i=1,\dots, r-1),
\\
   &&&
   J_i^j(n) \defeq \sum_{p=1}^r \sum_m :\! \psi^{ip}(n-m)\psi_{jp}(m)\!:,
\\
  & \widehat{a} & \qquad
  &
     J(n) \defeq \sum_{i=1}^r \sum_{p=1}^l \sum_m
     :\! \psi^{ip}(n-m)\psi_{ip}(m)\!:\,
     = \sum_{i=1}^r J_i^i(n)
     = \sum_{p=1}^l J_p^p(n),
\end{alignat*}
where} $:\!\ \ \!:$ denotes the normal ordering, def\/ined by
\(
   :\! x(m) y(n) \!: \defeq x(m) y(n)
\)
(resp.\ $y(n) x(m)$, $1/2(x(m)y(n)+y(n)x(m))$) if
$n > m$ (resp.\ $n < m$, $n=m$).
The degree operator $d$ gives the degree operators for
$\widehat{\algsl}(X)$, $\widehat{\algsl}(Y)$ and $\widehat{a}$.

The branching formula of $F$ is as follows \cite{Hasegawa}:
\begin{equation}\label{eq:branch}
   F \cong
   \bigoplus_{\substack{\lambda\in\mathscr Y^r_l}}
   V^{\widehat{\algsl}(X)}(\overline\lambda)
   \otimes V^{\widehat{\algsl}(Y)}(\overline{\lsp{t}\lambda})
   \otimes H^{\widehat{a}}_{|\lambda|}.
\end{equation}
We need to explain the notation: The set $\mathscr Y^r_l$ of
generalized Young diagram consists of sequences $\lambda =
(\lambda_1,\dots,\lambda_l)$ ($\lambda_p\in\Z$, $\lambda_1\ge\dots\ge
\lambda_l$) with the level $r$ constraint $\lambda_1 - \lambda_l \le
r$. The {\it size\/} of $\lambda$ is $|\lambda| = \sum_{p=1}^l
\lambda_i$.
We def\/ine the Maya diagram $M(\lambda)$ associated with $\lambda$ by
\begin{equation*}
  m_{ip}(n) =
  \begin{cases}
   \graysquare
    & \text{if $r(n-\frac12) + i \le \lambda_p$}, \\
    \square
    & \text{otherwise}.
  \end{cases}
\end{equation*}
A generalized Young diagram $\lambda$ def\/ines a dominant weight of
$\widehat{\algsl}(X)$ of level $r$ by
\begin{equation*}
   (r - \lambda_1 + \lambda_l) \Lambda_0
   + \sum_{p=1}^{l-1} (\lambda_p - \lambda_{p+1}) \Lambda_p,
\end{equation*}
where $\Lambda_p$ is the $p^{\mathrm{th}}$ fundamental weight. We
denote it by $\overline{\lambda}$.
Then $V^{\widehat{\algsl}(X)}(\overline\lambda)$ is the corresponding
irreducible integrable highest weight module of $\widehat{\algsl}(X)$.
Conversely a dominant weight $\overline{\lambda}$ with level~$r$ gives a generalized Young diagram unique up to shift:
$(\lambda_1,\dots,\lambda_l)\mapsto (\lambda_1+k,\dots,\lambda_l+k)$.

The {\it transposition \/} $\lsp{t}{\ }\colon \mathscr Y^r_l\to
\mathscr Y^l_r$ is def\/ined by $m_{pi}(n) \defeq m_{ip}(n)$, i.e., the
transposition of each $(l\times r)$ rectangle in the Maya diagram.
Then $\lsp{t}\lambda$ def\/ines a dominant weight
$\overline{\lsp{t}\lambda}$ of $\widehat{\algsl}(Y)$ of level $l$. The
corresponding representation is denoted by
$V^{\widehat{\algsl}(Y)}(\overline{\lsp{t}\lambda})$.

Finally $H^{\widehat{a}}_k$ denote the irreducible representation of
$\widehat{a}$ with charge $k$. Here the charge is the eigenvalue of
$J(0)$ which counts the number of $\graysquare$ in the
region $n > 0$ minus the number of $\square$ in $n < 0$.

The basis vector corresponding to the Maya diagram $M(\lambda)$
($\lambda\in\mathscr Y^r_l$) gives the highest weight vector of
$V^{\widehat{\algsl}(X)}(\overline\lambda)\otimes
V^{\widehat{\algsl}(Y)}(\overline{\lsp{t}\lambda}) \otimes
H^{\widehat{a}}_{|\lambda|}$ in the decomposition \eqref{eq:branch}.

We have $\widehat{\algsl}(Y)\oplus\widehat{a} \cong \widehat{\gl}(Y)$.
We denote $V^{\widehat{\algsl}(Y)}(\overline{\lambda}) \otimes
H^{\widehat{a}}_{|\lambda|}$ by
$V^{\widehat{\gl}(Y)}(\lsp{t}\lambda)$. Note that the ambiguity of the
shift disappears if we consider $\lsp{t}\lambda$ as a weight for
$\widehat{\gl}(Y)$.
Then \eqref{eq:branch} can be rewritten as
\begin{equation}
  \label{eq:branch'}
   F \cong
   \bigoplus_{\substack{\lambda\in\mathscr Y^r_l}}
   V^{\widehat{\algsl}(X)}(\lambda)
   \otimes V^{\widehat{\gl}(Y)}(\overline{\lsp{t}\lambda}).
\end{equation}

Let $\lambda\in\mathscr{Y}^r_l$ and $\overline{\mu}$ be a dominant
weight of $\widehat{\algsl}(X)$ of level $r$. Then the weight space
$V^{\widehat{\algsl}(X)}(\overline{\lambda})_{\overline\mu}$ is
isomorphic to
\begin{gather*}
  V^{\widehat{\algsl}(X)}(\overline\lambda)_{\overline\mu}
    \cong
  \Hom_{\widehat{\gl}(Y)_l\oplus \widehat{\mathfrak h}({\widehat{\algsl}(X)})}\big(
  V^{\widehat{\gl}(Y)}(\lsp{t}\lambda)\otimes\C_{\overline\mu}, F\big)
\\
\phantom{V^{\widehat{\algsl}(X)}(\overline\lambda)_{\overline\mu}}{}  \cong
  \bigoplus
  \Hom_{\widehat{\gl}(Y)_l\oplus \widehat{\mathfrak h}(\widehat{\gl}(X))}\big(
  V^{\widehat{\gl}(Y)}(\lsp{t}\lambda)\otimes\C_\mu, F\big),
\end{gather*}
where $\widehat{\mathfrak h}(\widehat{\algsl}(X))$ (resp.\
$\widehat{\mathfrak h}(\widehat{\gl}(X))$) is the Cartan subalgebra of
$\widehat{\algsl}(X)$ (resp.\ $\widehat{\gl}(X)$), $\C_{\overline\mu}$
(resp.\ $\C_\mu$) is its representation with weight $\overline{\mu}$
(resp.\ $\mu$), and $\bigoplus$ runs over all $\mu\in\mathscr Y^r_l$
whose corresponding $\widehat{\algsl}(X)$-weight is the given
$\overline\mu$.
According to $\mu$, the space $X$ decomposes into direct sum of
$1$-dimensional eigenspaces, and hence we have
\begin{equation}\label{eq:fund}
  \Hom_{\widehat{\mathfrak h}(\widehat{\gl}(X))}(\C_\mu, F)
  \cong
  \bigotimes_{p=1}^l V^{\widehat{\gl}(Y)}(\Lambda_{\mu_p}),
\end{equation}
where $\Lambda_{\mu_p}$ is the $\mu_p^{\mathrm{th}}$-fundamental
weight of $\widehat{\algsl}(Y)$ with $\mu_p$ understood modulo $r$,
and $V^{\widehat{\gl}(Y)}(\Lambda_{\mu_p})$ is
$V^{\widehat{\algsl}(Y)}(\Lambda_{\mu_p})\otimes
H^{\widehat{a}}_{\mu_p}$. Thus we get
\begin{equation*}
  V^{\widehat{\algsl}(X)}(\overline\lambda)_{\overline\mu}
  \cong \bigoplus
  \Hom_{\widehat{\gl}(Y)_l}\Bigg(V^{\widehat{\gl}(Y)}(\lsp{t}\lambda),
  \bigotimes_{p=1}^l V^{\widehat{\gl}(Y)}(\Lambda_{\mu_p})\Bigg).
\end{equation*}
However the right hand side is $0$ unless $|\lambda|=|\mu|$, so
we have at most one $\mu$ contributing to $\bigoplus$.
On the other hand $(|\lambda|\bmod l)$ is well-def\/ined for
$\overline\lambda$. If $\overline\mu$ is a weight of
$V^{\widehat{\algsl}(X)}(\overline\lambda)$, then we must have
$|\lambda| \equiv |\mu|\pmod l$.
Hence we have exactly one $\mu$ with $|\lambda|=|\mu|$ and
\begin{equation}\label{eq:duality}
  V^{\widehat{\algsl}(X)}(\overline\lambda)_{\overline\mu}
  \cong
  \Hom_{\widehat{\gl}(Y)_l}\Bigg(V^{\widehat{\gl}(Y)}(\lsp{t}\lambda),
  \bigotimes_{p=1}^l V^{\widehat{\gl}(Y)}(\Lambda_{\mu_p})\Bigg).
\end{equation}

Now we incorporate the degree operator. Recall that an irreducible
highest weight module $V^{\widehat{\algsl}(X)}(\overline\lambda)$ of
$\widehat{\algsl}(X)$ has a lift to $\algsl(X)_\aff =
\widehat{\algsl}(X)\oplus\C d^{\widehat{\algsl}(X)}$, which is unique
if we f\/ix the value of $\langle
d^X,\overline\lambda\rangle$ \cite[\S~9.10]{Kac}.
The same is true for $\widehat{\gl}(Y)$.
We f\/ix $\langle d^X,\overline\lambda\rangle$, $\langle
d^Y,\lsp{t}\lambda\rangle$ (and hence the lifts of
$V^{\widehat{\algsl}(X)}(\overline\lambda)$,
$V^{\widehat{\gl}(Y)}(\lsp{t}\lambda)$) so that we have
\begin{equation*}
   d = d^X\otimes 1
   + 1\otimes d^Y
\end{equation*}
in \eqref{eq:branch'}. This is possible since $d$ gives the degree
operator for $\widehat{\algsl}(X)$ and $\widehat{\gl}(Y)$.
Therefore
\(
   \langle d^X,\overline\lambda\rangle
   + \langle d^Y,\lsp{t}\lambda\rangle
\)
is equal to the eigenvalue of $d$ for
the highest weight vector $v_{\overline{\lambda}}\otimes
v_{\lsp{t}\lambda}$, i.e., the vector corresponding to the Maya
diagram $M(\lambda)$. Let us denote the eigenvalue by $\langle d,
M(\lambda)\rangle$.
We also have $d$ on \eqref{eq:fund}, the restriction from that on $F$.
Thus $\bigotimes_{p=1}^l V^{\widehat{\gl}(Y)}(\Lambda_{\mu_p})$ has a
natural lift to a $\gl(Y)_\aff$-module.
From the construction, it is the tensor product of lifts of factors
$V^{\widehat{\gl}(Y)}(\Lambda_{\mu_p})$ and the eigenvalue of $d$ for
the highest weight vector $\bigotimes v_{\Lambda_{\mu_p}}$ is equal to
that for the vector corresponding to the Maya diagram $M(\mu)$, i.e.,
$\langle d, M(\mu)\rangle$. Thus $\bigotimes_{p=1}^l
V^{\widehat{\gl}(Y)}(\Lambda_{\mu_p})$ is $\bigotimes_{p=1}^l
V^{{\gl}(Y)_\aff}(\Lambda_{\mu_p})\otimes e^{\langle d,
  M(\mu)\rangle\delta^Y}$ as a $\gl(Y)_\aff$-module, if we def\/ine
$\langle d^Y, \Lambda_{\mu_p}\rangle = 0$ and the second factor
$e^{\langle d, M(\mu)\rangle\delta^Y}$ is the trivial
$\widehat{\gl}(Y)$-module with $d^Y$ acting by the multiplication by
$\langle d, M(\mu)\rangle$.
Then the right hand side of \eqref{eq:duality} has the induced
operator $d$, which is equal to $d^X$ in the left
hand side from the def\/inition.
Thus we can decompose both sides of \eqref{eq:duality} into
eigenspaces of~$d$:
\begin{gather}\label{eq:duality'}
  V^{{\algsl}(X)_\aff}(\overline\lambda)_{\overline\mu}
  \cong
  \Hom_{{\gl}(Y)_\aff}\Bigg(V^{{\gl}(Y)_\aff}(\lsp{t}\lambda
  + t \delta^Y),
  \bigotimes_{p=1}^l V^{{\gl}(Y)_\aff}(\Lambda_{\mu_p})\Bigg),
\end{gather}
where we have chosen a lift of $\overline\mu$ to a weight of
$\algsl(X)_\aff$ and denoted it by the same notation, and
$t = \langle d^X, \overline\mu\rangle- \langle d, M(\mu)\rangle$.
We have the relation
\begin{equation}\label{eq:weights_rel}
  \langle d^X,\overline{\lambda} - \overline{\mu}\rangle
  = - \langle d^Y, \lsp{t}\lambda + t\delta^Y \rangle
  + \langle d,M(\lambda)\rangle  - \langle d,M(\mu)\rangle.
\end{equation}

The above weights are related to those in the main body of the paper
(\remsref{rem:check}) by
\begin{equation*}
   \bw = \sum_i w_i \Lambda_i = \sum_{p=1}^l \Lambda_{\mu_p}, \qquad
   \bw - \bv = \sum_i w_i \Lambda_i - v_i \alpha_i
    = \overline{\lsp{t}\lambda} + t \delta^Y.
\end{equation*}

\subsection{Tensor product multiplicities}

We decompose as $Y = Y_1\oplus Y_2$ with $\dim Y_\alpha = r_\alpha$
($\alpha=1,2$). Then we have
$\mathscr L(X\otimes Y) \cong \mathscr L(X\otimes Y_1)\oplus \mathscr
L(X\otimes Y_2)$ and hence
\begin{equation*}
   \Wedge^{\infty/2} \mathscr L(X\otimes Y)
   \cong
   \Wedge^{\infty/2} \mathscr L(X\otimes Y_1) \otimes
   \Wedge^{\infty/2} \mathscr L(X\otimes Y_2).
\end{equation*}
We apply the decomposition \eqref{eq:branch'} to both hand sides:
\begin{gather*}
  \bigoplus_{\lambda\in\mathscr Y^r_l}
  V^{\widehat{\algsl}(X)}(\overline\lambda)\otimes
  V^{\widehat{\gl}(Y)}(\lsp{t}\lambda)
 \cong
  \bigoplus_{\substack{\lambda_1\in\mathscr Y^{r_1}_l \\
      \lambda_2\in\mathscr Y^{r_2}_l}}
  V^{\widehat{\algsl}(X)}(\overline\lambda_1)\otimes
  V^{\widehat{\algsl}(X)}(\overline\lambda_2) \otimes
  V^{\widehat{\gl}(Y_1)}(\lsp{t}\lambda_1)\otimes
  V^{\widehat{\gl}(Y_2)}(\lsp{t}\lambda_2).
\end{gather*}
Hence
\begin{gather*}
  \Hom_{\widehat{\algsl}(X)_r}\big(V^{\widehat{\algsl}(X)}(\overline\lambda),
  V^{\widehat{\algsl}(X)}(\overline\lambda_1)\otimes
  V^{\widehat{\algsl}(X)}(\overline\lambda_2)\big)
\\
\qquad{}  \cong
  \Hom_{\widehat{\algsl}(X)_r\oplus\widehat{\gl}(Y_1)_l\oplus\widehat{\gl}(Y_2)_l}
  \big(V^{\widehat{\algsl}(X)}(\overline\lambda)\otimes
  V^{\widehat{\gl}(Y_1)}(\lsp{t}\lambda_1)
  \otimes V^{\widehat{\gl}(Y_2)}(\lsp{t}\lambda_2),
  \Wedge^{\infty/2} \mathscr L(X\otimes Y)\big)
\\
\qquad{}  \cong
  \bigoplus
    \Hom_{\widehat{\gl}(Y_1)_l\oplus\widehat{\gl}(Y_2)_l}
  \big(V^{\widehat{\gl}(Y_1)}(\lsp{t}\lambda_1)
  \otimes V^{\widehat{\gl}(Y_2)}(\lsp{t}\lambda_2),
  V^{\widehat{\gl}(Y)}(\lsp{t}\lambda)\big),
\end{gather*}
where the summation runs over all $\lambda$ whose corresponding
$\algsl(Y)$-weight is the given $\overline\lambda$. Since this is $0$
unless $|\lambda| = |\lambda_1|+|\lambda_2|$, we only have a single
summand as in the previous subsection, and hence
\begin{gather*}
  \Hom_{\widehat{\algsl}(X)_r}\big(V^{\widehat{\algsl}(X)}(\overline\lambda),
  V^{\widehat{\algsl}(X)}(\overline\lambda_1)\otimes
  V^{\widehat{\algsl}(X)}(\overline\lambda_2)\big)
\\
 \qquad{} \cong
    \Hom_{\widehat{\gl}(Y_1)_l\oplus\widehat{\gl}(Y_2)_l}
  \big(V^{\widehat{\gl}(Y_1)}(\lsp{t}\lambda_1)
  \otimes V^{\widehat{\gl}(Y_2)}(\lsp{t}\lambda_2),
  V^{\widehat{\gl}(Y)}(\lsp{t}\lambda)\big).
\end{gather*}

We can incorporate the degree operator after f\/ixing the values of
$\langle d^X, \overline{\lambda}\rangle$,
$\langle d^X, \overline{\lambda}_1\rangle$,
$\langle d^X, \overline{\lambda}_2\rangle$
as in the previous subsection. We have
\begin{gather*}
  \Hom_{{\algsl}(X)_\aff}\big(V^{{\algsl}(X)_\aff}(\overline\lambda),
  V^{{\algsl}(X)_\aff}(\overline\lambda_1)\otimes
  V^{{\algsl}(X)_\aff}(\overline\lambda_2)\big)
\\
 \qquad{} \cong
    \Hom_{({\gl}(Y_1)\oplus{\gl}(Y_2))_\aff}
  \big(V^{{\gl}(Y_1)_\aff}(\lsp{t}\lambda_1)
  \otimes V^{{\gl}(Y_2)_\aff}(\lsp{t}\lambda_2),
  V^{{\gl}(Y)_\aff}(\lsp{t}\lambda)\big)
\end{gather*}
with the relation
\begin{gather*}
   \langle d^X, \overline{\lambda}\rangle
   - \langle d^X, \overline{\lambda}_1\rangle
   - \langle d^X, \overline{\lambda}_2\rangle
\\
 \qquad{}  =
   - \langle d^Y, \lsp{t}{\lambda}\rangle
   + \langle d^Y, \lsp{t}{\lambda}_1\rangle
   + \langle d^Y, \lsp{t}{\lambda}_2\rangle
   + \langle d, M(\lambda)\rangle
   - \langle d, M(\lambda_1)\rangle
   - \langle d, M(\lambda_2)\rangle.
\end{gather*}

We need to re-write this isomorphism in terms of $\algsl(Y_1)_\aff$,
$\algsl(Y_2)_\aff$.
Let ${a}_1\subset{\gl}(Y_1)$, ${a}_2\subset{\gl}(Y_2)$ be the
central subalgebras generated by
\(
   \sum_{i=1}^{r_1} J_i^i(0)
\)
and
\(
   \sum_{i=r_1+1}^{r} J_i^i(0)
\)
respectively. Let $\widehat{a}_1$, $\widehat{a}_2$ be the
corresponding Heisenberg subalgebras of
$\widehat{\gl}(Y_1)$, $\widehat{\gl}(Y_2)$.
We consider the subalgebras
{\allowdisplaybreaks
\begin{gather*}
\begin{split}
  & {a}^0   \defeq
  \left\langle J^0(0) \defeq - r_2 \sum_{i=1}^{r_1} J_i^i(0)
   + r_1 \sum_{j=r_1+1}^{r} J_j^j(0)\right\rangle
  \subset
  {a}_1\oplus {a}_2,
\\
 & \widehat{a}^0  \defeq
  \left\langle J^0(n) \defeq - r_2 \sum_{i=1}^{r_1} J_i^i(n)
   + r_1 \sum_{j=r_1+1}^{r} J_j^j(n) \right\rangle_{n\in\Z}
  \subset
  \widehat{a}_1\oplus \widehat{a}_2.
\end{split}
\end{gather*}
We} also put $a^0_\aff \defeq \widehat a^0\oplus \C d$ with the usual
commutator relation.
We have
\(
  {a}_1\oplus {a}_2 = {a}^0 \oplus {a},
\)
where $a$ is the subalgebra generated by $J(0)$.
Therefore
\begin{gather}\label{eq:homsp}
     \Hom_{({\gl}(Y_1)\oplus{\gl}(Y_2))_\aff}
  \big(V^{{\gl}(Y_1)_\aff}(\lsp{t}\lambda_1)
  \otimes V^{{\gl}(Y_2)_\aff}(\lsp{t}\lambda_2),
  V^{{\gl}(Y)_\aff}(\lsp{t}\lambda)\big)
\\
 \quad{} \cong
  \Hom_{\widehat{\algsl}(Y_1)_l\oplus\widehat{\algsl}(Y_2)_l\oplus
    \widehat{a}^0 \oplus \C d}
  \big(V^{{\algsl}(Y_1)_\aff}(\overline{\lsp{t}\lambda_1})
  \!\otimes\! V^{{\algsl}(Y_2)_\aff}(\overline{\lsp{t}\lambda_2})
  \!\otimes\! H^{{a}^0_\aff}_{-r_2|\lambda_1|+r_1|\lambda_2|},
  V^{{\algsl}(Y)_\aff}(\overline{\lsp{t}\lambda})\big),
\nonumber
\end{gather}
where the charge of ${a}^0_\aff$ is def\/ined as the eigenvalue of
\(
  J^0(0)
\)
and the eigenvalue of $d$ for the vacuum vector of
$H^{{a}^0_\aff}_{-r_2|\lambda_1|+r_1|\lambda_2|}$ is set to be $0$.
This space is isomorphic to the space of vectors in
$V^{\algsl(Y)_\aff}(\overline{\lsp{t}\lambda})$ killed by
$U^+({\algsl}(Y_1)_\aff)\times U^+({\algsl}(Y_2)_\aff)$,
\(
   h
   - \langle h,
   \overline{\lsp{t}\lambda_1}+\overline{\lsp{t}\lambda_2}
   \rangle\id
\)
($h\in \mathfrak h(\algsl(Y_1))\oplus\mathfrak h(\algsl(Y_2))$),
$J^0(n) - \delta_{0n} (-r_2|\lambda_1|+r_1|\lambda_2|)$
($n\ge 0$) and
\(
  d - \langle d^Y, \lsp{t}\lambda_1\rangle
    - \langle d^Y, \lsp{t}\lambda_2\rangle
\)
by assigning the image of the tensor product of the highest weight vectors
and the vacuum vector of a homomorphism.

Consider
\begin{equation*}
  \Hom_{\widehat{\algsl}(Y_1)_l\oplus\widehat{\algsl}(Y_2)_l\oplus
    a^0\oplus \C d}
  \big(V^{{\algsl}(Y_1)_\aff}(\overline{\lsp{t}\lambda_1})
  \otimes V^{{\algsl}(Y_2)_\aff}(\overline{\lsp{t}\lambda_2}),
  V^{{\algsl}(Y)_\aff}(\overline{\lsp{t}\lambda})\big).
\end{equation*}
This is isomorphic to the space of vectors as above, except that the
condition for $J^0(n)$ is required only for $n=0$.
This space is an $\widehat{a}^0$-module, and is spanned by vectors
obtained from vectors in \eqref{eq:homsp} applying various $J^0(n)$
for $n < 0$. But since $[d, J^0(n)] = n J^0(n)$, the condition for $d$
also implies that they must lie in \eqref{eq:homsp}.
We get
\begin{gather}
  \Hom_{{\algsl}(X)_\aff}\big(V^{{\algsl}(X)_\aff}(\overline\lambda),
  V^{{\algsl}(X)_\aff}(\overline\lambda_1)\otimes
  V^{{\algsl}(X)_\aff}(\overline\lambda_2)\big)\nonumber
\\
 \qquad{} \cong
  \Hom_{\widehat{\algsl}(Y_1)_l\oplus\widehat{\algsl}(Y_2)_l \oplus
    {a}^0\oplus \C d}
  \big(V^{{\algsl}(Y_1)_\aff}(\overline{\lsp{t}\lambda_1})
  \otimes V^{{\algsl}(Y_2)_\aff}(\overline{\lsp{t}\lambda_2}),
  V^{{\algsl}(Y)_\aff}(\overline{\lsp{t}\lambda})\big).\label{eq:branchtensor}
\end{gather}

\subsection*{Acknowledgments}
This work is supported by the Grant-in-aid
for Scientif\/ic Research (No.19340006), JSPS.
This work was started while the author was visiting the
Institute for Advanced Study with supports by the Ministry of
Education, Japan and the Friends of the Institute.
The author would like to thank to A.~Braverman and M.~Finkelberg for discussion on the
subject, and
to the referees for their careful readings and comments.

 \LastPageEnding
\end{document}